\documentclass[moor]{informs11}              

\usepackage[sort]{natbib}
\NatBibNumeric
\def\bibfont{\small}%
\bibpunct[, ]{[}{]}{,}{n}{}{,}%

\usepackage[colorlinks=true,breaklinks=true,bookmarks=true,urlcolor=blue,
citecolor=blue,linkcolor=blue,bookmarksopen=false,draft=false]{hyperref}

\def\EMAIL#1{\href{mailto:#1}{#1}}

\TheoremsNumberedBySection
\EquationsNumberedBySection 

\usepackage{smile}
\usepackage{mkolar_definitions}
\usepackage{mathrsfs}
\usepackage{dutchcal}
\usepackage{enumitem}
\usepackage{mathtools}
\usepackage{xcolor,soul}

\usepackage{lipsum}
\usepackage{amsfonts}
\usepackage{graphicx}
\usepackage{epstopdf}
\usepackage{algorithm}
\usepackage{algpseudocode}
\usepackage{amsopn}
\usepackage[OT1]{fontenc}
\usepackage{mathrsfs}
\usepackage{hyperref}
\usepackage{float}
\usepackage{amsfonts,amsmath,amssymb,url,xspace}
\usepackage{verbatim}
\usepackage{mathtools}
\usepackage[bottom]{footmisc}
\usepackage{enumitem}
\usepackage{arydshln,leftidx,mathtools}
\usepackage{lscape}
\usepackage{chngcntr}
\allowdisplaybreaks[1]

\newcommand{\mcal}{\mathcal{m}}
\def\0{{\boldsymbol 0}}

\newcommand{\blue}{\color{blue}}

\newcommand{\black}{\color{black}}

\usepackage{xargs}
\usepackage[colorinlistoftodos,prependcaption,textsize=tiny]{todonotes}
\newcommandx{\unsure}[2][1=]{\todo[linecolor=red,backgroundcolor=red!25,bordercolor=red,#1]{#2}}
\newcommandx{\change}[2][1=]{\todo[linecolor=blue,backgroundcolor=blue!25,bordercolor=blue,#1]{#2}}
\newcommandx{\info}[2][1=]{\todo[linecolor=OliveGreen,backgroundcolor=OliveGreen!25,bordercolor=OliveGreen,#1]{#2}}
\newcommandx{\improvement}[2][1=]{\todo[linecolor=Plum,backgroundcolor=Plum!25,bordercolor=Plum,#1]{#2}}

\newlength{\leftstackrelawd}
\newlength{\leftstackrelbwd}
\def\leftstackrel#1#2{\settowidth{\leftstackrelawd}%
{${{}^{#1}}$}\settowidth{\leftstackrelbwd}{$#2$}%
\addtolength{\leftstackrelawd}{-\leftstackrelbwd}%
\leavevmode\ifthenelse{\lengthtest{\leftstackrelawd>0pt}}%
{\kern-.5\leftstackrelawd}{}\mathrel{\mathop{#2}\limits^{#1}}}

\RequirePackage[capitalize,nameinlink]{cleveref}[0.19]
\crefname{section}{section}{sections}

\begin{document}

\TITLE{Convergence Analysis of Accelerated Stochastic Gradient Descent under the Growth Condition}

\ARTICLEAUTHORS{
\AUTHOR{You-Lin Chen$^{\blue \textstyle *}$}
\AFF{Department of Statistics, The University of Chicago, \EMAIL{youlinchen@uchicago.edu}}
\AUTHOR{Sen Na\footnote{\black equal contribution}}
\AFF{International Computer Science Institute\\
Department of Statistics, The University of California, Berkeley, \EMAIL{senna@berkeley.edu}}
\AUTHOR{Mladen Kolar}
\AFF{Booth School of Business, The University of Chicago, \EMAIL{mladen.kolar@chicagobooth.edu}}
} 

\RUNAUTHOR{Chen, Na, and Kolar}
\RUNTITLE{Accelerated SGD under the Growth Condition}

\ABSTRACT{We study the convergence of accelerated stochastic gradient descent for strongly convex objectives under the \textit{growth condition}, which states that the variance of stochastic gradient is bounded by a multiplicative part that grows with the full gradient, and a constant additive part. Through the lens of the growth condition, we investigate four widely used accelerated methods: Nesterov's accelerated method (NAM), robust momentum method (RMM), accelerated dual averaging method (DAM+), and implicit DAM+ (iDAM+). While these methods are known to improve the convergence rate of SGD under the condition that the stochastic gradient has bounded variance, it is not well understood how their convergence rates are affected by the multiplicative noise. In this paper, we show that these methods all converge to a neighborhood of the optimum	with~accelerated convergence rates (compared to SGD) even under the growth condition. In particular, NAM, RMM, iDAM+ enjoy acceleration only with a mild multiplicative noise, while DAM+ enjoys acceleration even with a large multiplicative noise. Furthermore, we propose a generic tail-averaged scheme that allows the accelerated rates of DAM+ and iDAM+ to nearly attain the theoretical lower bound (up to a logarithmic factor in the variance term). We conduct numerical experiments to support our theoretical conclusions.
}


\maketitle

\section{Introduction}\label{sec:1}

In this paper, we consider an unconstrained optimization problem
\begin{equation}\label{eq:unconstrained-optimization}
\min_{x \in \RR^{d}} f(x),
\end{equation}
where $f$ is continuously differentiable, $\Lcal$-smooth, and $\mcal$-strongly convex with respect to the $\ell_2$~norm $\|\cdot\|$; that is, for any $x,y\in \RR^d$,
\begin{equation}\label{equ:class}
\begin{split}
f(y) \leq  f(x)+\nabla f(x)^\top(y-x)+\frac{\Lcal}{2}\|y-x\|^{2}, \ \ \
f(y) \geq  f(x)+\nabla f(x)^\top(y-x)+\frac{\mcal}{2}\|y-x\|^{2}.
\end{split}
\end{equation}
For solving Problem \eqref{eq:unconstrained-optimization}, the first-order methods utilize only the first-order information, that~is~the gradient $\nabla f(x)$. Due to their simplicity and light computational costs, the first-order methods are widely used in machine learning for solving large-scale optimization problems.

In deterministic setting, where the noiseless gradient $\nabla f(x)$ can be accessed, it is known that the vanilla gradient descent has a sub-optimal convergence rate for optimizing a strongly convex objective, while accelerated methods, such as Nesterov's accelerated method (NAM) \citep{Nesterov1983method} and Polyak's heavy ball (HB) \citep{Polyak1964Some}, enjoy the optimal convergence rate \citep{Nesterov2004Introductory}. However, in stochastic setting, the story is not so clear as the same type of analysis is inapplicable. The common belief is that stochastic accelerated methods mimic their deterministic counterparts, resulting in certain practical gains; strictly speaking, their theoretical underpinnings remain incomplete.

A stochastic method can only query a {\it stochastic oracle} (SO) to obtain a noisy gradient, i.e., \vskip-10pt
\begin{equation*}
g(x) = \nabla f(x) + \varepsilon(x),
\end{equation*}
where $\varepsilon(x)$ is a mean-zero noise that depends on the iterate $x\in\RR^d$. The assumption on the {\it stochastic oracle's variance}, $\EE\|\varepsilon(x)\|^2$,  is crucial for both the design and analysis of a practical stochastic method. Most of the literature is focused on a {\it SO with upper bounded variance} (OUBV), that is $\sup_x \EE \|\varepsilon(x)\|^2 \leq \sigma^2$ for some constant $\sigma^2$. The OUBV assumption is reasonable when the domain~of $f$ is compact or the noise $\varepsilon(\cdot)$ does not depend on the current iterate $x$ (that is, the noise is additive). However, this assumption can also be restrictive as it excludes many important problem instances, including the least-squares regression (LSR). More importantly, existing analyses under OUBV are inconsistent with some empirical observations. For example, \cite{Assran2020Convergence, Aybat2020Robust} proved under OUBV that both HB and NAM have a faster rate of convergence compared to stochastic gradient descent (SGD), while \cite{Kidambi2018Insufficiency, Liu2020Accelerating} constructed a least-squares problem for which HB and NAM \textit{cannot} outperform SGD even with the best choice of tuning parameters. The reason for this discrepancy, which~motivates our work, stems from the violation of OUBV in least-squares~\mbox{problems}.~For~such~problems, $\varepsilon(x)$ depends on $x$ so that the variance $\EE\|\varepsilon(x)\|^2$ can be large and grow with the distance~$\|x - x_\star\|^2$.~Here,~$x_\star = \arg \min f(x)$ denotes the unique minimizer of $f$.

In this paper, we aim to answer the following fundamental question:
\vskip2pt
{\em Do accelerated SGD methods enjoy accelerated (even optimal) convergence rates, as in deterministic case, without assuming OUBV?}
\vskip2pt
\noindent We provide an affirmative answer to this question. We prove that under a \textit{growth condition} on~the variance $\EE\|\varepsilon(x)\|^2$, which is a weaker and more realistic condition, (some) accelerated SGD methods achieve near-optimal accelerated convergence rates. Specifically, our analysis unifies the error~recursions of four accelerated methods into a single form in \eqref{equ:infor:result}, which requires carrying out novel~and sharper derivations upon the existing literature. By the unified recursion, we examine different~bias-variance trade-offs of the four methods and compare their ability of robustness to the multiplicative noise. The comparison results are summarized in \cref{sec:4.4}.

We first discuss LSR in detail in the next section to motivate our study. We will illustrate that OUBV is not a reasonable assumption for this type of problem, but a growth condition is satisfied.

\subsection{Motivation}\label{sec:1.1}

Consider the following LSR problem
\begin{equation*}
\min_{x\in\RR^d} f(x) =
\frac{1}{2} \EE_{a, b \sim \Dcal} \left[ a^\top x -b\right]^2,
\end{equation*}
where $(a, b)$ is a covariate-response pair sampled from a distribution $\Dcal$. To simplify our discussion, we assume a linear model $b = a^\top x_{\text{LSR}} + \varepsilon_{\text{LSR}}$, where $\varepsilon_{\text{LSR}}$ is a random variable independent from $a$ with $\EE[\varepsilon_{\text{LSR}}] = 0$ and $\EE[\varepsilon_{\text{LSR}}^2] = \sigma_{\text{LSR}}^2$, and $\EE[aa^\top]$ is invertible. Note that $x_{\text{LSR}}(=x_\star)$ is the unique minimizer of $f$ in this example. Given $n$ i.i.d. samples $\{a_i, b_i\}_{i=1}^n$ from $\Dcal$, \cite[Theorem 1]{Mourtada2022Exact} showed that any estimator $\hat{x}$ of $x_{\text{LSR}}$ has at least the following minimax risk
\begin{equation}\label{s:error}
\EE[ f(\hat{x})]-f(x_{\text{LSR}}) = \Omega \left( \frac{d \sigma_{\text{LSR}}^2}{n} \right),
\end{equation}
where $\alpha_n=\Omega(\beta_n)$ denotes that $\alpha_n \geq c \beta_n$ for some constant $c$ and large enough $n$. The above~minimax risk can be achieved by the empirical risk minimizer. When applying the first-order stochastic approximation (SA) methods on $f$, where at each round only few samples are accessed to approximate $f$ and $\nabla f$, \cite[Corollary 2]{Jain2018Parallelizing} showed that, with some moment conditions on $\Dcal$, the tail-averaged SGD~meets the risk bound ($\kappa = \Lcal/\mcal$ is the condition number)\vskip-8pt
\begin{equation} \label{eq:SGD-bound}
\EE[f(x_n^{\text{SGD}})]-f(x_{\text{LSR}})= \Ocal\left( \exp\left(- \frac{n}{ \kappa}\right) + \frac{d \sigma_{\text{LSR}}^2}{n} \right).
\end{equation}\vskip-3pt
\noindent Here, $\alpha_n=\Ocal(\beta_n)$ denotes that $\alpha_n \leq c \beta_n$ for some constant $c$ and large enough $n$. The above bound consists of two terms: the second term is the statistical rate or the variance term, which is optimal as seen from \eqref{s:error}; the first term is the algorithmic rate or the bias term, which is \textit{not} optimal. In fact, by standard complexity theory of convex programs, e.g., \cite[Sections 5.3.1, 7.2.6]{Nemirovskij1983Problem}, \cite[Theory 2.1.13]{Nesterov2004Introductory}, \cite[(1.3)]{Ghadimi2013Optimal}, \cite[Corollary B.5]{Cohen2018Acceleration}, \cite[(5), (6)]{Aybat2019Universally}, and \cite{Chen2012Optimal}, it is well known that the risk~lower bound of first-order SA methods under OUBV is \vskip-5pt
\begin{equation}\label{eq:lower-bound}
\EE[f(x_n)]-f(x_{\text{LSR}}) =
\Omega\bigg(\underbrace{ \exp \left(- \frac{n}{\sqrt{ \kappa}} \right)}_{\text{bias term}}+ \underbrace{ \frac{d\sigma_{\text{LSR}}^2}{n}}_{\text{variance term}}\bigg).
\end{equation}\vskip-3pt
\noindent Compared to \eqref{eq:SGD-bound}, the dependence on the condition number $\kappa$ in \eqref{eq:lower-bound} improves from $1/\kappa$ to $1/\sqrt{\kappa}$. The lower bound can be attained under OUBV by an accelerated method designed in \cite{Ghadimi2013Optimal}. We note that \eqref{eq:lower-bound} is a valid lower bound even if a weaker assumption is imposed on the stochastic oracle than OUBV. This is simply because that the constructed problem instances for lower bound analysis are still applicable when weaker conditions are required. However, there exist~only limited works that addressed the problem whether an algorithm can attain the bound \eqref{eq:lower-bound} under weaker conditions~on the oracle. A careful study of this problem helps us to explain the empirical observation that~HB and NAM fail to accelerate SGD when applied to LSR in \cite{Kidambi2018Insufficiency, Liu2020Accelerating}, as OUBV does not hold for LSR in general. To see this clearly, the SO $g(x)$ in LSR is \vskip-8pt
\begin{equation*}
g(x) = a(a^\top x - b) = a(a^\top x - a^\top x_{\text{LSR}} - \varepsilon_{\text{LSR}}) = aa^\top(x- x_{\text{LSR}}) - a\varepsilon_{\text{LSR}},
\end{equation*}
where the second equality follows from the linear model setup. The noise is then \vskip-8pt
\begin{equation*}
\varepsilon(x) = g(x) - \nabla f(x) =(a a^\top-\Sigmab) (x - x_{\text{LSR}}) - a\varepsilon_{\text{LSR}},
\end{equation*}
where $\Sigmab = \EE aa^\top$. Therefore, OUBV assumption is not satisfied unless $x$ is in a compact domain. In particular, denoting $\Ib_d$ as a $d\times d$ identity matrix, we have that \vskip-15pt
\begin{align}\label{equ:LSR:GC}
\EE\|\varepsilon(x)\|^2 &= (x- x_{\text{LSR}})^\top\EE[aa^\top - \Sigmab]^2(x- x_{\text{LSR}}) + \EE[\varepsilon_{\text{LSR}}^2\cdot a^\top a] \nonumber\\[2pt]
&=  \nabla f(x)^\top\EE\sbr{\Sigmab^{-1}(aa^\top - \Sigmab)^2\Sigmab^{-1}}\nabla f(x) + \sigma_{\text{LSR}}^2\cdot\EE a^\top a \nonumber\\[2pt]
&=  \nabla f(x)^\top\rbr{\Sigmab^{-1}\EE(aa^\top)^2\Sigmab^{-1} - \Ib_d}\nabla f(x) + \sigma_{\text{LSR}}^2\cdot\EE a^\top a,
\end{align}
which motivates us to study a weaker condition on the variance, called the \textit{growth condition}.

\subsection{Growth condition}

As shown in \eqref{equ:LSR:GC}, a weaker condition on the noise $\varepsilon(x)$ that~is satisfied at least by LSR is the growth condition \citep{Solodov1998Incremental, Tseng1998Incremental, Schmidt2013Fast, Vaswani2019Fast, Nguyen2019New, Stich2019Unified}. The formal definition~is~as~follows.

\begin{definition}[Growth Condition] \label{asp:GC}
We say that the error $\varepsilon(x)$ satisfies the growth condition with constants $\delta\geq0$ and $\sigma^2$ if \vskip-8pt
\begin{equation}\label{eq:GC}
\EE \|\varepsilon(x)\|^2 \leq \delta \cdot \|\nabla f(x)\|^2 + \sigma^2.
\end{equation}
\end{definition}

The growth condition states that the variance of the stochastic gradient is dominated by a~multiplicative part $\delta\|\nabla f(x)\|^2$ and a constant additive part $\sigma^2$. The multiplicative part shrinks relative to the true gradient. We note that when $\delta = 0$, the growth condition reduces to OUBV.

Our work provides comprehensive analyses of accelerated methods under the growth condition. The {\it accelerated rate} refers to a convergence rate, whose bias term is improved from $1-c(\delta)/\kappa$ as \eqref{eq:SGD-bound} to $1-c(\delta)/\sqrt{\kappa}$ as \eqref{eq:lower-bound} with some constant $c(\delta)$ depending on $\delta$. In this sense, we analyze~four prominent accelerated stochastic methods: Nesterov's accelerated method (NAM), robust momentum method (RMM), accelerated dual averaging method (DAM+), and implicit DAM+ (iDAM+). We show that all these methods attain accelerated rates under the growth condition, although they can tolerate different levels of $\delta$. In particular, we show the following results.
\begin{enumerate}[label=(\alph*),topsep=2pt]
\setlength\itemsep{0.2em}
\item NAM, RMM, iDAM+ enjoy accelerated rates but only for \textit{mild} multiplicative noise (e.g.,~$\delta<1$).

\item DAM+ enjoys an accelerated rate for any $\delta\geq 0$. Thus, it uniformly improves upon SGD.

\item The rates of DAM+ and iDAM+, with a proper scheme for averaging iterates and diminishing algorithms' parameters, nearly match \eqref{eq:lower-bound} (up to a logarithmic factor in the variance term).

\item When $\delta = 0$, the rates that we establish reduce to the known results under OUBV assumption; and all the four methods accelerate SGD.
\end{enumerate}
\vskip2pt
\noindent Based on the above results, we now revisit the variance of the gradient's noise of LSR in \eqref{equ:LSR:GC}.~It~is straightforward to see that the multiplicative noise of LSR for some distribution $\Dcal$ is not ``mild". For example, if $a$ follows multivariate normal distribution with mean zero and covariance $\Sigmab$, plugging the equality $\EE(aa^\top)^2 = \tr(\Sigmab)\Sigmab + 2\Sigmab^2$ into \eqref{equ:LSR:GC} shows that $\delta>1$. This calculation may reveal the reason why NAM, RMM (and iDAM+) do not outperform SGD on LSR, as observed in \cite{Kidambi2018Insufficiency, Liu2020Accelerating}. Although such a statement is not a rigorous validation on the failure of NAM, RMM, and iDAM+, which would require a lower bound analysis, we indeed identify an important characteristic that~an accelerated method should enjoy to ensure acceleration on some problem instances.

In addition, although \cite{Liu2020Accelerating} proved that stochastic NAM fails to accelerate for LSR, and \cite{Cyrus2018Robust} argued that stochastic NAM is fragile for multiplicative noise by numerical studies, we rigorously prove, by a finer analysis, that it accelerates SGD for \textit{mild} multiplicative noise, which is a novel and rather surprising result in stochastic analysis. We should mention that our analyses of NAM and RMM~are only for constant algorithms' parameters, confined by the dissipativity framework \citep{Lessard2016Analysis} that we rely on for proofs\footnote{We are aware that \cite{Hu2017Dissipativity} analyzed the convergence properties of a stochastic method with varying parameters via the dissipativity framework, but that work only showed a sublinear convergence rate.}. We leave the analyses of NAM and RMM with diminishing parameters for future works. However, for DAM+ and iDAM+, we design a novel scheme to average iterates and diminish parameters, and show that their rates nearly attain the lower bound under the growth~condition.

\subsection{Contribution and related work}

Our paper contributes to the convergence analysis of accelerated stochastic first-order methods under the growth condition \eqref{eq:GC}.~The~\mbox{existing}~\mbox{convergence} results on accelerated stochastic methods are either under the assumption of OUBV \citep{Assran2020Convergence, Aybat2020Robust}, or under the growth condition but with restrictive setups, e.g., over-parameterized models \citep{Vaswani2019Fast} and quadratic objectives \citep{Jain2018Accelerating}. Here, the over-parameterized model means $\sigma^2=0$ in \eqref{eq:GC}, so that $\nabla f(x) = 0$ implies $g(x) = 0$. Our paper studies general strongly convex objectives without further model assumptions.

The growth condition was first introduced for finite-sum problems, $f(x) = n^{-1} \sum_{i=1}^n f_i(x)$.~Assuming for some $c>0$ that \vskip-15pt
\begin{equation}\label{equ:old:GC}
\max_{i=1,\ldots, n}\|\nabla f_i(x)\| \leq c\|\nabla f(x)\|,
\end{equation}\vskip-5pt
\noindent \cite{Solodov1998Incremental, Tseng1998Incremental} analyzed a deterministic incremental gradient method; \cite{Guerbuezbalaban2015globally} analyzed an incremental Newton method; and \cite{Schmidt2013Fast} established the linear convergence of SGD with constant stepsize for strongly~convex objectives. The convergence of SGD for both convex and nonconvex objectives under the~growth condition \eqref{eq:GC} was studied in the review paper \cite{Bottou2018Optimization}, which we revisit for strongly convex objectives in \cref{sec:2}. \cite{Cevher2018linear} showed that \eqref{eq:GC} with $\sigma^2=0$ (i.e. over-parameterized setup) is necessary~for~constant-stepsize SGD to have linear convergence rate, while \cite{Vaswani2019Fast} designed an accelerated method and showed under the same setup that their method attains the accelerated rate for strongly convex objectives. Furthermore, \cite{Jofre2018variance} introduced a notion of \textit{local} growth condition: \vskip-8pt
\begin{equation*}
\sqrt{\EE\|\varepsilon(x)\|^2} - \sqrt{\EE\|\varepsilon(x')\|^2} \leq c\|x - x'\|, \quad \forall x, x'.
\end{equation*} \vskip-5pt
\noindent Under some extra assumptions, the authors also analyzed an accelerated SGD method~called~FISTA, and proved that the method attains the optimal iteration complexity for convex objectives.

As the important complement of SGD, the dual averaging method (DAM) updates the iterates by reusing the past gradients information via averaging \citep{Xiao2010Dual}. We analyze its two accelerated variates---accelerated DAM (DAM+) \cite{Diakonikolas2018Accelerated, Diakonikolas2019Approximate} and implicit DAM+ (iDAM+) \cite{Cohen2018Acceleration}---under the growth condition \eqref{eq:GC} in \cref{sec:4}, which complements the existing literature that often studies DAMs under OUBV. In particular, \cite{Chen2012Optimal} studied a regularized DAM that shares the same spirit as DAM+ in the sense that both methods take a convex combination of SGD and DAM. \cite{Chen2012Optimal} showed that the regularized DAM achieves the optimal rate \eqref{eq:lower-bound} for strongly convex objectives with a multi-stage technique from \cite{Ghadimi2013Optimal} and OUBV condition. Furthermore, by performing implicit Euler discretization for continuous-time accelerated mirror descent dynamics \cite{Xu2018Continuous}, \cite{Diakonikolas2018Accelerated} derived accelerated extra-gradient descent (AXGD) method, which was then enhanced by $\mu$AGD+ method in \cite{Cohen2018Acceleration}. However, \cite{Cohen2018Acceleration} only showed sub-optimal convergence rate $\Ocal(\log k/k^{\log k} + \sigma^2\log k/k)$ under OUBV condition, where $k$ is the iteration index. Instead, by properly averaging the iterates, we show that both DAM+ and iDAM+ attain the rate $\Ocal(\exp(-k/\sqrt{\kappa}) + \sigma^2\log k/k)$ in the presence of the multiplicative noise. This rate matches \eqref{eq:lower-bound}~up to a $\log k$ factor in front of $\sigma^2$, and is stronger than the one of \cite{Cohen2018Acceleration} under weaker noise conditions.

To sum up, the existing literature either analyzed SGD under different forms of growth conditions or analyzed accelerated SGD and DAMs under OUBV or over-parameterized regime. Throughout our analysis, we do not require $\sigma^2=0$ or mini-batch gradients, but study accelerated methods under standard growth condition setup. We consider four popular accelerated methods, which allows us~to see how the multiplicative part in \eqref{eq:GC} affects the final convergence rate. Informally, we show that NAM, RMM, DAM+, and iDAM+ have a unified error recursion 
\begin{equation}\label{equ:infor:result}
\EE V_{k} \leq \left(1-  \frac{c_1(\delta) }{\sqrt{\kappa}} \right) \EE V_{k-1} + c_2(\delta) \sigma^2,
\end{equation}
where $V_{k}$ is a (algorithm-dependent) potential function at the $k$-th step, and~$c_1(\cdot), c_2(\cdot)$~are~(algorithm-dependent) constants depending on $\delta$ with $c_1(0)=1$. As a consequence, our results indicate that all the four methods converge to a neighborhood of the minimizer with the accelerated rates under~the growth condition. Our results also reveal how the rate varies with $\delta$: in general, if $\delta$ is small, the rate is fast; if $\delta$ is large or approaching to a certain limit, $c_1(\delta)$ approaches to $0$ and the rate is~slow. Based on \eqref{equ:infor:result}, we further develop a unified tail-averaged framework in \cref{sec:5} and show that DAM+ and iDAM+ (nearly) attain the optimal convergence rate.

\subsection{Structure of the paper}
In \cref{sec:2}, we introduce preliminaries including some implications of strong convexity and smoothness, some basic properties of Bregman divergence, dissipativity theory, and convergence analysis of SGD. We study NAM and RMM in \cref{sec:3}; and~DAM+ and iDAM+ in \cref{sec:4}. In \cref{sec:5}, we design an averaging scheme for DAM+ and iDAM+ and study their convergence. Discussions and conclusions are presented in \cref{sec:6}.

\section{Preliminaries}\label{sec:2}

Throughout the presentation, we use $\|\cdot\|$ to denote the $\ell_2$ norm; $\otimes$ to~denote the Kronecker product; bold letters to denote matrices; $\Ib_d$ and $\0_d$ to denote the $d\times d$ identity and zero matrices, respectively. $\Fcal(\mcal, \Lcal)$ denotes the class of functions that are continuously differentiable, $\mcal$-strongly convex and $\Lcal$-smooth (see \eqref{equ:class}), and $\kappa=\Lcal/\mcal$ denotes the condition number~of $f$. For two positive sequences $\{\alpha_n, \beta_n\}_{n=1}^\infty$, $\alpha_n = \Ocal(\beta_n)$ and $\alpha_n = \Omega(\beta_n)$ mean $\alpha_n \leq c \beta_n$ and $\alpha_n \geq c \beta_n$, for some positive constant $c$, respectively. All matrix inequalities hold in the semidefinite sense.

The following two lemmas summarize some basic implications of strong convexity and smoothness conditions. We omit the proof but refer to \cite{Nesterov2004Introductory, Bubeck2015Convex} and references therein.

\begin{lemma}\label{lemma:convex}
For $f \in \Fcal(\mcal,\infty)$, the following statements hold.
\begin{enumerate}[label=(\alph*),topsep=2pt]
\setlength\itemsep{0.2em}
\item $f(y) \leq f(x)+\nabla f(x)^\top(y-x)+\|\nabla f(y)-\nabla f(x)\|^{2}/(2 \mcal) ,\ \forall x, y$;
\item $(\nabla f(y) - \nabla f(x))^\top(y-x) \geq \mcal \|y-x\|^2 ,\ \forall x, y$;
\item $\|\nabla f(y) - \nabla f(x)\|\geq \mcal \|y-x\| ,\ \forall x, y$.
\end{enumerate}
\end{lemma}

\begin{lemma}\label{lemma:smooth}
For $f \in \Fcal(0,\Lcal)$, the following statements hold.
\begin{enumerate}[label=(\alph*),topsep=2pt]
\setlength\itemsep{0.2em}
\item $\|\nabla f(y)-\nabla f(x)\| \leq \Lcal\|y-x\|,\ \forall x, y$;
\item $f(y) \geq f(x)+\nabla f(x)^\top(y-x)+\|\nabla f(y)-\nabla f(x)\|^{2}/(2 \Lcal),\ \forall x, y$;
\item $(\nabla f(y) - \nabla f(x))^\top(y-x) \geq \|\nabla f(y)-\nabla f(x)\|^2/\Lcal ,\ \forall x, y$;
\item $(\nabla f(y) - \nabla f(x))^\top(y-x) \leq \Lcal \|y-x\|^2 ,\ \forall x, y$.
\end{enumerate}
\end{lemma}

Given a continuously differentiable, strongly convex function $f$, the Bregman divergence associated with $f$ between two points $x, y$ is defined as the difference between the value of $f$ at $x$ and~the first-order Taylor expansion of $f$ around $y$ evaluated at $x$. Specifically, we let \vskip-8pt
\begin{equation} \label{eq:Bregman-divergence}
\Delta_{f}(x, y)=f(x)-f(y)- \nabla f(y)^\top(x-y).
\end{equation}\vskip-3pt
\noindent We also define the convex conjugate $f^\star$ of $f$\footnote{The convex conjugate is defined for any function and takes supremum over $x$. But our paper only considers continuously differentiable, strongly convex $f$ on $\RR^d$, which implies the closeness so that the supreme is attained.} as \vskip-10pt
\begin{equation} \label{eq:convex-conjugate}
f^\star(z) = \max_x \left\{ z^\top x - f(x) \right\}.
\end{equation} \vskip-3pt
\noindent The next lemma presents some properties of Bregman divergence and convex conjugate. See \cite{Zhou2018fenchel, Bauschke1997Legendre} for the proofs.

\begin{lemma} \label{lemma:BD}
The following statements hold for $\Delta_f$ and $f^\star$.
\begin{enumerate}[label=(\alph*),topsep=2pt]
\setlength\itemsep{0.2em}
\item $\nabla f^{\star}(z)=\arg \max_{x}\{z^\top x - f(x)\}$, $\nabla f(x)=\arg \max_{z}\{x^\top z - f^\star(z)\}$, $\left(\nabla f^\star \right)(\nabla f(x))=x$, and $(\nabla f)(\nabla f^\star(z)) = z$;

\item $\Delta_{f}(x, y)=\Delta_{f^{\star}}(\nabla f(y), \nabla f(x))$;

\item If $f(x) \leq g(x), \forall x$ for some function $g$, then $f^\star(z) \geq g^\star(z), \forall z$.
\end{enumerate}
\end{lemma}

Next, we introduce an analytical tool called integral quadratic constraints (IQCs), which is used to analyze NAM and RMM in \cref{sec:3}. IQC has a close relationship to dissipativity theory~in~control area, and is a popular technique used to provide unified analysis for different stochastic methods \citep{Lessard2016Analysis, Hu2017Dissipativity, Hu2018Dissipativity, Hu2020Analysis}. Formally, consider the stochastic iteration scheme $x_{k+1} = \Ab x_k + \Bb\omega_k$, where $x_k\in \RR^d$, $\omega_k\in \RR^{d_\omega}$, $\Ab \in \RR^{d\times d}$, $\Bb\in \RR^{d\times d_{\omega}}$, and $\{\omega_i\}_{i=1}^k$ is a stochastic process. Given $\Xb_j = \Xb_j^\top \in  \RR^{(d+d_\omega) \times (d+d_\omega)}$ for $j=1,\dots, J$ (usually defined based on the iteration scheme or problem conditions) and $0 \leq \rho\leq 1$, if we can find a positive semidefinite matrix $\Pb\in \RR^{d \times d}$ and non-negative scalars $\{\lambda_j\}_{j=1}^J$ such that \vskip-10pt
\begin{equation} \label{eq:LMI}
\begin{bmatrix}
\Ab^\top \Pb \Ab - \rho^2 \Pb\;\;\; & \;\;\;\Ab^\top \Pb \Bb \\
\Bb^\top \Pb \Ab \;\;\; & \;\;\; \Bb^\top \Pb \Bb
\end{bmatrix}
- \sum_{j=1}^J \lambda_j \Xb_j \leq \zero,
\end{equation}\vskip-5pt
\noindent then we can define $S_j(x, \omega)=(x^\top, \omega^\top)\Xb_j(x^\top, \omega^\top)^\top$ and $V(x) = x^\top \Pb x$, and immediately get \vskip-10pt
\begin{equation} \label{eq:exp-dissipation-expectation}
V(x_{k+1}) \leq \rho^2 V(x_{k}) + \sum_{j=1}^J \lambda_j S_j(x_k, \omega_k),
\end{equation}\vskip-5pt
\noindent where $V(x)$ is called either a Lyapunov function or a potential function that measures the optimality gap of $x$. Moreover, if $\sum_{j=1}^J\lambda_j \EE S_j(x_k, \omega_k) \leq c_k$ for some $c_k$, then \eqref{eq:exp-dissipation-expectation} leads to the error recursion~of $V_k = V(x_k)$, i.e. $\EE{V_{k+1}} \leq \rho^2 \EE{V_{k}} + c_k$. See Theorem 3 in \cite{Hu2018Dissipativity} for more details.

\subsection{SGD under the growth condition} \label{subsec:3-SGD}

We review the convergence of SGD for strongly convex objective $f\in \Fcal(\mcal,\Lcal)$ under the growth condition \citep{Bottou2018Optimization}. We first introduce some additional notation. Given a stochastic iteration sequence $\{x_i\}_{i=0}^\infty$, we let $\Gcal_0 \subseteq \Gcal_1\subseteq \Gcal_2\subseteq\ldots$ be a filtration of $\sigma$-algebras where $\Gcal_k = \sigma(\{x_i\}_{i=0}^k)$ contains all the randomness up to the iteration $k$. At the $k$-th iteration with the iterate $x_k$, we let $\varepsilon_k(x_k)$ be a realization that is drawn from $\varepsilon(x_k)$. Furthermore, the gradient realization is $g_k(x_k) = \nabla f(x_k) + \varepsilon_k(x_k)$. \textit{Throughout the paper and for all algorithms, the (stochastic) gradient is always evaluated at $x_k$, even if algorithms involve other variables in the scheme}. To ease notation, we let $g_k = g_k(x_k)$, $\nabla f_k = \nabla f(x_k)$, and $\varepsilon_k = \varepsilon_k(x_k)$. We further have
\begin{equation}\label{pequ:1}
\begin{aligned}
\EE[\varepsilon_k \mid \Gcal_k] &= 0, \hskip1.9cm \EE[g_k\mid \Gcal_k] = \nabla f_k,\\[2pt]
\EE[\|g_k^2\|\mid \Gcal_k] &= \|\nabla f_k\|^2 +  \EE[\|\varepsilon_k\|^2\mid \Gcal_k] \stackrel{\eqref{eq:GC}}{\leq} (1+\delta)\|\nabla f_k\|^2 + \sigma^2.
\end{aligned}
\end{equation}

The (constant-stepsize) SGD scheme takes the form $x_{k+1} = x_k - \eta g_k$. By Lemma \ref{lemma:convex}(b) and~\ref{lemma:smooth}(c), we get $\nabla f_k^\top(x_k - x_\star) \geq \mcal \|x_k-x_\star\|^2$ and $\nabla f_k^\top(x_k-x_\star) \geq \|\nabla f_k\|^2/\Lcal$,  which implies \vskip-5pt
\begin{equation}\label{pequ:2}
-2\nabla f_k^\top(x_k - x_\star) \leq - \mcal \|x_k -x_\star\|^2 - \frac{1}{\Lcal} \|\nabla f_k\|^2.
\end{equation}
Furthermore, using the above inequalities, we have \vskip-20pt
\begin{align*}
\EE\left[ \|x_{k+1} - x_\star\|^2 \mid \Gcal_k \right]
&\stackrel{\eqref{pequ:1}}{=}  \|x_k-x_\star\|^2 - 2 \eta \nabla f_k^\top (x_k - x_\star) + \eta^2 \EE \left[ \| g_k  \|^2 \mid \Gcal_k \right] \\[2pt]
&\hskip-0.3cm \stackrel{\eqref{pequ:1},\eqref{pequ:2}}{\leq} (1 - \mcal \eta) \|x_k-x_\star\|^2 - \frac{\eta}{\Lcal} \|\nabla f_k \|^2 + \eta^2 \rbr{(1+\delta)\|\nabla f_k \|^2 + \sigma^2} \\[2pt]
&\hskip4pt = (1 - \mcal \eta) \|x_k-x_\star\|^2  - \left( \frac{\eta}{\Lcal}- \eta^2 - \delta\eta^2\right) \|\nabla f_k \|^2 + \eta^2 \sigma^2.
\end{align*} \vskip-5pt
\noindent Letting $V_k = \|x_k - x_\star\|^2$ and taking full expectation, we further have \vskip-8pt
\begin{equation}\label{pequ:3}
\EE \left[ V_{k+1} \right] - (1-\mcal\eta) \EE \left[V_{k}\right] \leq  - \rbr{ \frac{\eta}{\Lcal}- \eta^2 - \delta\eta^2} \EE \|\nabla f_k \|^2 + \eta^2 \sigma^2.
\end{equation}\vskip-2pt
\noindent From \eqref{pequ:3}, we see that SGD has a term $-(\eta/\Lcal - \eta^2) \EE \|\nabla f_k \|^2$ to eliminate the multiplicative noise term $\delta \eta^2\|\nabla f_k\|^2$, brought by the growth condition. In particular, since \vskip-8pt
\begin{equation*}
\frac{\eta}{\Lcal}- \eta^2 \geq \delta\eta^2  \Longleftrightarrow \eta \leq \frac{1}{(1+\delta)\Lcal},
\end{equation*}
we can set $\eta = 1/((1+\delta)\Lcal)$ to achieve the fastest rate for SGD, and arrive at the error recursion \vskip-6pt
\begin{equation} \label{eq:vSGD}
\EE \left[ V_{k+1} \right]
\leq (1-\mcal\eta) \EE \left[V_{k}\right] + \eta^2 \sigma^2
=  \left(1- \frac{1}{(1+\delta)\kappa}\right) \EE \left[V_{k}\right] + \frac{\sigma^2}{(1+\delta)^2 \Lcal^2}.
\end{equation}\vskip-3pt
\noindent Applying \eqref{eq:vSGD} iteratively, we obtain \vskip-10pt
\begin{equation}\label{pequ:SGD}
\EE V_k \leq \rbr{1 - \frac{1}{(1+\delta)\kappa}}^k\EE V_0 + \frac{\sigma^2}{(1+\delta)\Lcal\mcal}, \quad \forall k\geq 0.
\end{equation}
The above derivation with $\delta = 0$ recovers the sharpest SGD analysis under OUBV \citep[Theorem 3.1]{Gower2019SGD}\citep[Theorem 4.6]{Bottou2018Optimization}. Given the recursion \eqref{pequ:SGD}, we rise the following question:
\vskip 2pt
\emph{Is it possible to apply accelerated SGD methods (with constant parameters) to improve the dependence on the condition number in bias term from $1/\kappa$ to $1/\sqrt{\kappa}$, under the growth condition?}
\vskip 2pt

We provide an affirmative answer for NAM, RMM in \cref{sec:3}; and DAM+, iDAM+ in \cref{sec:4}. In particular, we show that, for a certain (but different) range of $\delta$, the four methods all achieve~the accelerated rates.

\section{Nesterov's Accelerated Method and Robust Momentum Method} \label{sec:3}

We study NAM and RMM and show convergence properties under the growth condition. In particular, we establish the error recursions for NAM and RMM in terms of the multiplicative noise level $\delta$.

\subsection{NAM under the growth condition}

The stochastic NAM \citep{Scoy2018Fastest} has a scheme \vskip-20pt
\begin{subequations}\label{eq:SNV}
\begin{align}
x_{k} &= y_{k} + \beta (y_{k}-y_{k-1}), \label{eq:SNV:a}\\
y_{k+1} &= y_{k} + \beta (y_{k}-y_{k-1}) - \eta g_k, \label{eq:SNV:b}
\end{align}
\end{subequations}\vskip-5pt
\noindent where we recall that $g_k$ is evaluated at $x_k$. Note that \eqref{eq:SNV} recovers HB if we replace \eqref{eq:SNV:a} by~$x_k = y_k$. The following theorem characterizes NAM under the growth condition.

\begin{theorem} \label{thm:nam}
Consider NAM in \eqref{eq:SNV:a}-\eqref{eq:SNV:b}. Suppose that  the error $\varepsilon(x)$ satisfies the growth condition in Definition \ref{asp:GC} with constants $\delta$, $\sigma^2$. Suppose $f \in \Fcal( \mcal,\Lcal)$, $0<\eta\leq 1/\mcal$, and set \begin{equation*}
\beta = \frac{1-\sqrt{\eta \mcal}}{1+\sqrt{\eta \mcal}}, \quad
\rho^2 = 1-\sqrt{\eta \mcal}, \quad
\Pb = \left(
\begin{bmatrix}
	\sqrt{\frac{1}{2 \eta}} & \sqrt{\frac{\mcal}{2}} -     \sqrt{\frac{1}{2 \eta}}
\end{bmatrix}^\top
\begin{bmatrix}
\sqrt{\frac{1}{2 \eta}} & \sqrt{\frac{\mcal}{2}} -     \sqrt{\frac{1}{2 \eta}}
\end{bmatrix}
\right) \otimes \Ib_d.
\end{equation*}
Further, we let $\xi_k = ((y_k - x_\star)^\top, (y_{k-1} - x_\star)^\top)^\top$ and define the potential function as
\begin{equation}\label{equ:NAM:V}
V_k = \xi_{k}^\top\Pb\xi_{k} + f(y_{k})-f(x_\star).
\end{equation}
Then, for all $k\geq 1$, we have
\begin{equation} \label{eq:NAM-recursive}
\EE V_{k+1}\leq \rho^2\EE V_k  - \frac{\eta}{2} \cbr{(1-\Lcal\eta) - (1+\Lcal\eta)\delta  }\EE \|\nabla f_k\|^2 - \tau(\eta) \EE \left[ \|y_{k}- y_{k-1}\|^2 \right] + \frac{\eta(1+\Lcal\eta)}{2}\sigma^2,
\end{equation}
where $\tau(\eta) = \cbr{\mcal(1-\sqrt{\eta \mcal})^3}/\cbr{2\eta \mcal + 2\sqrt{\eta \mcal}}\geq0$.
\end{theorem}

\begin{proof}
The proof follows the dissipativity framework introduced in \cref{sec:2}. We derive two IQCs and combine them to obtain the recursion of the potential function $V_k$ in \eqref{equ:NAM:V}. For the scheme~\eqref{eq:SNV}, let us define
\begin{equation*}
\omega_k =
\begin{bmatrix}
\nabla f_k \\
\varepsilon_k
\end{bmatrix},\quad\quad  \Ab =
\begin{bmatrix}
1+\beta & -\beta \\
1 & 0
\end{bmatrix}
\otimes \Ib_d, \quad\quad \Bb =
\begin{bmatrix}
-\eta & -\eta \\
0 & 0
\end{bmatrix}
\otimes \Ib_d.
\end{equation*}
Then $\xi_{k+1} \stackrel{\eqref{eq:SNV}}{=} \Ab\xi_{k} + \Bb\omega_k$. By the convexity of $f$, we have
\begin{equation} \label{eq:class-1}
f(y_k) - f(x_k) \stackrel{\eqref{equ:class}}{\geq} \nabla f_k^\top (y_k-x_k) + \frac{\mcal}{2}\|y_k - x_k\|^2.
\end{equation}
By the smoothness of $f$ and the fact that $y_{k+1}\stackrel{\eqref{eq:SNV:b}}{=} x_k - \eta g_k$, we have
\begin{align} \label{eq:class-2}
f(x_k) - f( y_{k+1})
&\stackrel{\mathclap{\eqref{equ:class}}}{\geq} \eta \nabla f_k ^\top g_k - \frac{\Lcal \eta^2}{2} \|g_k\|^2
= \eta \|\nabla f_k\|^2 + \eta \nabla f_k^\top \varepsilon_k - \frac{\Lcal \eta^2}{2}(\|\nabla f_k\|^2+2\nabla f_k^\top \varepsilon_k+\|\varepsilon_k\|^2 ) \nonumber\\[2pt]
&= \frac{\eta}{2} (2-\Lcal \eta)\|\nabla f_k\|^2 - \frac{\Lcal \eta^2}{2} \|\varepsilon_k\|^2 + \eta(1-\Lcal \eta) \nabla f_k^\top \varepsilon_k,
\end{align}
where the second equality uses $g_k = \nabla f_k + \epsilon_k$. Noting that
\begin{equation} \label{eq:NAM-2}
y_k-x_k \stackrel{\eqref{eq:SNV:a}}{=} -\beta (y_k-x_\star) + \beta (y_{k-1}-x_\star),
\end{equation}
we sum up \eqref{eq:class-1} and \eqref{eq:class-2} and obtain
\begin{align} \label{eq:NAM-IQC-1}
f(y_k) - f( y_{k+1}) 
&\geq \nabla f_k^\top (y_k-x_k) + \frac{\mcal}{2}\|y_k-x_k\|^2 + \frac{\eta}{2} (2-\Lcal \eta)\|\nabla f_k\|^2 - \frac{\Lcal \eta^2}{2} \|\varepsilon_k\|^2 + \eta(1-\Lcal \eta) \nabla f_k^\top \varepsilon_k \nonumber\\[2pt]
&\stackrel{\mathclap{\eqref{eq:NAM-2}}}{=} \;\; \begin{bmatrix}
\xi_{k} \\
\omega_k
\end{bmatrix}^\top
\Xb_1
\begin{bmatrix}
\xi_{k} \\
\omega_k
\end{bmatrix}
- \frac{\Lcal \eta^2}{2} \|\varepsilon_k\|^2 + \eta(1-\Lcal \eta) \nabla f_k^\top \varepsilon_k,
\end{align}
where
\begin{equation*}
\Xb_1 =\frac{1}{2}
\begin{bmatrix}
\beta^{2} \mcal & -\beta^{2} \mcal & -\beta & 0 \\
-\beta^{2} \mcal & \beta^{2} \mcal & \beta & 0 \\
-\beta & \beta & \eta(2-L \eta) & 0 \\
0 & 0& 0 & 0
\end{bmatrix} \otimes \Ib_d.
\end{equation*}
The inequality \eqref{eq:NAM-IQC-1} is the first IQC that we derive for the proof. For deriving the second IQC, we again apply the convexity of $f$ and have
\begin{equation} \label{eq:class-3}
f(x_\star) - f(x_k) \stackrel{\eqref{equ:class}}{\geq} \nabla f_k^\top  (x_\star - x_k)  +\frac{\mcal}{2} \|x_\star - x_k\|^2.
\end{equation}
Noting that
\begin{equation} \label{eq:NAM-1}
x_\star-x_k \stackrel{\eqref{eq:NAM-2}}{=} (x_\star -y_k) -\beta (y_k-x_\star) + \beta (y_{k-1}-x_\star) = -(1+\beta) (y_k-x_\star) + \beta (y_{k-1}-x_\star),
\end{equation}
we sum up \eqref{eq:class-3} and \eqref{eq:class-2} and obtain
\begin{align}\label{eq:NAM-IQC-2}
f(x_\star) - f( y_{k+1})
&\geq \nabla f_k^\top  (x_\star - x_k) \ +\frac{\mcal}{2} \|x_\star - x_k\|^2  + \frac{\eta}{2} (2-\Lcal \eta)\|\nabla f_k\|^2 \nonumber\\[2pt]
& \ \ \ - \frac{\Lcal \eta^2}{2} \|\varepsilon_k\|^2 + \eta(1-\Lcal \eta) \nabla f_k^\top \varepsilon_k \nonumber\\[2pt]
&\stackrel{\mathclap{\eqref{eq:NAM-1}}}{=}\;\; \begin{bmatrix}
\xi_{k} \\
\omega_k
\end{bmatrix}^\top
\Xb_2
\begin{bmatrix}
\xi_{k} \\
\omega_k
\end{bmatrix}
- \frac{\Lcal \eta^2}{2} \|\varepsilon_k\|^2 + \eta(1-\Lcal \eta) \nabla f_k^\top\varepsilon_k,
\end{align}
where
\begin{equation*}
\Xb_{2}=\frac{1}{2}
\begin{bmatrix}
(1+\beta)^{2} \mcal & -\beta(1+\beta) \mcal & -(1+\beta) &0 \\
-\beta(1+\beta) \mcal & \beta^{2} \mcal & \beta &0  \\
-(1+\beta) & \beta & \eta(2-L \eta) & 0\\
0 & 0& 0 & 0
\end{bmatrix}
\otimes \Ib_d.
\end{equation*}
The inequality \eqref{eq:NAM-IQC-2} is the second IQC that we derive for the proof. Now we~\mbox{combine}~\eqref{eq:NAM-IQC-1}~and~\eqref{eq:NAM-IQC-2} to obtain the final result. Noting that $\eta \leq 1/\mcal$ implies $0\leq \rho^2\leq 1$, we multiply \eqref{eq:NAM-IQC-1} by $\rho^2$ and~\eqref{eq:NAM-IQC-2}  by $1-\rho^2$, sum them up, and obtain
\begin{multline} \label{eq:IQCs-NAM}
\rho^{2}\left(f(y_{k})-f(x_\star)\right)-\left(f(y_{k+1})-f(x_\star)\right) \\[2pt]
\geq
\begin{bmatrix}
\xi_{k} \\
\omega_k
\end{bmatrix}^\top
\left(\rho^{2} \Xb_{1}+\left(1-\rho^{2}\right) \Xb_{2}\right)
\begin{bmatrix}
\xi_{k} \\
\omega_k
\end{bmatrix}
-\frac{\Lcal \eta^{2}}{2} \|\varepsilon_k\|^{2}+\eta(1-L \eta) \nabla f_k^{\top} \varepsilon_{k}.
\end{multline}
With $\beta$, $\rho^2$, and $\Pb$ given in the statement of the theorem, by direct calculation we have
\begin{multline} \label{eq:NAM-DT}
\begin{bmatrix}
\Ab^\top \Pb \Ab - \rho^2 \Pb & \Ab^\top \Pb \Bb \\
\Bb^\top \Pb \Ab & \Bb^\top \Pb \Bb
\end{bmatrix}
- \left(\rho^{2} \Xb_{1}+\left(1-\rho^{2}\right) \Xb_{2}\right)  \\
=
\begin{bmatrix}
- \tau(\eta) & \tau(\eta)  & 0 & -\frac{1+\eta\mcal}{2(1+\sqrt{\eta \mcal})}  \\
\tau(\eta)& -\tau(\eta) & 0 & \frac{1-\sqrt{\eta \mcal}}{2(1+\sqrt{\eta \mcal})}  \\
0& 0 & -\frac{\eta(1-\Lcal \eta )}{2} &  \frac{\eta}{2}\\
-\frac{1+\eta\mcal}{2(1+\sqrt{\eta \mcal})}  & \frac{1-\sqrt{\eta \mcal}}{2(1+\sqrt{\eta \mcal})} & \frac{\eta}{2} & \frac{\eta}{2}
\end{bmatrix}
\otimes \Ib_d \eqqcolon \Cb,
\end{multline}
where  $\tau(\eta) = (\mcal(1-\sqrt{\eta \mcal})^3)/(2\eta \mcal + 2\sqrt{\eta \mcal})$.
Multiplying $(\xi_k^\top \omega_k^\top)$ and $(\xi_k^\top \omega_k^\top)^\top$ on the left and~right, respectively, we further have
\begin{align*}
\EE&\sbr{\xi_{k+1}^\top\Pb\xi_{k+1}   - \rho^{2} \xi_k^\top \Pb \xi_k } \;\;
\stackrel{\mathclap{\eqref{eq:NAM-DT}}}{=} \;\;\EE
\begin{bmatrix}
\xi_{k} \\
\omega_k
\end{bmatrix}^\top
\left( \rho^{2} \Xb_{1}+\left(1-\rho^{2}\right) \Xb_{2} +
\Cb \right)
\begin{bmatrix}
\xi_{k} \\
\omega_k
\end{bmatrix}\\[2pt]
& \stackrel{\mathclap{\eqref{eq:IQCs-NAM}, \eqref{pequ:1}}}{\leq} \quad\; \EE \sbr{ \frac{\Lcal \eta^{2}}{2} \|\varepsilon_k\|^{2}-\eta(1-L \eta) \nabla f_k^{\top} \varepsilon_{k}  - \rbr{f(y_{k+1})-f(x_\star) } + \rho^2\rbr{ f(y_{k})-f(x_\star)} } \\[2pt]
& \hskip 2cm + \frac{  \eta }{2}  \EE \|\varepsilon_k\|^2 -\frac{\eta(1-\Lcal \eta)}{2} \EE \|\nabla f_k\|^2- \tau(\eta) \EE \|y_{k}- y_{k-1}\|^2 \\[2pt]
& \stackrel{\mathclap{\eqref{pequ:1}}}{=}\;\;\; \EE \sbr{  - (f(y_{k+1})-f(x_\star) )+ \rho^2(f(y_{k})-f(x_\star)) } \\[2pt]
& \hskip 2cm + \frac{ \eta(1+\Lcal\eta) }{2}   \EE \|\varepsilon_k\|^2 -\frac{\eta(1-\Lcal \eta)}{2}  \EE \|\nabla f_k\|^2- \tau(\eta) \EE  \|y_{k}- y_{k-1}\|^2\\[2pt]
&\stackrel{\mathclap{\eqref{eq:GC}}}{\leq} \;\;\;\EE \sbr{  - (f(y_{k+1})-f(x_\star) )+ \rho^2(f(y_{k})-f(x_\star)) }\\[2pt]
& \hskip 2cm + \frac{\eta}{2}\cbr{(1+\Lcal\eta)\delta - (1-\Lcal\eta)}\EE\|\nabla f_k\|^2 - \tau(\eta) \EE  \|y_{k}- y_{k-1}\|^2 + \frac{\eta(1+\Lcal\eta)}{2}\sigma^2.
\end{align*}
Here, we recall that $\Gcal_k$ is the $\sigma$-algebra containing all randomness $\{g_i\}_{i=0}^{k-1}$, which is generated by $\{y_j, x_j\}_{j=0}^k$ for \eqref{eq:SNV}. Thus, $\EE[\nabla f_k^\top \varepsilon_k \mid \Gcal_k] = 0$ for the third equality. This completes the proof.
\end{proof}

\vskip -0.5cm
Our present proof of Theorem \ref{thm:nam} is partially aligned with the one from \cite{Aybat2020Robust}; however, our derivation is sharper than theirs. To be specific, the inequality \eqref{eq:IQCs-NAM} is consistent with \cite[Lemma 4.5]{Aybat2020Robust}. From here, the convergence result \cite[(4.25)]{Aybat2020Robust} uses a loose inequality \cite[(4.24)]{Aybat2020Robust}, and their arguments on $c_0$ and $c$ (in their notation) are only about the existence. In contrast, we refine the inequality \cite[(4.24)]{Aybat2020Robust} with equality \eqref{eq:NAM-DT}, specialize $c=1$, and obtain a matrix $\Cb$ (cf. \eqref{eq:NAM-DT}). We then examine the explicit form of $\Cb$ to construct inequalities (i.e., some entries of $\Cb$ are utilized in~the~derivation), as opposed to simply using $\Cb\preceq \0$ as in \cite[(4.26)]{Aybat2020Robust}. Ultimately, our recursion \eqref{eq:NAM-recursive} is sharper~than \cite[(4.25)]{Aybat2020Robust} (which stems from the recursion presented between (4.10) and (4.11) in \cite{Aybat2020Robust}), in the sense that we include a second term $-\{\eta(1-\Lcal\eta)/2\}\EE\|\nabla f_k\|^2$ on the right hand side of \eqref{eq:NAM-recursive}. This second term is the precise mechanism by which we weaken the OUBV condition to the growth condition.

We notice that \cite[Theorem K.1]{Aybat2019Universally} briefly addresses the convergence of NAM under the growth~condition. Our analysis is sharper and strictly improves their result. In particular, \cite{Aybat2019Universally} requires $\kappa \geq 4$ and $\delta  \leq c \kappa^{-3/2}$ for a small constant $c$. With a small enough stepsize $\eta$, \cite{Aybat2019Universally} shows that NAM~exhibits~$1-1/(3\sqrt{\kappa})$ accelerated rate. However, their condition on $\delta$ is restrictive in general. Our result does not require any restriction on $\kappa$, but improves the noise tolerance level on $\delta$ significantly.~The~recursion \eqref{eq:NAM-recursive} is sharper than theirs even when $\delta=0$. 

We have two observations from Theorem \ref{thm:nam}. First, under the growth condition, the effect of the multiplicative noise on the convergence is $\delta\cbr{\eta(1+\Lcal\eta)/2} \EE\|\nabla f_k\|^2$. Second, NAM can inherently tolerate some amount of multiplicative noise, as revealed by the term $-\cbr{\eta(1-\Lcal\eta)/2} \EE\|\nabla f_k\|^2$ in \eqref{eq:NAM-recursive}. Since \vskip-10pt
\begin{equation*}
\frac{\eta(1-\Lcal\eta)}{2} \geq \frac{\eta(1+\Lcal\eta)\delta}{2} \Longleftrightarrow \eta \leq \frac{1-\delta}{(1+\delta)\Lcal},
\end{equation*} \vskip-5pt
\noindent the fastest rate of convergence can be obtained by setting $\eta = (1-\delta)/((1+\delta)\Lcal)$, which requires~$\delta\in [0, 1)$ though. We summarize the convergence rate in the next corollary.

\begin{corollary}\label{cor:nam}
Consider NAM in \eqref{eq:SNV} under the growth condition with constants $\delta\in[0, 1)$ and $\sigma^2$. Let $\eta = (1-\delta)/((1+\delta)\Lcal)$ and $\beta = (1-\sqrt{\eta\mcal})/(1+\sqrt{\eta\mcal})$ for NAM. Then \vskip-5pt
\begin{equation*}
\EE   V_{k+1}   \leq  \rbr{ 1- \sqrt{\frac{1-\delta}{(1+\delta)\kappa}} } \EE  V_{k}  + \frac{1-\delta}{(1+\delta)^2\Lcal}  \sigma^2,\quad \forall k\geq 1,
\end{equation*} \vskip-3pt
\noindent where $V_k$ is defined in \eqref{equ:NAM:V}. Furthermore, the iterate convergence is \vskip-5pt
\begin{equation*}
\EE \|y_{k}-x_\star \|^2
= \Ocal \left( \left( 1- \sqrt{\frac{1-\delta}{(1+\delta)\kappa}} \right)^{k} +\sigma^2 \right), \quad \quad \forall k \geq 0.
\end{equation*}
\end{corollary}

\begin{proof}
Under the setup of $\eta$, we have $\eta\cbr{(1-\delta) - (1+\delta)\Lcal\eta}/2 \geq  0$, $\tau(\eta) \geq  0$, and $\eta(1+\Lcal\eta)/2 = (1-\delta)/[(1+\delta)^2\Lcal]$. Applying Theorem \ref{thm:nam} immediately gives us the first statement. For the second part, we apply the error recursion in Theorem \ref{thm:nam} iteratively and obtain \vskip-7pt
\begin{equation}\label{eq:NAM-potential-ieq}
\EE\xi_{k+1}^\top\Pb\xi_{k+1} \leq\EE V_{k+1}\leq  \rho^{2k} \EE V_1 + \frac{\eta(1+\Lcal\eta)\sigma^2}{2(1-\rho^2)},\quad\quad\forall k\geq 0,
\end{equation}
where the first inequality is due to $f(y_{k+1}) - f(x_\star)\geq 0$. By the definition of $\Pb$ in Theorem \ref{thm:nam}, \vskip-15pt
\begin{equation*}
\xi_{k+1}^\top\Pb\xi_{k+1} = \|s_{k+1}\|^2, 
\quad \text{ where } \quad
s_{k+1}\coloneqq \sqrt{\frac{1}{2\eta}} (y_{k+1}-x_\star) + \left( \sqrt{\frac{\mcal}{2}} - \sqrt{\frac{1}{2\eta}}\right) (y_{k}-x_\star),
\end{equation*} \vskip-5pt
\noindent which implies for all $k\geq 0$ that \vskip-20pt
\begin{align*}
y_{k+1}-x_\star &=  (1-\sqrt{\mcal \eta})(y_{k}-x_\star) + \sqrt{2 \eta} s_{k+1} = (1-\sqrt{\mcal \eta})^k(y_1-x_\star) + \sqrt{2\eta} \sum_{t=1}^{k}(1-\sqrt{\mcal\eta})^{k-t}s_{t+1}\\
&= \rho^{2k}(y_1-x_\star) + \sqrt{2\eta} \sum_{t=1}^{k}\rho^{2(k-t)}s_{t+1}.
\end{align*}
Thus, for all $k\geq 0$, we have \vskip-15pt
\begin{align*}
\|y_{k+1}-x_\star\|^2
&= \nbr{\rho^{2k}(y_1-x_\star) + \sqrt{2\eta} \sum_{t=1}^{k}\rho^{2(k-t)}s_{t+1}}^2
\leq  2\rho^{4k}\|y_1-x_\star\|^2 + 4\eta \nbr{ \sum_{t=1}^{k}\rho^{2(k-t)}s_{t+1} }^2\\[2pt]
&\leq 2\rho^{4k}\|y_1-x_\star\|^2 + 4\eta \sum_{t=1}^{k}\rho^{0.5(k-t)}\cdot \sum_{t=1}^{k}\rho^{0.5(k-t)}\rho^{1.5^2(k-t)^2}\|s_{t+1}\|^2\\[2pt]
&\leq 2\rho^{4k}\|y_1-x_\star\|^2 + 4\eta \sum_{t=1}^{k}\rho^{0.5(k-t)}\cdot \sum_{t=1}^{k}\rho^{(0.5 + 1.5^2)(k-t)}\|s_{t+1}\|^2\\[2pt]
&\leq 2\rho^{4k}\|y_1-x_\star\|^2 + \frac{4\eta}{1-\sqrt{\rho}} \sum_{t=1}^{k}\rho^{2.75(k-t)}\|s_{t+1}\|^2,
\end{align*} \vskip-5pt
\noindent where the second and the third inequalities use the convexity of $\|\cdot\|^2$ and Jensen's inequality, while the fourth inequality uses the fact that $\rho^{1.5^2(k-t)^2}\leq \rho^{1.5^2(k-t)}$. Taking the full expectation on both sides, we further get for all $k\geq 0$, \vskip-20pt
\begin{align*}
\EE\|y_{k+1}-x_\star\|^2
&\leq  2\rho^{4k}\EE\|y_1-x_\star\|^2 + \frac{4\eta}{1-\sqrt{\rho}}\sum_{t=1}^{k}\rho^{2.75(k-t)}\EE\|s_{t+1}\|^2\\[2pt]
&\stackrel{\mathclap{\eqref{eq:NAM-potential-ieq}}}{\leq}\;\;2\rho^{4k}\EE\|y_1-x_\star\|^2 + \frac{4\eta}{1-\sqrt{\rho}}\sum_{t=1}^{k}\rho^{2.75(k-t)} \cbr{\rho^{2t}\EE V_1 + \frac{\eta(1+\Lcal\eta)\sigma^2}{2(1-\rho^2)}}\\[2pt]
&= 2\rho^{4k}\EE\|y_1-x_\star\|^2 + \frac{4\eta\rho^{2k}\EE V_1}{1-\sqrt{\rho}}\sum_{t=1}^{k}\rho^{0.75(k-t)} + \frac{2\eta^2(1+\Lcal\eta)\sigma^2}{(1-\sqrt{\rho})(1-\rho^2)}\sum_{t=1}^{k}\rho^{2.75(k-t)}\\[2pt]
&\leq  2\rho^{4k}\EE\|y_1-x_\star\|^2 + \frac{4\eta\EE V_1}{(1-\sqrt{\rho})(1-\rho^{0.75})}\rho^{2k} + \frac{2\eta^2(1+\Lcal\eta)}{(1-\sqrt{\rho})(1-\rho^2)(1-\rho^{2.75})}\sigma^2\\[2pt]
&= \Ocal \left( \rho^{2(k+1)} +\sigma^2 \right).
\end{align*} \vskip-5pt
\noindent This completes the proof by noting that $\EE\|y_0 - x_{\star}\|^2$ satisfies the above inequality trivially.
\end{proof}

Corollary \ref{cor:nam} states that NAM can accelerate constant-stepsize SGD under the growth condition when $\delta < 1$, by improving the dependence on the condition number in the rate from $1/\kappa$ to $1/\sqrt{\kappa}$.~To our knowledge, the recursion \eqref{eq:NAM-recursive} is the first result that shows NAM can tolerate~a~(mild)~multiplicative noise. Compared to the existing analyses \citep{Hu2017Dissipativity, Aybat2019Universally, Aybat2020Robust} that did not include the extra term $-\cbr{\eta(1-\Lcal\eta)/2} \EE\|\nabla f_k\|^2$ in the recursion to tolerate the multiplicative noise, our improvement is achieved by carefully designing the potential function and conducting a sharper analysis.

It is worth mentioning that Corollary \ref{cor:nam} does not contradict observations in \cite{Liu2020Accelerating},~where~the~authors constructed a LSR problem on which NAM cannot outperform SGD. For most LSR problems with Gaussian designs, as in \cite{Liu2020Accelerating}, we have $\delta \geq 1$, which contradicts assumptions of  Corollary \ref{cor:nam}.

\subsection{RMM under the growth condition}

We study the convergence of the robust momentum method (RMM) proposed in \cite{Cyrus2018Robust}, which has the following updating scheme \vskip-15pt
\begin{subequations}\label{eq:NAM-20}
\begin{align}
x_{k} &= y_{k} + \frac{\beta}{\eta } z_{k}, \label{eq:NAM-20a}\\
z_{k+1} &= \beta z_{k} - \frac{\eta}{\Lcal} g_k, \label{eq:NAM-20b}\\
y_{k+1} &= y_{k} + z_{k+1}. \label{eq:NAM-20c}
\end{align}
\end{subequations}
The work \cite{Cyrus2018Robust} proved that the deterministic RMM has linear convergence for strongly convex~objectives, and empirically showed that RMM is robust to the noise satisfying $\|\varepsilon_k\|\leq \delta\|\nabla f_k\|$ with small $\delta$. We study stochastic RMM and provide a theoretical guarantee for their empirical observations, under a more general growth condition \eqref{eq:GC} on the noise. Note that their noise model is only applicable for over-parameterized regime. We present the results in the following theorem and~corollary. The proofs are similar to Theorem \ref{thm:nam} and Corollary \ref{cor:nam}; thus we defer them to Appendix \ref{sup:pf:sec:3}.

\begin{theorem} \label{thm:rmm}
Consider RMM in \eqref{eq:NAM-20}. Suppose that the error $\varepsilon(x)$ satisfies the growth condition in Definition \ref{asp:GC} with constants $\delta$, $\sigma^2$. Suppose $f \in \Fcal( \mcal,\Lcal)$. For any $0<\theta\leq 1$, we denote $\tilde{\kappa} \coloneqq \kappa/\theta = \Lcal/(\theta \mcal)$ and let $\rho \in [1 - 1/\sqrt{\tilde{\kappa}},\sqrt{1-1/\tilde{\kappa}}]$ be any fixed, target convergence rate. We set $\eta, \beta$ in \eqref{eq:NAM-20} as
\begin{equation*}
\eta = \tilde{\kappa}(1-\rho)(1-\rho^2), \quad\quad\quad \beta=\frac{\tilde{\kappa}\rho^3}{\tilde{\kappa}-1},
\end{equation*}
and define quantities
\begin{equation*}
\lambda = \frac{\theta^2 \mcal^2(\tilde{\kappa}-\tilde{\kappa}\rho^2-1)}{2\rho(1-\rho)(1-\rho^2)^2 }, \quad  \nu=\frac{(1+\rho)(1-\tilde{\kappa}(1 -\rho)^2)}{2\rho}, \quad \Pb =
\begin{bmatrix}
(1-\rho^2)^2 & \rho^2(1-\rho^2) \\
\rho^2(1-\rho^2) & \rho^4
\end{bmatrix} \otimes \Ib_d.
\end{equation*}
Further, we let $\xi_k = ((y_k - x_\star)^\top, z_k^\top)^\top$ and define the potential function as
\begin{equation}\label{pequ:13}
V_k = \lambda \xi_k^\top\Pb\xi_k + h_\theta(x_{k-1}) + \frac{\delta(\Lcal - \theta\mcal)(\eta+\beta)}{\Lcal}\|\nabla f_{k-1}\|^2,
\end{equation}
where
\begin{equation}\label{equ:h}
h_\theta(x) = (\Lcal- \theta \mcal) \rbr{f(x) -f(x_\star) -\frac{\theta \mcal}{2}\|x-x_\star\|^2} - \frac{1}{2} \|\nabla f(x)- \theta \mcal (x-x_\star)\|^2.
\end{equation}
Then, for all $k\geq 1$, we have
\begin{equation} \label{eq:RMM-rec}
\EE V_{k+1}\leq \rho^2\EE V_k - \cbr{\nu \rbr{1 - 2\theta + \frac{1}{\tilde{\kappa}^2}} - \delta(1 - \nu + \eta + \beta)}\EE\|\nabla f_k\|^2 + \rbr{\rho^2(\eta+\beta) + 1-\nu}\sigma^2.
\end{equation}
\end{theorem}

\begin{proof}
See Appendix \ref{pf:thm:rmm}.
\end{proof}

Similar to Theorem \ref{thm:nam}, Theorem \ref{thm:rmm} reveals two things about RMM: under the growth condition, the effect of the multiplicative noise on the convergence is $\delta(1-\nu+\eta+\beta)\EE\|\nabla f_k\|^2$; however, RMM can tolerate this noise due to the term $-\nu(1-2\theta+1/\tilde{\kappa}^2) \EE\|\nabla f_k\|^2$. While these multipliers are more complex than those in Theorem \ref{thm:nam}, the next corollary simplifies the result and clearly illustrates that RMM is robust to a small $\delta$.

\begin{corollary}\label{cor:rmm}
Consider RMM in \eqref{eq:NAM-20} under the growth condition with constants $\delta\in[0, \\ 1/4)$ and $\sigma^2$. We consider two cases:
\begin{enumerate}[label=(\alph*),topsep=0pt]
\setlength\itemsep{0.0em}
\item if $\delta = 0$, then we set $\theta = 1$ and $\rho = 1 - 1/\sqrt{\kappa}$;
\item if $\delta\in (0, 1/4)$, then we set $\theta = 1/2 - \sqrt{\delta}$ and $\rho = 1 - \sqrt{2}\theta/\sqrt{\kappa}$.
\end{enumerate}
Then, we define other quantities including $\eta, \beta, \lambda, \Pb$ as in Theorem \ref{thm:rmm}. It holds that
\begin{equation}\label{pequ:12}
\EE V_{k+1}\leq \rho^2\EE V_k  + (3\rho^2+1)\sigma^2, \quad \forall k\geq 1,
\end{equation}
where $V_k$ is defined in \eqref{pequ:13}. Furthermore, the iterate convergence is
\begin{equation*}
\EE \|y_{k}-x_\star \|^2\leq \Ocal \left( \rho^{2k} +\sigma^2 \right),\quad \quad \forall k \geq 0.
\end{equation*}
\end{corollary}

\begin{proof}
See Appendix \ref{pf:cor:rmm}.
\end{proof}

Compared to deterministic analysis of RMM in \cite{Cyrus2018Robust}, our stochastic analysis has four major~differences. (i) Our analysis starts from the linear matrix equality \eqref{eq:LME-NAM-20-2}, which is similar to our analysis of NAM in \eqref{eq:NAM-DT}. The analysis of \cite{Cyrus2018Robust} starts from the smoothness and strong convexity conditions, and does not clearly specify the control matrices $\Xb_1, \Xb_2, \Pb$ (cf. \cite[Theorem 1]{Cyrus2018Robust}). (ii) We introduce~a relaxation factor $\theta\in(0, 1]$ to address technical challenges brought by the growth condition. Such~a factor is original and critical to show the accelerated rate: when $\delta = 0$, we set $\theta=1$, while~when $\delta>0$, we set $\theta = 0.5-\sqrt{\delta}$ depending on $\delta$. (iii) Our potential function \eqref{pequ:13} also depends on the multiplicative noise factor $\delta$, and the recursion \eqref{eq:RMM-rec} is significantly different from \cite[(10)]{Cyrus2018Robust}. (iv) More technically, our stochastic analysis of the recursion at the $k$-th step relies on the past gradient noise $\epsilon_{k-1}$. Thus, we have to handle product terms like $(x_k-x_\star)^T\epsilon_{k-1}$ and $\nabla f_k^T\epsilon_{k-1}$, which do not mean zero when taking conditional expectation. We refer to the derivation~between~\eqref{eq:IQC-stochastic} and~\eqref{pequ:8} for more details. Such a difficulty does not appear in the existing~NAM or RMM  analysis. Also,~\cite{Cyrus2018Robust} employs deterministic analysis that does not involve any noise terms.

Corollary \ref{cor:rmm} states that RMM can accelerate constant-stepsize SGD under the growth condition when $\delta < 1/4$, by improving the dependence on the condition number in the rate from $1/\kappa$ to~$1/\sqrt{\kappa}$. Our result strongly coincides with the empirical observation in \cite[Figure 3]{Cyrus2018Robust}, where the authors observed that $\delta \approx 0.26$ is a transition point of RMM between the rates $1 - \Ocal(1/\kappa)$ and $1-\Ocal(1/\sqrt{\kappa})$.

An interesting observation is that we use two different sets of parameters for $\delta = 0$ and $\delta\in(0, 1/4)$ to let $\rho$ be as small as possible. In particular, when $\delta = 0$, we obtain $\rho = 1 - 1/\sqrt{\kappa}$ convergence rate (same for NAM), which recovers the deterministic analysis in \cite{Cyrus2018Robust}. However, when $\delta\in(0, 1/4)$, we require a different set of parameters, which does not approach to the set of $\delta = 0$ when $\delta\rightarrow 0$. In fact, letting $\delta \rightarrow 0$, we should set $\theta = 1/2$ and $\rho = 1 - 1/\sqrt{2\kappa}$ for $\delta = 0$. Such a setup is also valid and enjoys the recursion \eqref{pequ:12}. However, this rate is not as good as $\rho = 1 - 1/\sqrt{\kappa}$. The inconsistency comes from the extra parameter $\theta\in(0, 1]$, which we deliberately introduce into the analysis as it allows us to establish robustness of RMM to the multiplicative noise. Without introducing $\theta$, that is, $\theta = 1$, the term $-\nu (1-2\theta + \theta^2/\kappa^2)\EE\|\nabla f_k\|^2$ in \eqref{eq:RMM-rec} has a nonnegative multiplier, since $\nu\geq 0$. As a result, we would not be able to show that the scheme could tolerate $\delta>0$.

We further illustrate the inconsistent setups of $\delta=0$ and $\delta\in(0,1/4)$ as follows. Consider offsetting the multiplicative noise with $\kappa>1$. By \eqref{eq:RMM-rec}, we require \vskip-8pt
\begin{equation}\label{npeqeu:1}
\nu g(\theta) \coloneqq \nu \left(1-2\theta + \frac{\theta^2}{\kappa^2} \right)\geq (1-\nu +\eta+\beta)\delta.
\end{equation} \vskip-5pt
\begin{enumerate}[label=(\alph*),topsep=0pt]
\setlength\itemsep{0.0em}
\item When $\delta = 0$, we only need $\nu g(\theta) = 0$.
\begin{itemize}
	\setlength\itemsep{0.0em}
	\item If $\nu = 0$, we obtain $\rho = 1 - 1/\sqrt{\tilde{\kappa}} = 1-\sqrt{\theta}/\sqrt{\kappa}$. To let $\rho$ be the smallest, we choose the largest $\theta$, which is $\theta = 1$, and get $\rho = 1 - 1/\sqrt{\kappa}$.
	\item If $g(\theta) = 0$, we can easily see that $g(\theta)$ within $(0, 1]$ has a single root $\theta_0<1$. Then, the smallest rate in this case is $\rho =  1-\sqrt{\theta_0}/\sqrt{\kappa}$, which is greater than $1 - 1/\sqrt{\kappa}$.
\end{itemize}
Thus, if $\delta = 0$, we should let $\nu = 0$ instead of $g(\theta) = 0$ to have the smallest $\rho$.

\item When $\delta>0$, observing that $\nu\leq 1$ and $\eta, \beta >0$, the right hand side of \eqref{npeqeu:1} is positive. Thus, we must have $\nu g(\theta)>0$, which implies $\nu>0$ and $g(\theta)>0$. This already implies the inconsistent setup with (a) as $\delta\rightarrow 0$. In particular, $g(\theta)>0 \Longleftrightarrow 0<\theta<\theta_0<1$; hence, $\delta\rightarrow 0$ cannot lead to $\theta\rightarrow 1$ because $g(1)<0$.	
\end{enumerate}

Finally, we remark on the potential function \eqref{pequ:13}, which depends on $h_{\theta}(\cdot)$ and $\|\nabla f\|^2$. While this makes the potential function look more complicated compared to \eqref{equ:NAM:V} of NAM, the terms $h_{\theta}(\cdot)$ and $\|\nabla f\|^2$ are positive and, thus, are negligible when we study the iterate convergence. The~iterate convergence only relies on the quadratic form $\lambda\xi_k^\top\Pb\xi_{k}$.

\subsection{Comparison of SGD, NAM, and RMM}

We have shown that SGD, NAM, and RMM all converge under the growth condition, with the rates that are the same as theirs under the OUBV condition. However, the methods allow different levels of $\delta$. Table \ref{tab:SGD-NAM-RMM} summarizes the error recursions of (constant-parameters) SGD, NAM, and RMM. Note that the rates in Table \ref{tab:SGD-NAM-RMM} are consistent with the convergence rates of the iterates, but different methods have different~potential~functions~$V_k$.

NAM and RMM are two accelerated methods that use different momentum designs. NAM uses the same parameter $\beta$ of momentum $y_k - y_{k-1}$ for $x_k$ and $y_{k+1}$; while RMM employs distinct~parameters $\beta/\eta$ and $\beta$ of momentum $y_k - y_{k-1}$ for $x_k$ and $y_{k+1}$, respectively. Both methods achieve the accelerated rate $1-\Ocal(1/\sqrt{\kappa})$ for small $\delta$. The rate of SGD holds for any $\delta\geq 0$, but is not optimal in terms of $\kappa$. Since $\sqrt{(1-\delta)/(1+\delta)} < 1/(1+\delta)$, the rate of NAM slows down faster than that of SGD as $\delta$ increases. When $\delta= 0$, RMM has a better rate than NAM: the former is $(1-1/\sqrt{\kappa})^2$, while the latter is $1-1/\sqrt{\kappa}$. However, when $\delta\in(0, 1/4)$, NAM has a better rate than RMM since $1/\sqrt{2} - \sqrt{2\delta} <\sqrt{(1-\delta)/(1+\delta)}$.

When $\sigma^2=0$, \cite{Cyrus2018Robust} empirically showed that RMM is robust to small multiplicative noise while NAM is less robust. However, as shown in Table \ref{tab:SGD-NAM-RMM}, the variance of RMM is the worst among the three methods. Together with the observation in \cite{Cyrus2018Robust}, we reasonably believe that RMM is~more~suitable for over-parameterized models. In addition, the experiments of \cite{Cyrus2018Robust} did not choose algorithm's parameters such as $\beta, \eta$ based on the noise level $\delta$; while Corollaries \ref{cor:nam} and \ref{cor:rmm} suggest to adjust parameters to mitigate multiplicative noise and obtain better convergence behavior. Thus, our two papers conduct the comparisons in different senses. They fixed parameters setup, and compared~the performance of the setup on different $\delta$; while we fix $\delta$ and set parameters accordingly, and compare the performance of different methods with optimal parameters. 

We examine the dual averaging-based methods in the next section, which have received less~attention in the literature on stochastic analysis. Our analysis of the dual averaging methods follows~the same structure as the present section, where we first establish an error recursion in the form of~\eqref{equ:infor:result}, and then specify proper algorithm's parameters to attain the accelerated rate. We note~that the~unified form of recursion enables us to examine different bias-variance trade-offs of all the methods~and facilitates the comparison of their robustness ability to the multiplicative noise (see \cref{sec:4.4}).

\begin{table}
\TABLE
{Summary of convergence rates of (constant-stepsize) SGD, NAM and RMM. \label{tab:SGD-NAM-RMM}}
{\vspace{-4mm}\begin{tabular}{cccc}
\hline
\up \down Algorithm & Bias term ($\EE V_k$)   & Variance term ($\sigma^2$)  & Requirement for $\delta$  \vspace{2pt} \\
\hline
\up \down SGD   & $1- 1/[(1+\delta)\kappa]$  & $1/[(1+\delta)^2\Lcal^2]$  & $0\leq \delta $ \vspace{2pt} \\
\hline
\up \down NAM & $1- \sqrt{(1-\delta)/[(1+\delta)\kappa]}$   & $ (1-\delta)/[(1+\delta)^2\Lcal]$  & $0\leq \delta<1$  \vspace{2pt}\\
\hline
\up \down \multirow{2}{*}{RMM}    & $\rbr{1 - 1/\sqrt{\kappa}}^2$ & $4- 3/\sqrt{\kappa}$ & $\delta=0$\\[-4pt]
\up \down & $ [1 - (1/\sqrt{2} - \sqrt{2\delta})/\sqrt{\kappa} ]^2$ & $4 - [3(1/\sqrt{2}-\sqrt{2 \delta})]/\sqrt{\kappa}$ & $ 0< \delta < 1/4$  \vspace{2pt} \\
\hline
\end{tabular} \vspace{-3mm}
}
{}
\end{table}

\section{Dual Averaging-Based Methods} \label{sec:4}

We study two different schemes to \mbox{accelerate}~Dual~Averaging Method (DAM): accelerated dual averaging method (DAM+) and implicit DAM+ (iDAM+). Let us first briefly introduce DAM \citep{Xiao2010Dual, Cohen2018Acceleration, Diakonikolas2018Accelerated, Diakonikolas2019Approximate}. Given the past iterates $\{x_i\}_{i=1}^k$ and corresponding stochastic gradients $\{g_i\}_{i=1}^k$, DAM updates the iterate as \vskip-5pt
\begin{equation} \label{eq:DAM}
x_{k+1} = \arg \min_x \left\{ \sum_{i=1}^k a_i \left[g_i^\top(x-x_i) +\frac{\mcal}{2} \|x-x_i\|^2  + \frac{m_0}{2} \|x-x_1\|^2 \right] \right\},
\end{equation}\vskip-3pt
\noindent where $\mcal$ is the strong convexity parameter, $m_0\geq0$ is a constant, and $\{a_i\}_{i=1}^k$ is a positive weight sequence. When $f$ is only convex ($\mcal=0$), it is important to have $m_0>0$ for regularization. However, since $f$ is assumed to be strongly convex in this paper, we let $m_0 =0$ to simplify the presentation. When $a_1 = \dots = a_{k-1}=0$ and $a_k=1$, \eqref{eq:DAM} with $m_0=0$ recovers SGD.

The following preliminary lemma will be used in the analysis of DAM+ and iDAM+. We slightly abuse notation $h_{\theta}(x)$ from Theorem \ref{thm:rmm}.

\begin{lemma}\label{lemma:lower-bd}
Given the iterate sequence $\{x_i\}_{i=1}^k$ and the weight sequence $\{a_i\}_{i=1}^k$, let us define \vskip-15pt
\begin{equation}\label{eq:hk}
A_k = \sum_{i=1}^{k} a_i,\;\; z_{k+1} = -\sum_{i=1}^{k} a_i g_i, \;\; \phi_k(x) = \sum_{i=1}^{k} \frac{a_i\mcal}{2} \|x-x_i\|^2, \;\; h_k(x) = \sum_{i=1}^{k}a_{i}g_{i}^{\top}(x-x_{i}) + \phi_k(x).
\end{equation}
The following statements hold for any $k\geq 1$:
\begin{enumerate}[label=(\alph*),topsep=2pt]
\setlength\itemsep{0.2em}
\item $\nabla \phi_k^\star(z_{k+1}) = \arg \min_x h_k(x)$;
\item $ h_k(\nabla \phi_k^\star(z_{k+1})) = - \phi_k^\star(z_{k+1}) -  \sum_{i=1}^{k} a_i g_i^\top x_i$;
\item we have $f(x_\star) \geq L_{k+1}$, where \vskip-10pt
\begin{equation}\label{eq:DA-ld2}
L_{k+1}= \frac{1}{A_k} \sum_{i=1}^k a_i \left[f(x_i) - \varepsilon_{i}^\top (x_\star -x_i) \right]  - \frac{1}{A_k} \phi_k^\star(z_{k+1}) - \frac{1}{A_k} \sum_{i=1}^{k} a_i g_i^\top x_i.
\end{equation}
\end{enumerate}
\end{lemma}

\begin{proof}
By Lemma \ref{lemma:BD}(a), we have
\begin{align*}
\nabla \phi_k^\star(z_{k+1})
& =  \argmax_x \{  z_{k+1}^\top x - \phi_k(x)\} = \argmin_x \{  - z_{k+1}^\top x + \phi_k(x)\} \\[2pt]
&\stackrel{\mathclap{\eqref{eq:hk}}}{=} \;\; \argmin_x \cbr{  \sum_{i=1}^k  a_i g_i^\top x  + \phi_k(x)  } = \argmin_x \cbr{ \sum_{i=1}^k  a_i g_i^\top x  + \phi_k(x)  -  \sum_{i=1}^k a_i g_i^\top x_i } \\[2pt]
&= \argmin_x h_k(x),
\end{align*}
which proves (a). For (b), we have  \vskip-15pt
\begin{align*}
h_k(\nabla \phi_k^\star(z_{k+1})) 
&\stackrel{\mathclap{(a)}}{=}  \min_x h_k(x)  =  \min_x \cbr{ \sum_{i=1}^k a_i g_i^\top (x-x_i) + \phi_{k}(x) } \\[2pt]
&= \min_x\left\{ -z_{k+1}^\top x +\phi_k(x) \right\} - \sum_{i=1}^{k} a_i g_i ^\top x_i = - \max_x\left\{ z_{k+1}^\top x -\phi_k(x) \right\} - \sum_{i=1}^{k} a_i g_i ^\top x_i  \\[2pt]
& \stackrel{\mathclap{\eqref{eq:convex-conjugate}}}{=}\;\; -\phi_k^\star(z_{k+1})  - \sum_{i=1}^{k} a_i g_i ^\top x_i.
\end{align*} \vskip-5pt
\noindent Finally, by the strong convexity of $f$, we have
\begin{equation} \label{eq:lower-bd}
f(x_\star) \geq  f(x_k) + \nabla f_k^\top (x_\star - x_k) + \frac{ \mcal}{2} \|x_\star - x_k\|^2,\quad \forall k\geq 1.
\end{equation} \vskip-5pt
\noindent Taking a weighted average on the right hand side, we have \vskip-15pt
\begin{align*}
f(x_\star)
& \stackrel{\mathclap{\eqref{eq:lower-bd}}}{\geq} \; \frac{1}{A_k} \sum_{i=1}^k a_i \sbr{f(x_i) + \nabla f_i^\top (x_\star-x_i)+\frac{\mcal}{2} \| x_\star - x_i \|^2 } \\[2pt]
&=\frac{1}{A_k} \sum_{i=1}^k a_i \left[f(x_i) - \varepsilon_i^\top (x_\star-x_i) \right] + \frac{1}{A_k}  \sum_{i=1}^k a_i \left[g_i^\top (x_\star-x_i)+ \frac{\mcal}{2} \| x_\star - x_i \|^2  \right] \\[2pt]
&\geq \frac{1}{A_k} \sum_{i=1}^k a_i \left[f(x_i) - \varepsilon_i^\top (x_\star -x_i)  \right]  + \frac{1}{A_k} \min_x \left\{ \sum_{i=1}^k a_i \left[g_i^\top (x-x_i)+ \frac{\mcal}{2} \| x - x_i \|^2  \right] \right\} \\[2pt]
&\stackrel{\mathclap{\eqref{eq:hk}}}{=}  \;\frac{1}{A_k} \sum_{i=1}^k a_i \left[f(x_i) - \varepsilon_i^\top (x_\star -x_i)  \right]  + \frac{1}{A_k} \min_x h_k(x).
\end{align*} \vskip-5pt
\noindent Then, (c) follows by applying (a) and (b). This completes the proof.
\end{proof}

By Lemma \ref{lemma:lower-bd}(a), the updating rule \eqref{eq:DAM} with $m_0=0$ can be rewritten as
\begin{equation} \label{eq:DAM-1}
x_k = \nabla \phi_{k-1}^\star (z_{k}), \quad\quad\quad z_{k+1} = z_k - a_k g_k.
\end{equation}
\noindent We further present an analytical framework used in this section.

\vskip5pt
\noindent\textbf{Analysis framework.}
The accelerated DAMs are often studied based on the duality gap technique, developed in \cite{Diakonikolas2019Approximate}. Different from the IQC technique used in \cref{sec:3}, the duality gap~technique focuses on the function value convergence by constructing upper and lower bounds of $f$, and allows varying parameters. There are three main steps:\\
\noindent \underline{Step 1}:
given an \textit{output iterate} $y_k$ of the algorithm, we let the upper bound be $U_k = f(y_k)$, while~we construct a suitable lower bound $L_k$ of $f(x_\star)$, i.e., $f(x_\star) \geq L_k$. Thus, the gap $V_k = U_k - L_k$ satisfies $f(y_k) - f(x_\star) \leq V_k$.\\
\noindent \underline{Step 2}: given the sequence $A_k$ (recall the definition of $A_k$ in \eqref{eq:hk}), we establish the upper bound~of $A_k U_{k+1} -A_{k-1} U_{k} \leq \text{UpperBound}$ and the lower bound of $A_k L_{k+1} - A_{k-1} L_{k} \geq \text{LowerBound}$.\\
\noindent \underline{Step 3}: Combining the above two steps, we get
\begin{equation*}
A_{k} V_{k+1} - A_{k-1} V_{k} = (A_k U_{k+1} -A_{k-1} U_{k} )- (A_k L_{k+1} - A_{k-1} L_{k} ) \leq \text{UpperBound} - \text{LowerBound},
\end{equation*}
which leads to a recursion $f(y_{k+1}) - f(x_\star) \leq V_{k+1} \leq V_k A_{k-1}/ A_{k} + (\text{UpperBound} - \text{LowerBound})/A_k$; and the convergence rate is $A_{k-1}/ A_{k} = 1-a_k/ A_{k}$.

We emphasize that we aim at bounding $f(y_k) - f(x_\star)$ rather than $f(x_k) - f(x_\star)$ in the following presentation, where $y_k$ is also the algorithm variable and usually regarded as the algorithm output. We take DAM+ as an example and show the details of three steps, while we only present the result of iDAM+ for conciseness and defer its proofs to Appendix \ref{sup:pf:sec:4}.

\subsection{DAM+ under the growth condition}

DAM+ combines the SGD and DAM~\mbox{iterates}.~In particular, DAM+ is given by the following scheme:
\begin{subequations} \label{eq:ADA}
\begin{align}
x_{k} &= \frac{A_{k-1}}{A_{k}} y_{k} + \frac{a_{k}}{A_{k}} \nabla \phi_{k}^\star (z_{k}), \label{eq:ADAa}\\
z_{k+1} &= z_{k} - a_k g_k, \label{eq:ADAb}\\
y_{k+1} &= x_{k} - \eta_k g_k. \label{eq:ADAc}
\end{align}
\end{subequations}
We see that $x_k$ is the convex combination of the SGD iterate $y_k$ and the DAM iterate $\nabla \phi_{k}^\star (z_{k})$,~which leads to a better prediction of $x_\star$. Compared to \eqref{eq:DAM-1}, $x_k$ is defined using $\nabla \phi_{k}^\star (z_{k})$ that depends on $x_k$ itself, instead of using $\nabla \phi_{k-1}^\star (z_{k})$. By Lemma \ref{lemma:BD}(a), $\nabla \phi_{k}^{\star}(z)$ has a closed form:
\begin{equation} \label{eq:closed-form}
\nabla \phi_{k}^{\star}(z)=\frac{z + \mcal \sum_{i=1}^{k} a_{i} x_i}{A_k \mcal} = \frac{z + \mcal \sum_{i=1}^{k-1} a_{i} x_i}{m_k} + \frac{a_k}{A_k}x_k, \quad \text{ where }\; m_k \coloneqq A_k \mcal.
\end{equation}
Since $a_k/A_k <1$, \eqref{eq:ADAa} is always solvable for $x_k$. The similar modification replacing $\nabla \phi_{k-1}^\star (z_{k})$ by $\nabla \phi_{k}^\star (z_{k})$ is adopted in \cite[Appendix C]{Diakonikolas2019Approximate} to improve the factor in the convergence bound of DAM+.

We mention that DAM+ is a popular algorithm, and its formulation in different papers is slightly different. Our formulation is adapted from \cite[(ASC)]{Diakonikolas2019Approximate} and \cite[(AGD)]{Diakonikolas2018Accelerated}, while similar formulations without involving $\phi_k^\star$ can be found in \cite{Jain2018Accelerating, Vaswani2019Fast}. In particular, using the closed form of $\nabla \phi_{k}^{\star}(z)$ in~\eqref{eq:closed-form}, we establish a relationship between DAM+ and the accelerated methods in \cite{Jain2018Accelerating, Vaswani2019Fast}.

\begin{proposition} \label{prop:ADAM-other}
Consider the following two schemes for $k\geq 1$:
\begin{equation} \label{eq:ADA-1}
\begin{split}
x_{k} = \frac{A_{k-1}}{A_{k}} y_{k} + \frac{a_{k}}{A_{k}} \nabla \phi_{k-1}^\star (z_{k}),   \ \ \
z_{k+1} = z_{k} - a_k g_k,  \ \ \
y_{k+1} = x_{k} - \eta_k g_k,
\end{split}
\end{equation}
and
\begin{equation} \label{eq:SGD-N}
\begin{split}
x_{k}' = \frac{A_{k-1}}{A_{k}} y_{k}' + \frac{a_{k}}{A_{k}} z^\prime_k,   \ \ \
z^\prime_{k+1} =  \frac{\mcal a_k }{m_k} x_k' +  \frac{m_{k-1}}{m_k} z^\prime_{k} - \frac{a_k}{m_k} g_{k},  \ \ \
y_{k+1}' = x_{k}' - \eta_k g_k,
\end{split}
\end{equation}
where $m_k = A_k \mcal$. For consistency, we let $m_0 = a_0 = A_0=0$. If $z_1 = z_1' = 0$, then the two sequences $(x_1, z_2, y_2, x_2, z_3, y_3, \ldots)$ and $(x_1', z_2', y_2', x_2', z_3', y_3', \ldots)$ generated by \eqref{eq:ADA-1} and \eqref{eq:SGD-N} have the relation
\begin{equation*}
x_j = x_j', \quad \forall j\geq 1; \quad\quad\quad y_j = y_j', \quad \nabla\phi_{j-1}^\star (z_{j}) = z_{j}', \quad \forall j\geq 2,
\end{equation*}
where $\nabla\phi_j^\star(\cdot)$ is defined in \eqref{eq:closed-form}.
\end{proposition}

\begin{proof}
We prove the result by induction. Letting $0/0 = 0$ for consistency, we have $x_1 \stackrel{\eqref{eq:ADA-1}}{=} \nabla \phi_0^\star(z_1) = \nabla \phi_0^\star(0) \stackrel{\eqref{eq:closed-form}}{=} 0 = z_1' \stackrel{\eqref{eq:SGD-N}}{=} x_1'$. Suppose $x_i = x_i'$, $\forall i\leq k$ for some $k \geq 1$. Then it suffices to show \vskip-8pt
\begin{equation*}
y_{k+1} = y_{k+1}', \quad \nabla\phi_k^\star(z_{k+1}) = z_{k+1}',
\end{equation*}
since, by \eqref{eq:ADA-1} and \eqref{eq:SGD-N}, we would have $x_{k+1} = x_{k+1}'$, which completes the induction step. Certainly, \vskip-8pt
\begin{equation*}
y_{k+1} \stackrel{\eqref{eq:ADA-1}}{=} x_k - \eta_kg_k = x_k' - \eta_kg_k \stackrel{\eqref{eq:SGD-N}}{=} y_{k+1}'.
\end{equation*}
Further, applying the recursion of $z_k$ in \eqref{eq:ADA-1}, we have $z_{k+1} = -\sum_{i=1}^{k}a_ig_i$, $\forall k\geq 1$. Then, by \eqref{eq:closed-form}, \vskip-8pt
\begin{equation*}
\nabla \phi_{k}^\star (z_{k+1})
= \frac{ - \sum_{i=1}^{k} a_i g_i + \mcal \sum_{i=1}^{k} a_{i} x_i}{m_k}.
\end{equation*}
Applying the recursion of $z_k'$ in \eqref{eq:SGD-N}, \vskip-18pt
\begin{align*}
z^\prime_{k+1}
&= \frac{\mcal a_k }{m_k} x_k'  - \frac{a_k}{m_k} g_{k} + \frac{m_{k-1}}{m_k} z^\prime_{k}
= \frac{\mcal a_k}{m_k} x_k'  - \frac{a_{k}}{m_k} g_{k} + \frac{m_{k-1}}{m_k} \left(  \frac{\mcal a_{k-1}}{m_{k-1}} x_{k-1}'  - \frac{a_{k-1}}{m_{k-1}} g_{k-1} + \frac{m_{k-2}}{m_{k-1}} z^\prime_{k-1}  \right) \\[2pt]
&= \frac{\mcal \sum_{i=k-1}^{k} a_i x_i' - \sum_{i=k-1}^{k} a_i g_i }{m_k} +  \frac{m_{k-2}}{m_{k}} z^\prime_{k-1}
= \frac{\mcal \sum_{i=1}^{k} a_i x_i' - \sum_{i=1}^{k} a_i g_i }{m_k} +  \frac{m_0}{m_{k}} z^\prime_{1} \\[2pt]
&= \frac{\mcal \sum_{i=1}^{k} a_i x_i' - \sum_{i=1}^{k} a_i g_i }{m_k},
\end{align*}
where the last equality uses $m_0 = 0$. Comparing the above two displays and noting that $x_i = x_i'$,~$\forall i\leq k$, we have $\nabla \phi_k^\star(z_{k+1}) = z_{k+1}'$. This completes the proof.
\end{proof}

The difference between \eqref{eq:ADA-1} and \eqref{eq:ADA} is that $\nabla \phi_{k-1}^\star (z_{k})$ in \eqref{eq:ADA-1} replaces $\nabla \phi_{k}^\star (z_{k})$ in \eqref{eq:ADAa}.~From Proposition \ref{prop:ADAM-other}, we know that the term $\nabla\phi_k^\star(z_k)$ with the closed form \eqref{eq:closed-form} that accumulates all the past gradients can be computed in an online fashion, by defining a variable similar to $z_k'$ in~\eqref{eq:SGD-N}. In particular, by \eqref{eq:closed-form}, we have \vskip-22pt
\begin{equation*}
\nabla \phi_{k}^\star (z_{k}) \stackrel{\eqref{eq:closed-form}}{=} \frac{a_k}{A_k}x_k + z_k'',\quad\quad \text{with }\quad  z_k'' \stackrel{\eqref{eq:closed-form}}{\coloneqq} \frac{z_k + \mcal\sum_{i=1}^{k-1}a_ix_i}{m_k} \stackrel{\eqref{eq:ADA}}{=} \frac{-\sum_{i=1}^{k-1}a_ig_i + \mcal\sum_{i=1}^{k-1}a_ix_i}{m_k}. 
\end{equation*}
Thus, \eqref{eq:ADAa} and \eqref{eq:ADAb} of DAM+ can be rewritten as \vskip-8pt
\begin{equation*}
\rbr{1 - \frac{a_k^2}{A_k^2}}\cdot x_k = \frac{A_{k-1}}{A_k} y_k + \frac{a_k}{A_k}z_k'', \quad\quad z_{k+1}'' = \frac{m_k}{m_{k+1}}z_k'' + \frac{\mcal a_k}{m_{k+1}}x_k - \frac{a_k}{m_{k+1}}g_k.
\end{equation*}
Combing the above display with \eqref{eq:ADAc}, we know that DAM+ is also the \mbox{first-order}~SA~method~where a single noisy gradient $g_k$ is inquired at each step (similar for DAM and iDAM+).

Now we set the stage to establish the convergence of DAM+. We follow the three steps outlined before. For the first step, we let the upper bound $U_k = f(y_k)$, and let the lower bound $L_k$ be given in \eqref{eq:DA-ld2}. The next lemma completes the second step by bounding the differences $A_k U_{k+1} - A_{k-1} U_{k}$~and  $A_k L_{k+1} - A_{k-1} L_{k}$.

\begin{lemma} \label{lemma:ADAM}
Given $U_k=f(y_k)$ and $L_k$ in \eqref{eq:DA-ld2}, the following results hold for the scheme \eqref{eq:ADA}.
\begin{enumerate}[label=(\alph*),topsep=2pt]
\setlength\itemsep{0.2em}
\item  $A_k U_{k+1} - A_{k-1} U_{k} \leq  a_k f(x_k) + A_{k} (f(y_{k+1})- f(x_k)) + a_k \nabla f_k^\top (\nabla \phi_k^\star (z_{k}) - x_k)$;
\item $A_k L_{k+1} - A_{k-1} L_{k} = a_k f(x_k)   - a_k \varepsilon_{k}^\top x_\star   - a_k \nabla f_k^\top x_k - \phi_k^\star (z_{k+1}) + \phi_{k-1}^\star (z_{k})$.
\end{enumerate}
\end{lemma}

\begin{proof}
We have that \vskip-18pt
\begin{equation} \label{eq:ADA1}
(A_{k-1}+a_k) x_k=A_k x_k \stackrel{\eqref{eq:ADAa}}{=} A_{k-1} y_k + a_k\nabla \phi_{k}^\star (z_{k}) \Longrightarrow  A_{k-1} (y_k - x_k)  = a_k(x_k-\nabla \phi_k^\star (z_{k})).
\end{equation}\vskip-5pt
\noindent By the convexity of $f$, it holds that \vskip-12pt
\begin{equation} \label{eq:ada-conv}
A_{k-1} (f(y_{k})-f(x_{k})) \stackrel{\eqref{equ:class}}{\geq} A_{k-1} \nabla f_k^\top (y_{k} - x_{k}) \stackrel{\eqref{eq:ADA1}}{=} -a_k \nabla f_k^\top (\nabla \phi_k^\star (z_{k}) - x_k).
\end{equation}\vskip-5pt
\noindent Thus, we get \vskip-15pt
\begin{align*}
A_{k} U_{k+1} - A_{k-1} U_{k}
&= a_k f(x_{k}) + A_{k} (f(y_{k+1})- f(x_{k})) + A_{k-1} (f(x_{k})-f(y_{k})) \\[2pt]
&\stackrel{\mathclap{\eqref{eq:ada-conv}}}{\leq}\; a_k f(x_k) + A_{k} (f(y_{k+1})- f(x_k)) + a_k \nabla f_k^\top (\nabla \phi_k^\star (z_{k}) - x_k).
\end{align*}
This proves (a). For (b), Lemma~\ref{lemma:lower-bd}(b) and (c) lead to \vskip-15pt
\begin{align*}
A_{k} L_{k+1} - A_{k-1} L_{k}
& = a_k f(x_k) - a_k \varepsilon_{k}^\top (x_\star -x_k) - a_k g_k^\top x_k - \phi_k^\star (z_{k+1}) + \phi_{k-1}^\star (z_{k})\\[2pt]
& = a_k f(x_k)   - a_k \varepsilon_{k}^\top x_\star   - a_k \nabla f_k^\top x_k - \phi_k^\star (z_{k+1}) + \phi_{k-1}^\star (z_{k}),
\end{align*} \vskip-5pt
\noindent which completes the proof.
\end{proof}

The next theorem finishes the third step and shows the error recursion for DAM+.

\begin{theorem}\label{thm:ADA-robust}
Consider DAM+ in \eqref{eq:ADA}. Suppose that the error $\varepsilon(x)$ satisfies the growth condition in Definition \ref{asp:GC} with constants $\delta$, $\sigma^2$. Suppose $f \in \Fcal(\mcal, \Lcal)$. Given a positive weight sequence $\{a_i\}_{i=1}^k$ and $A_k = \sum_{i=1}^k a_i$, we let the convergence rate be $\rho_k = A_{k-1}/A_k$ and the potential function be $V_k = U_k - L_k$, where $U_k=f(y_k)$ and $L_k$ is defined in \eqref{eq:DA-ld2}. Then, for all $k \geq 1$, we have
\begin{equation*}
\EE V_{k+1}
\leq\rho_k \EE V_{k}
- \left(  \eta_k - \frac{(1+\delta)\Lcal \eta_k^2}{2} - \frac{(1+\delta)a_k^2}{2 A_k^2 \mcal} \right) \EE \|\nabla f_k\|^2
+ \left(  \frac{\Lcal \eta_k^2}{2} + \frac{a_k^2}{2 A_k^2 \mcal} \right) \sigma^2 .
\end{equation*}
\end{theorem}

\begin{proof}
We begin by establishing some intermediate results. We have 
\begin{equation} \label{eq:ADA-2}
a_k \nabla f_k = - a_k\varepsilon_k  + a_k g_k\stackrel{\eqref{eq:ADAb}}{=} - a_k\varepsilon_k + z_{k} - z_{k+1}
\end{equation}
and 
\begin{equation}\label{npequ:9}
\phi_{k-1}(x) \stackrel{\eqref{eq:hk}}{\leq} \phi_{k} (x), \;\; \forall x \stackrel{\text{Lemma } \ref{lemma:BD}(c)}{\Longrightarrow} \phi_{k-1}^\star (z)\geq \phi_{k}^\star (z), \;\; \forall z.
\end{equation}
Thus,
\begin{align}\label{eq:BD2}
\phi_k^\star (z_{k+1}) - \phi_{k-1}^\star (z_{k}) &- (z_{k+1}-z_{k})^\top \nabla \phi_k^\star (z_{k}) \stackrel{\eqref{npequ:9}}{\leq} \phi_k^\star (z_{k+1}) - \phi_{k}^\star (z_{k}) - (z_{k+1}-z_{k})^\top \nabla \phi_k^\star (z_{k}) \nonumber\\[2pt]
&\stackrel{\mathclap{\eqref{eq:Bregman-divergence}}}{=}\;\Delta_{\phi_k^\star} (z_{k+1},z_{k})  \stackrel{\text{Lemma }~\ref{lemma:BD}(b)}{=}  \Delta_{\phi_k}(\nabla \phi_k^\star(z_{k+1}), \nabla \phi_k^\star(z_{k})) \nonumber\\[2pt]
&\leq  \frac{m_k}{2} \|\nabla \phi_k^\star(z_{k+1  }) - \nabla \phi_k^\star(z_{k})\|^2 \stackrel{\eqref{eq:closed-form}, \eqref{eq:ADAb}}{=}  \frac{a_k^2}{2 m_k} \|g_k\|^2,
\end{align}
where $\Delta_{f}(x, y)$ is defined in \eqref{eq:Bregman-divergence} and the third equality uses the fact that $(\phi_k^\star)^\star = \phi_k$, since $\phi_k$ is convex and continuous (Fenchel\textendash Moreau theorem). Furthermore, by the smoothness of $f$, we have
\begin{multline} \label{eq:gd}
f(y_{k+1}) - f(x_k) \stackrel{\eqref{equ:class}}{\leq} \nabla f_k^\top(y_{k+1}-x_k) + \frac{\Lcal}{2}  \|y_{k+1}-x_k\|^2  \\[2pt]
=  \frac{\Lcal}{2} \|y_{k+1}-x_k\|^2 + g_k^\top(y_{k+1}-x_k) - \varepsilon_{k}^\top(y_{k+1} - x_k)
\stackrel{\eqref{eq:ADAc}}{=}\left(\frac{\Lcal \eta_k^2}{2} - \eta_k \right) \|g_k\|^2 + \eta_k \varepsilon_{k}^\top g_k.
\end{multline}
Combining the results above, we have 
\begin{align} \label{eq:ADA-pot}
 A_k &V_{k+1} - A_{k-1} V_{k} \nonumber\\[2pt]
&=  (A_k U_{k+1} - A_{k-1} U_{k}) - (A_k L_{k+1} - A_{k-1} L_{k}) \nonumber\\[2pt]
&\leftstackrel{\text{Lemma } \ref{lemma:ADAM}}{\leq} A_{k} (f(y_{k+1})- f(x_k)) + a_k \nabla f_k^\top (\nabla \phi_k^\star (z_{k}) - x_k)   +a_k \varepsilon_{k}^\top x_\star + a_k \nabla f_k^\top x_k + \phi_k^\star (z_{k+1}) - \phi_{k-1}^\star (z_{k}) \nonumber\\
&\leftstackrel{\eqref{eq:gd}}{\leq}  A_k \rbr{ -\eta_k \left(1 - \frac{\Lcal \eta_k}{2}\right) \|g_k\|^2  +  \eta_k \varepsilon_{k}^\top g_k }  + a_k \nabla f_k^\top \nabla \phi_k^\star (z_{k})  +a_k \varepsilon_{k}^\top x_\star + \phi_k^\star (z_{k+1}) - \phi_{k-1}^\star (z_{k}) \nonumber\\[2pt]
&\leftstackrel{\eqref{eq:ADA-2}}{=}   - A_k \eta_k \left(1 - \frac{\Lcal \eta_k}{2}\right) \|g_k\|^2 + A_k \eta_k \varepsilon_{k}^\top g_k + a_k \varepsilon_{k}^\top (x_\star- \nabla \phi_k^\star (z_{k})) - (z_{k+1}-z_{k})^\top \nabla \phi_k^\star (z_{k}) \nonumber\\
& \quad\quad + \phi_k^\star (z_{k+1}) - \phi_{k-1}^\star (z_{k}) \nonumber\\[2pt]
&\leftstackrel{\eqref{eq:BD2}}{\leq}  -A_k \eta_k \left(1 - \frac{\Lcal \eta_k}{2}\right)\|g_k\|^2 + A_k \eta_k \varepsilon_{k}^\top g_k + a_k \varepsilon_{k}^\top (x_\star- \nabla \phi_k^\star (z_{k})) + \frac{a_k^2}{2 m_k} \|g_k\|^2.
\end{align}
Recall that $\Gcal_k$ is the $\sigma$-algebra containing all the randomness of $\{g_i\}_{i=1}^{k-1}$, which is generated by $\{x_j, z_j, y_j\}_{j=1}^k$ for \eqref{eq:ADA}. Thus, $\EE[\varepsilon_k^\top x_i]=0$ for all $i \leq k$, $\EE[\varepsilon_k^\top g_i]=0$ for all $i<k$, and
\begin{equation*}
\EE \left[ \varepsilon_k^\top \nabla \phi_k^\star (z_k)  \middle|\Gcal_{k} \right] \stackrel{\eqref{eq:closed-form}}{=}  \EE \left[ \varepsilon_k^\top \left(-\sum_{i=1}^{k-1} a_i g_i+\mcal\sum_{i=1}^k a_i x_i\right)/m_k \middle| \Gcal_{k} \right] \stackrel{\eqref{pequ:1}}{=}0, \ \
\EE \left[ \varepsilon_k^\top g_k |\Gcal_{k} \right] \stackrel{\eqref{pequ:1}}{=} \EE \left[ \|\varepsilon_k\|^2 | \Gcal_k \right].
\end{equation*}
Multiplying \eqref{eq:ADA-pot} by $1/A_k$ and taking full expectation on both sides, we obtain
\begin{align*}
\EE \left[ V_{k+1} - \frac{A_{k-1}}{A_{k}} V_{k} \right]
& \leq - \left( \eta_k \left(1 - \frac{\Lcal \eta_k}{2}\right) - \frac{a_k^2}{2 A_k^2 \mcal} \right) \EE \| g_k\|^2 + \eta_k \EE \|\varepsilon_{k}\|^2 \\[2pt]
& \leftstackrel{\mathclap{\eqref{pequ:1}}}{=} - \left( \eta_k \left(1 - \frac{\Lcal \eta_k}{2}\right) - \frac{a_k^2}{2 A_k^2 \mcal} \right) \EE \|\nabla f_k\|^2  + \left( \eta_k - \eta_k \left(1 - \frac{\Lcal \eta_k}{2}\right) + \frac{a_k^2}{2 A_k^2 \mcal} \right) \EE \|\varepsilon_{k}\|^2 \\[2pt]
& = - \left( \eta_k \left(1 - \frac{\Lcal \eta_k}{2}\right) - \frac{a_k^2}{2 A_k^2 \mcal} \right) \EE \|\nabla f_k\|^2  + \left(  \frac{\Lcal \eta_k^2}{2}+ \frac{a_k^2}{2 A_k^2 \mcal} \right) \EE \|\varepsilon_{k}\|^2 \\[2pt]
& \leftstackrel{\mathclap{\eqref{eq:GC}}}{\leq }  -\left( \eta_k  - \frac{(1+\delta)\Lcal \eta_k^2}{2} - \frac{(1+\delta)a_k^2}{2 A_k^2 \mcal} \right) \EE \|\nabla f_k\|^2 + \left(  \frac{\Lcal \eta_k^2}{2} + \frac{a_k^2}{2 A_k^2 \mcal} \right) \sigma^2.
\end{align*}
This completes the proof.
\end{proof}

Theorem \ref{thm:ADA-robust} suggests that, under the growth condition, the multiplicative noise affects the~convergence of DAM+ through the term $\delta \{\Lcal\eta_k^2/2 + a_k^2/(2A_k^2\mcal)\}\EE \|\nabla f_k\|^2$, which can be negated by~the term $-\{\eta_k - \Lcal\eta_k^2/2 - a_k^2/(2A_k^2\mcal)\}\EE \|\nabla f_k\|^2$. Specifically, under the condition
\begin{equation*}
\delta \cbr{\frac{\Lcal\eta_k^2}{2} + \frac{a_k^2}{2A_k^2\mcal}} \leq\eta_k - \frac{\Lcal\eta_k^2}{2} - \frac{a_k^2}{2A_k^2\mcal} \Longleftrightarrow \frac{(1+\delta)a_k^2}{2A_k^2\mcal} \leq \eta_k - \frac{(1+\delta)\Lcal \eta_k^2}{2},
\end{equation*}
to let $a_k/A_k$ (i.e., the convergence rate) be as large as possible, we maximize over $\eta_k\geq 0$ on the~right hand side and obtain $\eta_k = 1/\{(1+\delta)\Lcal\}$. Plugging into the above condition, we require
\begin{equation*}
\frac{(1+\delta)a_k^2}{2A_k^2\mcal} \leq \frac{1}{2(1+\delta)\Lcal} \Longleftrightarrow \frac{a_k}{A_k} \leq \frac{1}{(1+\delta)\sqrt{\kappa}}.
\end{equation*}
Note that any $\delta\geq 0$ is allowed here. We summarize the convergence rate in the next corollary.

\begin{corollary}\label{cor:1}
Consider DAM+ in \eqref{eq:ADA} under the growth condition with constants $\delta\geq 0$~and $\sigma^2$. Under the same definitions in Theorem \ref{thm:ADA-robust}, we let $\eta_k = 1/((1+\delta)\Lcal)$ and $a_k/A_k = 1/((1+\delta)\sqrt{\kappa})$, that is, $a_k = A_{k-1}/[(1+\delta)\sqrt{\kappa}-1]$. Then,
\begin{equation} \label{eq:rate-ADA}
\EE V_{k+1} \leq \left(1-\frac{1}{(1+\delta) \sqrt{\kappa}} \right) \EE V_k + \frac{\sigma^2}{\Lcal (1+\delta)^2}, \quad\quad \forall k \geq 1.
\end{equation}
\end{corollary}

\begin{proof}
By Theorem~\ref{thm:ADA-robust} and noting that
\begin{align*}
\rho_k =  \frac{A_{k-1}}{A_k} = 1-\frac{a_k}{A_k} = 1 - \frac{1}{(1+\delta)\sqrt{\kappa}}\; \text{ and }\;
\frac{\Lcal \eta_k^2}{2} + \frac{a_k^2}{2 A_k^2 \mcal} =  \frac{1}{2(1+\delta)^2\Lcal} + \frac{1}{2(1+\delta)^2\Lcal} = \frac{1}{(1+\delta)^2\Lcal},
\end{align*}
we immediately obtain the result.	
\end{proof}

Since $f(y_k) - f(x_\star) \leq V_k$, from \eqref{eq:rate-ADA} we have that
\begin{equation*}
\frac{\mcal}{2} \EE \|y_{k} -  x_\star \|^2 
\stackrel{\eqref{equ:class}}{\leq} \EE [f(y_{k}) -  f(x_\star)] =
\Ocal \left( \left(1- \frac{1}{(1+\delta)\sqrt{\kappa}}\right)^k + \sigma^2 \right).
\end{equation*}
When $\delta =0$, DAM+ recovers the accelerated rate $1- 1/\sqrt{\kappa}$. Furthermore, for any $\delta\geq 0$, DAM+ accelerates SGD by improving the dependence on the condition number from $1/\kappa$ to $1/\sqrt{\kappa}$. This result is significantly better than the results of NAM and RMM, where a small $\delta$ is required~to~attain the accelerated rate. Compared to \cite{Diakonikolas2018Accelerated, Diakonikolas2019Approximate, Vaswani2019Fast}, our analysis has the following differences. (i) \cite{Vaswani2019Fast}~analyzed a DAM under the over-parameterized regime (i.e., $\sigma^2=0$). That method differs from DAM+ as discussed in Proposition \ref{prop:ADAM-other}, and our analysis employs the duality gap technique from \cite{Diakonikolas2018Accelerated, Diakonikolas2019Approximate} without requiring $\sigma^2=0$. (ii) \cite{Diakonikolas2018Accelerated, Diakonikolas2019Approximate} analyzed DAMs that are similar to DAM+ under \textit{deterministic setting}. We generalize their analyses to \textit{stochastic setting with the growth condition}, leading~to~quite different derivation and results. See \cite[Theorem 4.6]{Diakonikolas2019Approximate} for example. (iii) We design a unified scheme to average DAM-based iterates in \mbox{\cref{sec:5}}. We show that DAM+ enjoys the rate that can nearly match the bound \eqref{eq:lower-bound}. Such a result is missing in the literature \cite{Diakonikolas2018Accelerated, Diakonikolas2019Approximate, Vaswani2019Fast}.

\subsection{iDAM+ under the growth condition}

The implicit DAM+ is an alternative way to accelerate DAM. The iteration scheme of iDAM+ is \vskip-15pt
\begin{subequations} \label{eq:iADA}
\begin{align}
x_k &= \frac{A_{k-1}}{A_k} y_{k} + \frac{a_k}{A_k} \nabla \phi_k^\star (z_{k}), \label{eq:iADAa}\\[2pt]
z_{k+1} &= z_{k} - a_k g_k, \label{eq:iADAb}\\[2pt]
y_{k+1} &= \frac{A_{k-1}}{A_k} y_{k} + \frac{a_k}{A_k} \nabla \phi_k^\star (z_{k+1}), \label{eq:iADAc}
\end{align}
\end{subequations}\vskip-5pt
\noindent where $x_k, z_{k+1}$ in \eqref{eq:iADAa}, \eqref{eq:iADAb} are the same as in DAM+; while instead of using SGD step, $y_{k+1}$ in \eqref{eq:iADAc} is a convex combination of $y_k$ and $\nabla \phi_k^\star (z_{k+1})$. The iDAM+ scheme can be obtained by~a backward Euler discretization of the continuous dynamic induced by DAM+, and is equivalent to $\mu$AGD+ proposed by \cite{Cohen2018Acceleration}. The backward Euler discretization is also applied in \cite[(AXGD)]{Diakonikolas2018Accelerated}, \cite[(MP)]{Diakonikolas2019Approximate} to derive the accelerated extra-gradient descent (AXGD) methods. As explained after~Proposition \ref{prop:ADAM-other}, iDAM+ is also the first-order SA method.

The next theorem shows the error recursion for iDAM+.

\begin{theorem}\label{thm:iADA-robust}
Consider iDAM+ in \eqref{eq:iADA}. Suppose that the error $\varepsilon(x)$ satisfies the growth~condition in Definition \ref{asp:GC} with constants $\delta$, $\sigma^2$.	Suppose $f \in \Fcal(\mcal, \Lcal)$. Given a positive weight sequence $\{a_i\}_{i=1}^k$ and $A_k = \sum_{i=1}^k a_i$, we let the convergence rate be $\rho_k = A_{k-1}/A_k$ and the potential function~be $V_k = U_k - L_k$, where $U_k=f(y_k)$ and $L_k$ is defined in \eqref{eq:DA-ld2}. Then, for all $k \geq 1$, we have \vskip-12pt
\begin{equation*}
\EE  V_{k+1} \leq \rho_k \EE V_{k}  - \left( (1 - \delta )  \mcal  -(1+\delta)\frac{a_k^2 \Lcal}{A_k^2}\right)\frac{a_k^2}{2 A_k^2 \mcal^2}   \EE \|\nabla f_{k}\|^{2} + \left(  \mcal  + \frac{a_k^2 \Lcal}{ A_k^2}  \right) \frac{a_k^2}{2 A_k^2 \mcal^2}  \sigma^2.
\end{equation*}
\end{theorem}
\vskip-1pt
\begin{proof}
See Appendix \ref{pf:thm:iADA-robust}.
\end{proof}

Theorem~\ref{thm:iADA-robust} is similar to Theorem~\ref{thm:ADA-robust}. It shows that the multiplicative noise affects the convergence through the term $\delta\cbr{(\mcal + a_k^2\Lcal/A_k^2)a_k^2/(2A_k^2\mcal^2)}\EE \|\nabla f_{k}\|^{2}$; meanwhile, it shows that~iDAM+ can inherently tolerate this noise due to the term $\cbr{(\mcal - a_k^2\Lcal/A_k^2)a_k^2/(2A_k^2\mcal^2)}\EE \|\nabla f_{k}\|^{2}$. Since \vskip-5pt
\begin{equation*}
(1 - \delta )  \mcal  \geq (1+\delta)\frac{a_k^2 \Lcal}{A_k^2}\Longleftrightarrow \frac{a_k}{A_k} \leq \sqrt{\frac{1-\delta}{(1+\delta)\kappa}},
\end{equation*}
we can set $a_k/A_k = \sqrt{(1-\delta)/((1+\delta)\kappa)} $ to obtain the fastest rate of convergence rate $\rho_k$, requiring $\delta\in[0, 1)$ though. We summarize the convergence result in the next corollary.

\begin{corollary} \label{cor:iADA}
Consider iDAM+ in \eqref{eq:iADA} under the growth condition with constants $\delta\in[0, 1)$ and $\sigma^2$. Under the same definitions in Theorem \ref{thm:iADA-robust}, we let $a_k/A_k = \sqrt{(1-\delta)/((1+\delta)\kappa)}$, that is, $a_k = A_{k-1}/[\sqrt{(1+\delta)\kappa/(1-\delta)} - 1]$. Then \vskip-7pt
\begin{equation}\label{eq:rate-iADA}
\EE V_{k+1} \leq  \left(1- \sqrt{\frac{1-\delta}{ (1+\delta)\kappa}} \right) \EE V_{k} + \frac{1-\delta}{(1+\delta)^2\Lcal}  \sigma^2, \quad\quad \forall k\geq 1.
\end{equation}
\end{corollary}
\begin{proof}
See Appendix \ref{pf:cor:iADA}.
\end{proof}

Since $f(y_k) - f(x_\star) \leq V_k$, from  \eqref{eq:rate-iADA} we have \vskip-15pt
\begin{equation*}
\frac{\mcal}{2} \EE \|y_{k+1} -  x_\star \|^2 
\stackrel{\eqref{equ:class}}{\leq} \EE [f(y_{k+1}) -  f(x_\star)] 
=\Ocal \left( \left(1- \sqrt{\frac{1-\delta}{\kappa(1+\delta)}}\right)^k + \sigma^2 \right).
\end{equation*} \vskip-5pt
\noindent When $\delta =0$, iDAM+ recovers the accelerated convergence rate $1- 1/\sqrt{\kappa}$ and \cite[Corollary B.5]{Cohen2018Acceleration}.~We see that iDAM+ accelerates the rate of SGD for any $\delta\in[0, 1)$. Like NAM and RMM, we observe that iDAM+ can only tolerate a mild multiplicative noise, which is worse than DAM+. The recursion~of iDAM+ in \eqref{eq:rate-iADA} is the same as the recursion of NAM in Corollary \ref{cor:nam}, but for a different potential function $V_k$. Compared to the analysis of \cite{Cyrus2018Robust}, we do not impose OUBV and~bounded~domain~condition. In addition, \cite{Cyrus2018Robust} established a sub-optimal rate $\Ocal(\log k/k^{\log k} + \sigma^2\log k/k)$ when the weights~$a_k$ are specified suitably, while we strengthen this rate in \cref{sec:5} to $\Ocal(\exp(-k/\sqrt{\kappa}) + \sigma^2\log k/k)$ by proposing a unified analysis framework.

\subsection{Comparison of DAM+ and iDAM+}\label{sec:4.4}

We have studied the convergence of DAM+ and iDAM+ under the growth condition, and have shown that they converge with the rates that~are~the same as theirs under OUBV condition. Table \ref{tab:ADAM-iADAM} summarizes the results of DAM+ and iDAM+.

DAM+ and iDAM+ leverage different techniques to accelerate DAM and both achieve the accelerated rate $1-\Ocal(1/\sqrt{\kappa})$. In terms of the convergence rate, DAM+ is faster than iDAM+; while in terms of the variance, iDAM+ has smaller variance than DAM+. The rate of DAM+ holds~for any $\delta\geq 0$, while the rate of iDAM+ only holds for $\delta\in[0, 1)$. On the other hand, the tuning parameters of iDAM+ are the weights $a_k$, while the tuning parameters of DAM+ are the weights $a_k$ and stepsizes $\eta_k$, as it is a hybrid of DAM and SGD.

We notice that, to enjoy linear convergence, $a_k/A_k$ of both methods should be constant, that is, $a_k$ should be proportional to $A_{k-1} = \sum_{i=1}^{k-1}a_i$ with different multipliers. Thus, $a_k$ needs to increase exponentially (i.e. $a_k = \Ocal(e^k)$). However, if $a_k$ increases too fast, it amplifies the noise even when~the iterates are close to the minimizer $x_\star$, which counteracts the benefit of averaging the past gradients. Indeed, under bounded domain assumption, \cite{Cohen2018Acceleration} also showed that a constant ratio $a_k/A_k$ results in a $\Ocal((1-1/\sqrt{\kappa})^k)$ bias term  (i.e. linear rate) and a constant variance term; while if $a_k =\Ocal(k^p)$,~then one has $\Ocal(1/k^{p+1})$ for bias (i.e. sublinear rate) and $\Ocal(\sigma^2/k)$ for variance. In other words, slowing down the increase of $a_k$ decreases the variance term and even makes the variance converge to zero, but decelerates the convergence rate as well.

Combining Tables \ref{tab:SGD-NAM-RMM} and \ref{tab:ADAM-iADAM} together, we see clearly that all the studied accelerated methods,~NAM, RMM, DAM+, and iDAM+, can accelerate SGD under the growth condition, as how they behave under OUBV condition. However, their ability to tolerate the multiplicative noise is different. The accelerated rate of DAM+ holds for any $\delta\geq 0$, while the accelerated rates of the other three methods hold for a moderate $\delta$ (e.g., $\delta<1$). When $\delta=0$ (i.e., OUBV condition), RMM enjoys the best bias term with the worst variance term. Furthermore, iDAM+ and NAM enjoy similar~bias~and variance terms; DAM+ enjoys the same variance term as SGD but uniformly improves the bias~term. Our findings are greatly aligned with the existing results on least-squares regressions: NAM and HB fail to accelerate SGD \citep{Kidambi2018Insufficiency, Liu2020Accelerating} while DAM+ can accelerate SGD \citep{Jain2018Accelerating}. Additionally, our results provably suggest that NAM can tolerate a \textit{mild} multiplicative noise ($\delta<1$), which is somewhat surprising to the existing understanding where people tend to believe NAM is always fragile for~the~multiplicative  noise, based on its performance on LSR \citep{Liu2020Accelerating}. Our results do not contradict the observations in \cite{Liu2020Accelerating} since LSR has a large multiplicative noise; see \cref{sec:1.1}.

\begin{table}
\TABLE
{Summary of convergence rates for DAM+ and iDAM+. \label{tab:ADAM-iADAM}}
{\begin{tabular}{cccc}
		\hline
		\up \down  Algorithm & Bias term ($\EE V_k$) & Variance term($\sigma^2$)   & Requirement for $\delta$  \vspace{2pt} \\
		\hline
		\up \down  DAM+ & $1- 1/[(1+\delta)\sqrt{\kappa}]$ & $1/[(1+\delta)^2\Lcal]$  & $0\leq \delta$     \vspace{2pt}   \\
		\hline
		\up \down  iDAM+ & $1- \sqrt{(1-\delta)/[(1+\delta)\kappa]}$ & $(1-\delta)/[(1+\delta)^2\Lcal]$  & $0\leq \delta <1$    \vspace{2pt}      \\
		\hline
\end{tabular}\vspace{-5pt}}
{}
\end{table}

\section{Unified Near-optimal Analysis Framework} \label{sec:5}

In this section, we propose a~generic~scheme for averaging iterates and diminishing parameters for DAM+ and iDAM+. With this scheme,~we show that their convergence rates \textit{nearly} match the theoretical lower bound of first-order SA methods under OUBV (or compact domain, which implies OUBV). For convenience, we rewrite \eqref{eq:lower-bound}~as \vskip-8pt
\begin{equation} \label{eq:SA-LB}
f(x_k)-f(x_\star)=\Omega \rbr{ \exp \rbr{- \frac{k}{\sqrt{\kappa}}} + \frac{\sigma^2}{k}}.
\end{equation} \vskip-2pt
\noindent Since an oracle that satisfies OUBV also satisfies the growth condition, this lower bound is also~valid under the growth condition.

Let us consider the following recursion \vskip-8pt
\begin{equation} \label{eq:recursion}
e_{k+1} \leq (1- \gamma_1 r_k) e_k -  \gamma_2 r_k s_k + \gamma_3 r_k^2
\end{equation} \vskip-2pt
\noindent for a positive sequence $\{e_k, s_k, r_k\}_{k>0}$ and non-negative constants $\gamma_1, \gamma_2, \gamma_3$ such that $r_k \leq 1/\gamma_4$ and $\gamma_1/\gamma_4\leq1$ for some $\gamma_4$. The recursions established for DAM+ and iDAM+ in \cref{sec:4} can be cast as \eqref{eq:recursion}. For example, we consider the recursion of DAM+ in Theorem \ref{thm:ADA-robust}. Let $a_k^2/A_k^2 = \mcal \Lcal \eta_k^2$ and restrict $\eta_k \leq 1/(2\Lcal(1+\delta))$, then we obtain \vskip-8pt
\begin{equation*}
\EE V_{k+1} \leq (1- \sqrt{\mcal \Lcal} \eta_k)\EE V_k - \frac{\eta_k}{2} \EE \|\nabla f(x_k)\|^2 + \Lcal \sigma^2\eta_k^2 .
\end{equation*}\vskip-2pt
\noindent Comparing the above display with \eqref{eq:recursion}, we have the correspondence \vskip-8pt
\begin{equation*}
(e_k, s_k, r_k) = (\EE V_k, \EE \| \nabla f(x_k) \|^2, \eta_k), \ \  (\gamma_1,\gamma_2, \gamma_3, \gamma_4)= \left( \sqrt{\mcal \Lcal}, 1/2,  \Lcal \sigma^2, 2 \Lcal (1+\delta)\right).
\end{equation*} \vskip-2pt
\noindent If we choose $r_k$ to be as large as possible, that is, $r_k=1/\gamma_4$, then solving \eqref{eq:recursion} leads to \vskip-8pt
\begin{equation} \label{eq:trivial-bd}
e_{k+1}\leq \left(1- \frac{\gamma_1}{\gamma_4} \right)^{k} e_1 + \frac{\gamma_3}{\gamma_4 \gamma_1},
\end{equation} \vskip-2pt
\noindent which shows that the constant rate $r_k$ implies a linear convergence to a neighborhood of the~optimum. The term $\gamma_3/(\gamma_4 \gamma_1)$ corresponds to the variance term of the algorithm that does not vanish when $k$ increases. A scheme to diminish the algorithm's parameters is required to eliminate the~constant variance term and achieve the bound in \eqref{eq:SA-LB}.

We design the following scheme to diminish $r_k$ and to average the first $k$ iterates with weights $\{w_{k, t}\}_{t\leq k}$: given a tuning parameter $\alpha>0$, \vskip-8pt
\begin{equation} \label{eq:suffix-averaging}
r_k =\begin{cases}
1/\gamma_4 & \text{if } \alpha k<\gamma_4 \\
1/(\gamma_1 k) & \text{otherwise}
\end{cases}
, \quad\quad\quad
w_{k, t} =
\begin{cases}
0 &\text{if } \alpha k<\gamma_4 \text{ and } t<k\\
1 & \text{if } \alpha k<\gamma_4 \text{ and } t=k\\
0 & \text{if } \alpha k \geq \gamma_4 \text{ and } \alpha t < \gamma_4\\
1 & \text{if } \alpha k \geq \gamma_4 \text{ and } \alpha t\geq \gamma_4
\end{cases}.
\end{equation}
In what follows, we assume $\gamma_4/\alpha$ is an integer for simplifying our presentation, but $\alpha=\gamma_1$ is a more natural choice. If $\gamma_3=0$, that is $\sigma^2=0$, we can simply choose $\alpha=0$ and \eqref{eq:suffix-averaging} recovers the analysis under the over-parameterized regime. The setup \eqref{eq:suffix-averaging} consists of two phases. In the first phase, we use a constant stepsize to reduce the bias term exponentially, and assign all weights on the current iterate, i.e. $w_{k,t}=0$ for $t<k$ and $w_{k,t}=1$ for $t=k$. Thus, the iteration sequence approaches to a neighborhood of the optimum. In the second phase, we diminish the convergence rate $1-r_k$ and average the last several iterates to reduce the variance term.

The scheme \eqref{eq:suffix-averaging} is similar to the $\alpha$-suffix averaging \citep{Rakhlin2012Making}. Given the total number of iterations~$K$ and a constant $0<\alpha<1$, the $\alpha$-suffix averaging discards the first $\alpha K$ iterates and averages the last $(1-\alpha) K$ iterates. The major difference of \eqref{eq:suffix-averaging} to the $\alpha$-suffix averaging is that \eqref{eq:suffix-averaging} uses constant parameters for finite $k$ iterations until $k=\gamma_4/\alpha$, which does not depend on the total number of iterations $K$. In contrast, the $\alpha$-suffix averaging employs constant parameters for $k$ iterations until $ k= \alpha K$. In other words, when $K$ varies, the $\alpha$-suffix averaging results in different schemes to decrease parameters in different stages. Our scheme is consistent regardless of $K$.

The following proposition shows how to nearly achieve \eqref{eq:SA-LB} with the scheme \eqref{eq:suffix-averaging}.

\begin{proposition}\label{prop:suffix-averaging}
Let $W_k = \sum_{t=1}^k w_{k, t}$. Given $r_k$ and weights $w_{k,t}$ defined in \eqref{eq:suffix-averaging}, we have \vskip-10pt
\begin{equation} \label{eq:sa-bd}
\frac{1}{W_{k}} \sum_{t=1}^{k} w_{k,t} s_{t} = \Ocal \left( \exp\left[ -\frac{\gamma_1 k}{\gamma_4}\right] + \frac{\gamma_3 \log k}{k} \right), \quad \forall k \geq 1.
\end{equation}
\end{proposition}

\vskip-2pt

\begin{proof}
The recursion \eqref{eq:recursion} implies that \vskip-10pt
\begin{equation} \label{eq:rec-1}
e_{k+1} \leq (1- \gamma_1 r_k) e_k -  \gamma_2 r_k s_k + \gamma_3 r_k^2 \leq (1- \gamma_1 r_k) e_k + \gamma_3 r_k^2.
\end{equation} \vskip-2pt
\noindent For $k$ such that $\alpha k <\gamma_4$, by \eqref{eq:suffix-averaging} we have $\sum_{t=1}^{k} w_{k, t} s_{t}/W_k = s_k$ and $r_k = 1/\gamma_4$. Thus, \vskip-15pt
\begin{align} \label{eq:phase1}
\frac{\gamma_2}{W_k \gamma_4} \sum_{t=1}^{k} w_{k, t} s_{t}
&= \gamma_2 r_k s_k \stackrel{\eqref{eq:recursion}}{\leq} (1- \gamma_1 r_k) e_k + \gamma_3 r_k^2 \;\;
\stackrel{\mathclap{\eqref{eq:rec-1}}}{\leq} \;\; \rbr{1 - \frac{\gamma_1}{\gamma_4}}^k e_1 + \frac{\gamma_3}{\gamma_4^2} \sum_{t=0}^{\infty}\rbr{1 - \frac{\gamma_1}{\gamma_4}}^t \nonumber\\[2pt]
&= \left(1- \frac{\gamma_1}{\gamma_4}\right)^k e_1 +  \frac{\gamma_3}{\gamma_1 \gamma_4}
\leq  e_1 \exp\left[ -\frac{\gamma_1 k}{\gamma_4}\right] + \frac{\gamma_3}{\gamma_1 \gamma_4} \leq  e_1 \exp\left[ -\frac{\gamma_1 k}{\gamma_4}\right] + \frac{\gamma_3}{\gamma_1 \alpha k}.
\end{align} \vskip-4pt
\noindent Here, the fifth inequality is due to the fact that $0\leq \gamma_1/\gamma_4\leq 1$, coming from the condition of the recursion, and the fact that $(1-x)^k \leq \exp(-kx)$ for any $x\in[0,1]$; the sixth inequality is due to $\alpha k <\gamma_4$. We see that \eqref{eq:phase1} is consistent with the statement.

For $k$ such that $\alpha k \geq \gamma_4$, we consider $t\leq k$ such that $\alpha t \geq \gamma_4$. From \eqref{eq:suffix-averaging} we have that $w_{k,t}=1$. Arranging \eqref{eq:recursion} implies that \vskip-18pt
\begin{align} \label{eq:wktst}
\gamma_2 w_{k,t} s_{t}
&\leq  \;\frac{w_{k,t} \left(1-\gamma_1 r_{t}\right) e_{t}}{r_{t}}-\frac{ w_{k,t} e_{t+1}}{r_{t}}+ \gamma_3 w_{k,t} r_{t}
\stackrel{\eqref{eq:suffix-averaging}}{=} \frac{ e_{t}}{r_{t}} - \gamma_1 e_{t} -\frac{ e_{t+1}}{r_{t}}+\gamma_3 r_{t} \nonumber \\[2pt]
&\leftstackrel{\mathclap{\eqref{eq:suffix-averaging}}}{=}\;\; \gamma_1 t  e_{t}- \gamma_1 e_{t} -  \frac{e_{t+1}}{r_t} +\frac{\gamma_3 }{\gamma_1 t} = \gamma_1(t-1)  e_{t}- \gamma_1t e_{t+1} +\frac{\gamma_3 }{\gamma_1 t}.
\end{align} \vskip-3pt
\noindent Thus, if $\alpha k \geq \gamma_4$, we note that $\gamma_4/\alpha$ is a positive integer and $w_{k, t} = 0$ for $\alpha t<\gamma_4$, and further have \vskip-17pt
\begin{align} \label{eq:second-iter}
\gamma_2 \sum_{t=1}^{k} w_{k,t} s_{t}
&=\; \gamma_2 \sum_{t=1}^{\gamma_4/\alpha-1} w_{k,t} s_{t} + \gamma_2 \sum_{t=\gamma_4/\alpha}^{k} w_{k,t} s_{t} =  \gamma_2 \sum_{t=\gamma_4/\alpha}^{k} w_{k,t} s_{t} \nonumber\\[2pt]
&\leftstackrel{\mathclap{\eqref{eq:wktst}}}{\leq}\; \sum_{t=\gamma_4/\alpha}^k \left( \gamma_1(t-1)e_t - \gamma_1te_{t+1}  +\frac{\gamma_3 }{\gamma_1 t}  \right) = \gamma_1(\gamma_4/\alpha-1)e_{\gamma_4/\alpha} - \frac{e_{k+1}}{r_{k}} + \frac{\gamma_3}{\gamma_1}\sum_{k=\gamma_4/\alpha}^k \frac{1 }{t} \nonumber \\[2pt]
&\leq  \frac{\gamma_1 \gamma_4 e_{\gamma_4/\alpha}}{\alpha}  + \frac{\gamma_3(\log (k) + 1)}{\gamma_1},
\end{align}\vskip-4pt
\noindent where the last inequality uses $\sum_{t=1}^{ k} 1/t \leq \log\left(  k \right) + 1$. Furthermore, by  \eqref{eq:trivial-bd}, \vskip-5pt
\begin{equation} \label{eq:first-iter}
e_{\gamma_4/\alpha} \leq \rbr{1 - \frac{\gamma_1}{\gamma_4}}^{\gamma_4/\alpha-1}e_1 + \frac{\gamma_3}{\gamma_4 \gamma_1}\leq e_1 \exp\left[ -\frac{\gamma_1 (\gamma_4/\alpha-1)}{\gamma_4} \right] + \frac{\gamma_3}{\gamma_1 \gamma_4}.
\end{equation}\vskip-2pt
\noindent Since $W_k = \sum_{t=1}^{k} w_{k,t} = k -  \gamma_4/\alpha+1$, combining \eqref{eq:second-iter} and \eqref{eq:first-iter} leads to \vskip-16pt
\begin{align} \label{eq:phase2}
& \frac{\gamma_2}{W_{k}} \sum_{t=1}^{k} w_{k,t} s_{t} \;
\leftstackrel{\mathclap{\eqref{eq:second-iter}}}{\leq}\; \frac{\gamma_1 \gamma_4 e_{\gamma_4/\alpha}}{(k -  \gamma_4/\alpha+1)\alpha} + \frac{\gamma_3(\log k + 1)}{(k -  \gamma_4/\alpha+1)\gamma_1} \nonumber\\[2pt]
&\leftstackrel{\mathclap{\eqref{eq:first-iter}}}{ \leq}\; \frac{\gamma_1 \gamma_4}{(k -  \gamma_4/\alpha+1)\alpha} \cbr{e_1 \exp\left[ -\frac{\gamma_1 (\gamma_4/\alpha-1)}{\gamma_4} \right] + \frac{\gamma_3}{\gamma_1 \gamma_4} } + \frac{\gamma_3(\log k + 1)}{(k -  \gamma_4/\alpha+1)\gamma_1} = \Ocal \left( \frac{\gamma_3 \log k}{k}\right).
\end{align}\vskip -2pt
\noindent The proof follows by combining \eqref{eq:phase1} and \eqref{eq:phase2}.
\end{proof}

Proposition \ref{prop:suffix-averaging} is different from the existing results. Given a fixed total number of iterations~$K$, \cite{Rakhlin2012Making, Jain2018Accelerating} used tail averaging and showed the convergence for all $k\leq K$. As mentioned earlier, a finite $K$ is required to make their scheme and results well defined. In contrast, the result in Proposition~\ref{prop:suffix-averaging}  holds for any $k\geq 1$, even if $k$ goes to infinity. In addition, \cite[Lemma 2]{Stich2019Unified} showed that there exists a constant stepsize $r$ such that 
\begin{equation*}
\frac{1}{W_{K}} \sum_{k=1}^{K} w_{k} s_{k} = \Ocal \left( \exp\left[ -\frac{\gamma_1 K}{\gamma_4}\right] + \frac{\gamma_3 \log K}{K} \right),
\end{equation*}
where $w_k=(1-\gamma_1 r)^{-(k+1)}$ and $W_k = \sum_{t=1}^{k} w_{t}$. However, their scalar $\gamma_4$ relies on $K$ and their~argument fails for varying $K$. Finally, we note that if $r_k$ is chosen more carefully, one may remove $\log k$ term in \eqref{eq:sa-bd}. See \cite{Grimmer2019Convergence} for example. However, we leave it for future work.

Applying Proposition \ref{prop:suffix-averaging}, we are able to show the convergence of DAM+ and iDAM+ in the~next two theorems. In particular, we show that DAM+ and iDAM+ achieve the near-optimal convergence rate using the scheme \eqref{eq:suffix-averaging}, although they can tolerate different noise levels of $\delta$.

\begin{theorem} \label{thm:ADA-suffix-averaging}
Consider DAM+ in \eqref{eq:ADA} under the growth condition in Definition \ref{asp:GC} with~$\delta\geq 0$. Let $r_k = \eta_k$ and $w_{k, t}$ be specified by \eqref{eq:suffix-averaging} with $\gamma_1 = \sqrt{\mcal\Lcal}$, $\gamma_4 = 2\Lcal(1+\delta)$ and any $\alpha>0$. Let $a_k/A_k = \sqrt{\mcal \Lcal}\eta_k$. Then, for any $k\geq 1$, it holds that \vskip-4pt
\begin{equation} \label{eq:ADAM-lower-bd}
\EE [f(\bar{x}_k) - f(x_\star)]
= \Ocal \left( \exp\left[ - \frac{k }{2(1+\delta)\sqrt{\kappa}}\right] + \frac{\log\left(  k \right)}{ k}  \sigma^2 \right),
\end{equation}
where $\bar{x}_k = \sum_{t=1}^k w_{k,t} x_t /W_k$ and $W_k = \sum_{t=1}^{k} w_{k,t}$.
\end{theorem}

\begin{proof}
Plugging $a_k^2/A_k^2 = \mcal \Lcal \eta_k^2$ into Theorem \ref{thm:ADA-robust} leads to \vskip-6pt
\begin{equation*}
\EE V_{k+1}
\leq  (1- \sqrt{\mcal \Lcal} \eta_k) \EE V_{k}
- \left(  \eta_k - (1+\delta)\Lcal \eta_k^2\right) \EE \|\nabla f_k\|^2
+ \Lcal \eta_k^2 \sigma^2.
\end{equation*}
Since
\begin{equation*}
\frac{\eta_k}{2}\geq  (1+\delta)\Lcal \eta_k^2 \Longleftrightarrow \eta_k \leq \frac{1}{ 2(1+\delta) \Lcal} = \frac{1}{\gamma_4},
\end{equation*}
under the above condition we further have \vskip-6pt
\begin{equation*}
\EE V_{k+1} \leq (1- \sqrt{\mcal \Lcal} \eta_k)\EE V_k - \frac{\eta_k}{2} \EE \|\nabla f_k\|^2 + \Lcal \eta_k^2 \sigma^2.
\end{equation*}
Comparing with the recursion \eqref{eq:recursion}, we define \vskip-6pt
\begin{equation} \label{eq:ADA-rec-3}
(e_k, s_k) = (\EE V_k, \EE \| \nabla f(x_k) \|^2), \; \; \;\; (\gamma_1,\gamma_2, \gamma_3)= ( \sqrt{\mcal \Lcal}, 1/2,  \Lcal \sigma^2) .
\end{equation}
Then, by the convexity of $f$ and $\gamma_4=2\Lcal(1+\delta)$, we get \vskip-15pt
\begin{align*}
&\EE [f(\bar{x}_k) - f(x_\star)]
\leq \frac{1}{W_k} \sum_{t=1}^k w_{k,t} \EE \left[ f(x_t) -f(x_\star) \right]
\stackrel{\text{Lemma}~\ref{lemma:convex}(a)}{\leq} \frac{1}{2 \mcal W_k} \sum_{t=1}^k w_{k,t} \EE \|\nabla f(x_t)\|^2 
\\[2pt]
&\stackrel{\eqref{eq:ADA-rec-3}}{=}  \Ocal \left( \frac{1}{ W_k} \sum_{t=1}^{k} w_{k, t} s_t \right)
\stackrel{\eqref{eq:sa-bd}}{=}  \Ocal \left( \exp\left[ -\frac{\gamma_1 k}{\gamma_4}\right] + \frac{\gamma_3 \log k}{k} \right)
\stackrel{\eqref{eq:ADA-rec-3}}{=}  \Ocal \left( \exp\left[ - \frac{k }{2(1+\delta)\sqrt{\kappa}}\right] + \frac{\log\left(  k \right)}{ k}  \sigma^2 \right),
\end{align*}
which completes the proof.
\end{proof}

\vspace{-10pt}
\begin{theorem} \label{thm:iADA-suffix-averaging}
Consider iDAM+ in \eqref{eq:iADA} under the growth condition in \mbox{Definition}~\ref{asp:GC}~with~$\delta\in[0,1)$. Let $r_k = a_k/A_k$ and $w_{k, t}$ be specified by \eqref{eq:suffix-averaging} with $\gamma_1 = 1/2$, $\gamma_4 = \sqrt{(1+\delta)\kappa/(1-\delta)}$ and any $\alpha>0$. Then, for any $k\geq 1$, it holds that
\begin{equation} \label{eq:iADAM-lower-bd}
\EE [f(\bar{y}_k) - f(x_\star)]
= \Ocal \left( \exp\left[ - \frac{k }{2} \sqrt{\frac{1-\delta}{(1+\delta)\kappa} }\right] + \frac{\log\left(  k \right)}{ k}  \sigma^2 \right),
\end{equation}
where $\bar{y}_k = \sum_{t=1}^k w_{k,t} y_t /W_k$ and $W_k = \sum_{t=1}^{k} w_{k,t}$.
\end{theorem}

\begin{proof}
Letting $a_k/A_k=r_k$, Theorem \ref{thm:iADA-robust} leads to \vskip-10pt
\begin{equation*}
\EE  V_{k+1} \leq (1-r_k) \EE V_{k} -\frac{r_k^2}{2 \mcal^2}  \left( (1 - \delta )  \mcal  -(1+\delta)r_k^2 \Lcal\right) \EE \|\nabla f_{k}\|^{2} +\frac{r_k^2}{2 \mcal^2}\left(  \mcal  + r_k^2 \Lcal \right)   \sigma^2.
\end{equation*} \vskip-3pt
\noindent Since \vskip-15pt
\begin{equation*}
(1 - \delta )  \mcal\geq (1+\delta)r_k^2 \Lcal  \Longleftrightarrow r_k \leq \sqrt{\frac{1-\delta}{(1+\delta)\kappa} } = \frac{1}{\gamma_4},
\end{equation*}\vskip-5pt
\noindent under the above condition we further have \vskip-10pt
\begin{equation*}
\EE  V_{k+1}\leq (1-r_k)  \EE V_{k}+\frac{r_k^2}{2 \mcal^2}\left(  \mcal  + r_k^2 \Lcal \right)   \sigma^2 \leq (1-r_k)  \EE V_{k}+\frac{r_k^2}{(1+\delta)\mcal}  \sigma^2.
\end{equation*} \vskip-5pt
\noindent Comparing with the recursion \eqref{eq:recursion}, we define \vskip-10pt
\begin{equation} \label{eq:iADA-rec}
(e_k, s_k) = (\EE V_k, \EE V_k), \; \; \;\; (\gamma_1,\gamma_2, \gamma_3)= ( 1/2, 1/2,  \sigma^2/[(1+\delta)\mcal]) .
\end{equation}\vskip-3pt
\noindent Then, by the convexity of $f$, we get \vskip-15pt
\begin{align*}
&\EE [f(\bar{y}_k) - f(x_\star)]
\leq  \frac{1}{W_k} \sum_{t=1}^k w_{k,t} \EE \left[ f(y_t) -f(x_\star) \right] \leq \frac{1}{W_k} \sum_{t=1}^k w_{k,t} \EE V_t \\[2pt]
&\stackrel{\eqref{eq:iADA-rec}}{=}   \frac{1}{ W_k} \sum_{t=1}^{k} w_{k, t} s_t  \stackrel{\eqref{eq:sa-bd}}{=}  \Ocal \left( \exp\left[ -\frac{\gamma_1 k}{\gamma_4}\right] + \frac{\gamma_3 \log k}{k} \right)
\stackrel{\eqref{eq:iADA-rec}}{=}  \Ocal \left( \exp\left[ - \frac{k }{2} \sqrt{\frac{1-\delta}{(1+\delta)\kappa} }\right] + \frac{\log\left(  k \right)}{ k}  \sigma^2 \right),
\end{align*}
where the second inequality uses $f(y_t) - f(x_\star) \leq V_t$. This completes the proof.
\end{proof}

We end this section by discussing Theorems \ref{thm:ADA-suffix-averaging} and \ref{thm:iADA-suffix-averaging}. The theorems suggest that DAM+ and iDAM+ can improve the dependence on the condition number in the convergence rate from $1/\kappa$ to $1/\sqrt{\kappa}$, under the growth condition. However, the ability to tolerate the noise level $\delta$ is different. There is existing literature that shows similar results to Theorems \ref{thm:ADA-suffix-averaging} and \ref{thm:iADA-suffix-averaging} for accelerated~SGD methods under stronger conditions. \cite{Cohen2018Acceleration} assumed OUBV and bounded domain conditions for {\small iDAM+}, and showed the following convergence rate \vskip-10pt
\begin{equation*}
\Ocal \left(\frac{\log (k)}{k^{\log (k)}} \cdot \frac{(\Lcal-\mcal)\left\|x_{\star}-x_{0}\right\|^{2}}{2}+\frac{\log (k)}{k} \cdot \frac{\sigma^{2}}{\mcal}\right).
\end{equation*} \vskip-3pt
\noindent Our analysis improves the above bias term to $\exp (-k)$ in a more general setting. \cite{Ghadimi2013Optimal} introduced an accelerated stochastic approximation algorithm (AC-SA). Using a multistage trick that changes the parameters of AC-SA with respect to the iteration $k$, they proved that AC-SA achieved \eqref{eq:sa-bd} under OUBV. In contrast, our theorem provides a unified framework where the error recursion~of~any algorithm satisfying \eqref{eq:recursion} can meet \eqref{eq:sa-bd} without OUBV. Our scheme in \eqref{eq:suffix-averaging} is also more straightforward than the multistage trick.

\section{Numerical Experiments}\label{sec:ns}
	
\hskip -3pt In this section, we conduct numerical experiments to~\mbox{validate}~our theoretical findings for the four considered accelerated methods: NAM, RMM, DAM+,~and~iDAM+. We optimize a simple quadratic loss, consider a stochastic oracle that satisfies the growth condition but not OUBV, and also implement SGD for the sake of a complete comparison. We set the parameters of all the methods as suggested in the corresponding theorems, so that we can investigate the exact relationship between convergence rate, condition number, and multiplicative noise. From the experiments, we observe that the numerical results are consistent with the conclusions drawn from Corollaries \ref{cor:nam}, \ref{cor:rmm}, \ref{cor:1}, and \ref{cor:iADA}. Furthermore, the numerical results reveal that our theorems not only present upper bounds of the error recursions, but also accurately and sharply characterize the true behavior of each method (at least on a quadratic objective). Our code is publicly available at \url{https://github.com/youlinchen/Convergence-Analysis-of-AcceleratedSGD}.

\subsection{Experimental setups}

We let $x\in \RR^2$ and consider the problem of optimizing a quadratic loss $f(x) = x^\top\Wb x/2$, where $\Wb = \diag(\kappa, 1)$ for some varying condition number $\kappa\geq 1$ specified~later. We also consider a stochastic oracle of noisy gradients given by $g(x) = \Zb\nabla f(x)$ with $\Zb = \diag(z_1, z_2)$. Here, $z_1$ and $z_2$ are draw independently from a normal distribution with mean $1$ and variance $\delta/2$.

\newpage
\noindent With the above setup, we note that \vskip-6pt
\begin{equation*}
\EE\|\epsilon(x)\|^2 = \EE\|g(x) - \nabla f(x)\|^2 = \EE\|(\Zb - \Ib)\nabla f(x)\|^2 = \delta\|\nabla f(x)\|^2. 
\end{equation*}
Thus, our stochastic oracle satisfies the growth condition (with $\sigma^2=0$) but does not satisfy OUBV.

We aim to examine the relationship between convergence rate, condition number $\kappa$, and multiplicative noise level $\delta$. We take DAM+ as an example. By Corollary \ref{cor:1}, we know that by setting $\eta_k = 1/((1+\delta)\kappa)$ and $a_k = A_{k-1}/((1+\delta)\sqrt{\kappa}-1)$, DAM+ should exhibit the following behavior~(note that $f(x_\star) = 0$): \vskip-6pt
\begin{equation}\label{equ:1}
\EE[f(y_k)] \leq \EE[f(y_0)]\rbr{1 -\frac{1}{(1+\delta)\sqrt{\kappa}} }^k \leq \EE[f(y_0)]\cdot \exp\rbr{- \frac{k}{(1+\delta)\sqrt{\kappa}}}.
\end{equation}
To approximate the expectation on the left-hand side, we perform 1000 independent runs with~random initializations for each setup of $\delta$ and $\kappa$, and let $\EE[f(y_k)]\approx \sum_{j=1}^{1000}f(y_k^{(j)})/1000$. Here, $y_k^{(j)}$ is the iterate $y_k$ in the $j$-th run. Furthermore, we let $K$ be the least iteration index that triggers the following stopping criterion:\vskip-6pt
\begin{equation*}
\frac{1}{1000}\sum_{j=1}^{1000}f(y_K^{(j)}) \leq 10^{-8}\cdot \frac{1}{1000}\sum_{j=1}^{1000}f(y_0^{(j)}).
\end{equation*} 
The definition of $K$ is motivated by \eqref{equ:1}, based on which we expect to see a linear convergence of the initial error. In other words, we expect to see that for DAM+, $K/((1+\delta)\sqrt{\kappa})$ is a constant.~We conduct experiments under the same setup for the other four methods. For all methods, we should have the relationship $K\times \text{convergence rate}= \text{constant}$. Thus, we fix either $\delta$ or $\kappa$ and plot $K$ versus the other quantity. The plots enable us to explore if the accelerated methods exhibit accelerated rates and if the established accelerated rates are sharp or not.

\subsection{Experimental results}

\begin{figure}[h]
\centering
\includegraphics[width=0.325\textwidth]{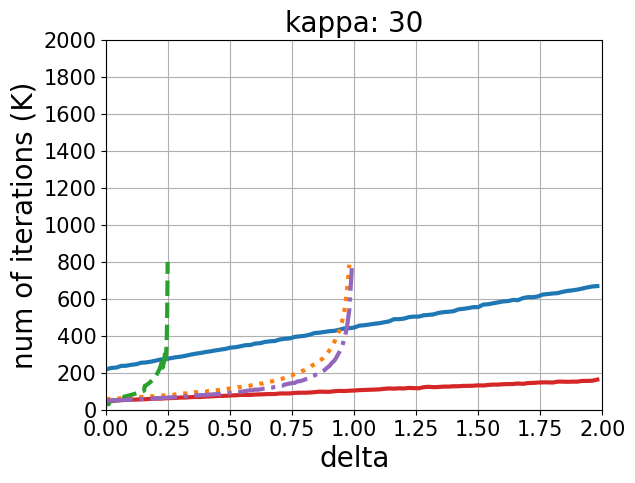}
\includegraphics[width=0.325\textwidth]{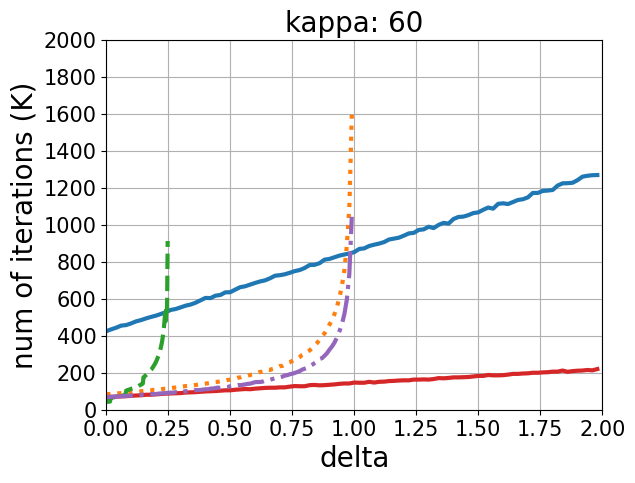}
\includegraphics[width=0.325\textwidth]{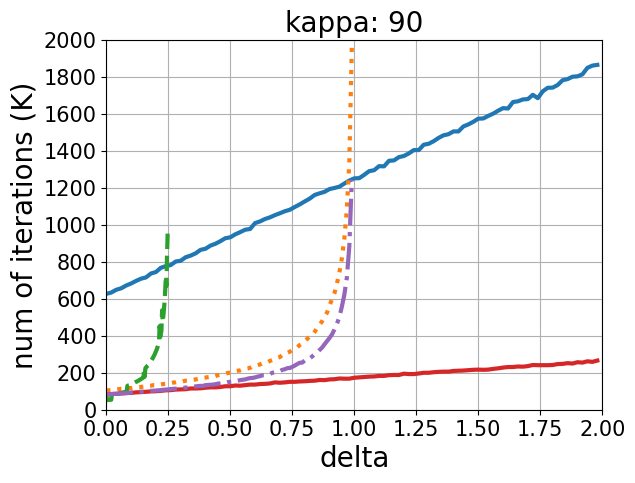}
\includegraphics[width=0.6\textwidth]{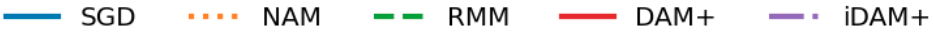}
\caption{The relationship between $K$ and $\delta$ given a fixed $\kappa$. We observe that $K \asymp_\kappa \sqrt{(1+\delta)/(1-\delta)}$ for NAM and iADAM+, $K\asymp_\kappa 1+\delta$ for SGD and DAM+, and $K\asymp_\kappa 1/(1-2\sqrt{\delta})$ for RMM. These relationships are consistent with our theorems.} \label{fig:1}	
\end{figure}

\begin{figure}[h]
\centering
\includegraphics[width=0.325\textwidth]{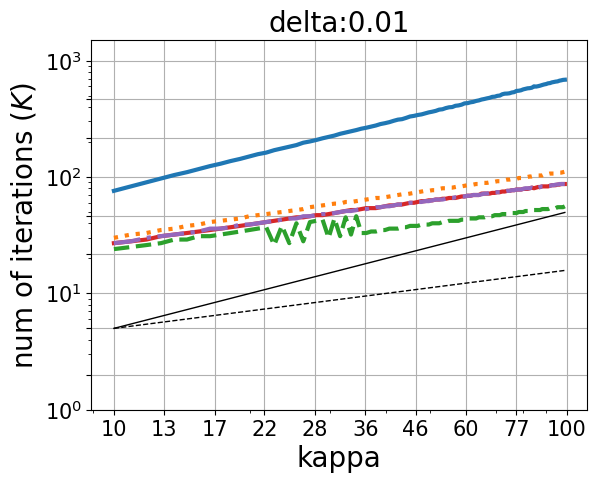}
\includegraphics[width=0.325\textwidth]{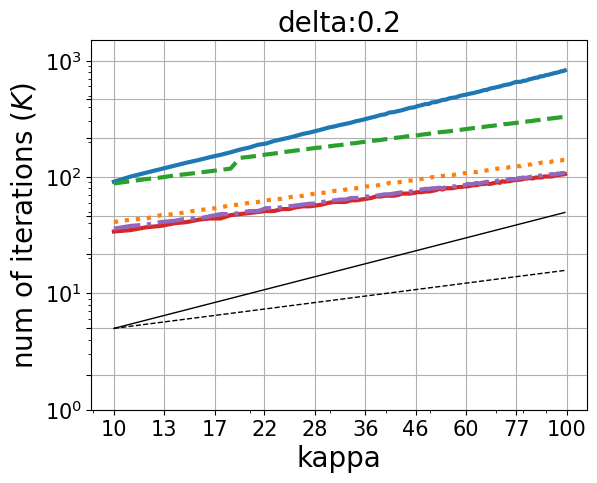}
\includegraphics[width=0.325\textwidth]{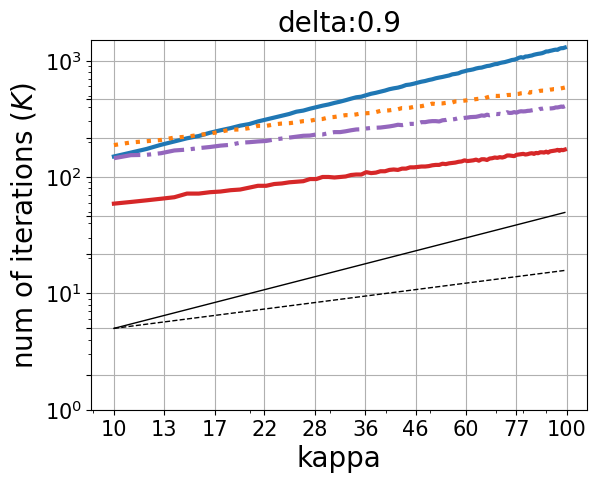}
\includegraphics[width=0.6\textwidth]{Figures/legend1} \ \ 
\includegraphics[width=0.3\textwidth]{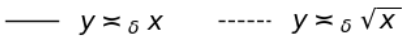}
\caption{The relationship between $K$ and $\kappa$ given a fixed $\delta$ (in log scale). The solid and dashed black lines are~two guided lines to reveal the dependence of $K$ on $\kappa$. We observe that $K\asymp_\delta \kappa$ for SGD and $K\asymp_\delta\sqrt{\kappa}$ for the other four accelerated methods. These relationships are consistent with our theorems.} \label{fig:2}	
\end{figure}

Our results are summarized in Figures \ref{fig:1} and \ref{fig:2}. For simplicity, let us denote $a\asymp_\kappa b$ if $a/b$ is a constant depending on $\kappa$. Then, by \eqref{pequ:SGD}, Corollaries \ref{cor:nam}, \ref{cor:rmm}, \ref{cor:1}, and \ref{cor:iADA}, we expect to observe $K \asymp_\kappa \sqrt{(1+\delta)/(1-\delta)}$ for NAM and iDAM+, $K\asymp_\kappa 1+\delta$ for DAM+~and SGD, and $K\asymp_\kappa 1/(1-2\sqrt{\delta})$ for RMM. These relationships are precisely discovered in Figure \ref{fig:1}. Furthermore, from Figure \ref{fig:1}, we observe that the convergence rates of all methods become slower~as $\kappa$ increases, while the relationships between $K$ and $\delta$ remain unchanged. It is also worth noting~that although NAM and iDAM+ have the same tolerance of $\delta$ and similar convergence rates based on our theorems, iDAM+ slightly outperforms NAM under different setups in our experiments.

Analogously, when fixing $\delta$, the theorems imply that $K\asymp_\delta \kappa$ for SGD and $K\asymp_\delta\sqrt{\kappa}$ for the other four accelerated methods. These relationships are also discovered in Figure \ref{fig:2}. Thus, we conclude that the accelerated methods indeed achieve accelerated rates by improving the dependence on the condition number from $\kappa$ to $\sqrt{\kappa}$. In addition, Figure \ref{fig:2} suggests that different methods have~different robustness ability to the multiplicative noise. To be specific, as $\delta$ increases, SGD and DAM+ exhibit robust performance, while NAM, RMM, and iDAM+ gradually perform worse, meaning that these three methods are more sensitive to the multiplicative noise. This observation also aligns with the presented theorems (cf. \cref{sec:4.4}). Overall, the numerical results support our sharp analysis and demonstrate that the established error bounds of each method are not merely upper bounds, but provide rather accurate descriptions of the true convergence behavior.

\section{Conclusion}\label{sec:6}

Momentum-based accelerated SGD has been successfully applied for training neural networks. However, the theoretical understanding of the acceleration in stochastic setting is less well understood. Under OUBV assumption, different accelerated methods have been proven to enjoy a faster rate than SGD. However, experiments have also shown that some accelerated methods cannot outperform SGD in a variety of problem instances when OUBV assumption fails. Thus, there exists an evident gap between the theory and empirical observations.

Our paper serves as a step towards understanding the behavior of different accelerated methods under a weaker and more realistic condition on the stochastic oracle---the growth condition. The growth condition assumes that the variance of stochastic noise is bounded by two parts: the additive part, which is bounded just like under OUBV assumption, and the multiplicative part, which is proportional to the square of the magnitude of the gradient. Under this weaker condition, we provide a comprehensive investigation of a variety of accelerated methods, including NAM, RMM, DAM+, and iDAM+. We show that all these methods converge under the growth condition with rates that are the same as their rates under OUBV assumption. However, their ability to tolerate noise is different. Among the methods that we studied, DAM+ achieves the accelerated rate for any $\delta\geq 0$, while the other three methods can only tolerate a small $\delta$ (e.g., $\delta < 1$). We do not provide a lower bound analysis for NAM, RMM, iDAM+, and, therefore, cannot conclude that these methods do not improve the convergence over SGD when $\delta\geq 1$. However, our results~indeed highly coincide with recent empirical observations that illustrate the lack of benefits of these methods over the simple SGD on LSR problems.

One of the limitations of our work is the sharpness of the analysis. Due to the lack of lower~bound analysis on $\delta$, it is unclear whether the condition of having a small $\delta$ for different methods is~necessary for acceleration; and how sharp such a condition is. Considering our results are aligned with experiments in \cite{Cyrus2018Robust}, we believe that a small $\delta$ is necessary for NAM and RMM to accelerate SGD. However, a rigorous lower bound analysis is required and is left for future work. In addition, we plan to extend the growth condition analysis to non-strongly convex and non-convex objectives.

\begin{APPENDICES}

\section{Proofs of Section \ref{sec:3}}\label{sup:pf:sec:3}

\subsection{Proof of Theorem \ref{thm:rmm}}\label{pf:thm:rmm}
The proof follows the dissipativity framework introduced in~\cref{sec:2}. We require the following preparation lemma, which generalizes \cite[Lemma 3]{Cyrus2018Robust}.

\begin{lemma}\label{lem:main-IQC}
Suppose $f \in \Fcal( \mcal,\Lcal)$. Let $\theta \in (0, 1]$, $\{x_i\}_{i=1}^k$ be any sequence, and $h_\theta(x)$ is defined in \eqref{equ:h}. Then, for all $k \geq 0$, we have
\begin{enumerate}[label=(\alph*),topsep=2pt]
\setlength\itemsep{0.2em}
\item $h_\theta(x_k) \geq 0$;
\item $(\nabla f_k - \theta \mcal (x_k - x_\star))^\top (\Lcal (x_k - x_\star) - \nabla f_k) \geq h_\theta(x_k)$;
\item $(\nabla f_k - \theta \mcal(x_k - x_\star))^\top \cbr{\Lcal (x_{k} - x_{k-1}) - (\nabla f_k - \nabla f_{k-1} ) }\geq h_\theta(x_k) - h_\theta(x_{k-1})$.
\end{enumerate}
\end{lemma}

\begin{proof}
Define $h_\theta^\prime(x) = f(x) - f(x_\star) -\theta \mcal \|x-x_\star\|^2/2$. Since $f \in \Fcal(\mcal, \Lcal)$ and $0<\theta\leq 1$, $h_\theta^\prime$ is convex and $(\Lcal-\theta\mcal)$-smooth. Since $h_\theta^\prime(x_\star) = 0$ and $\nabla h_\theta^\prime(x_\star) = 0$, Lemma \ref{lemma:smooth}(b) applied on $h_\theta^\prime$ gives us
\begin{equation*}
h_\theta^\prime(x) \geq h_\theta^\prime(x_\star) + \nabla h_\theta^\prime(x_\star)^\top (x - x_\star) + \frac{\|\nabla h_\theta^\prime(x) -  \nabla h_\theta^\prime(x_\star)\|^2}{2(\Lcal-\theta\mcal)}=  \frac{\|\nabla f(x) - \theta\mcal (x - x_\star)\|^2}{2(\Lcal-\theta\mcal)}.
\end{equation*}
Thus, (a) follows from rearranging terms in the inequality above. Since $\nabla h_\theta^\prime (x) = \nabla f(x) - \theta \mcal (x - x_\star)$, applying Lemma \ref{lemma:smooth}(b) gives us
\begin{align*}
	(\nabla f(x) - \theta \mcal (x-x_\star))^\top & (\Lcal (x - x_\star) - \nabla f(x)) = \nabla h_\theta^\prime (x)^\top((\Lcal-\theta\mcal)(x-x_\star) - \nabla h_\theta^\prime (x)) \\[2pt]
	&\geq (\Lcal-\theta\mcal)\cbr{h_\theta'(x) - h_\theta'(x_\star) + \frac{\|\nabla h_{\theta}'(x) - \nabla h_{\theta}'(x_\star)\|^2}{2(\Lcal - \theta\mcal)}} - \|\nabla h_\theta^\prime (x)\|^2 \\[2pt]
	&= (\Lcal-\theta\mcal)h_\theta'(x) - \frac{1}{2}\|\nabla h_\theta^\prime (x)\|^2 = h_\theta(x).
\end{align*}
This shows (b). Similarly, for any two points $x, x'$, applying Lemma~\ref{lemma:smooth}(b),
\begin{align*}
	(\nabla f(x) -& \theta \mcal (x-x_\star))^\top  \cbr{\Lcal (x - x') - (\nabla f(x) - \nabla f(x') ) } \nonumber\\[2pt]
	&= \nabla h_\theta^\prime (x)^\top \cbr{(\Lcal-\theta\mcal)(x - x') - (\nabla h_\theta^\prime (x) - \nabla h_\theta^\prime (x')) }  \\[2pt]
	&\geq (\Lcal-\theta\mcal)\cbr{h_\theta'(x) - h_\theta'(x') + \frac{\|\nabla h_\theta'(x) - \nabla h_\theta'(x')\|^2}{2(\Lcal - \theta\mcal) } } - \nabla h_\theta^\prime (x)^\top(\nabla h_\theta^\prime (x) - \nabla h_\theta^\prime (x'))  \\[2pt]
	&=  (\Lcal-\theta\mcal) h_\theta^\prime (x) - \frac{1}{2} \|\nabla h_\theta^\prime (x)\|^2  - \cbr{ (\Lcal-\theta\mcal) h_\theta^\prime (x') - \frac{1}{2} \|\nabla h_\theta^\prime (x')\|^2 }  \\[2pt]
	&= h_\theta(x) - h_\theta(x'),
\end{align*}
which proves (c). The proof is  now complete.
\end{proof}

We are now ready to prove Theorem~\ref{thm:rmm}. We define
\begin{equation*}
 \xi_k =
 \begin{bmatrix}
 y_k - x_\star\\
 z_k
 \end{bmatrix}, \quad 
 \omega_k =
 \begin{bmatrix}
 \nabla f_k \\
 \varepsilon_k
 \end{bmatrix},\quad  \Ab =
 \begin{bmatrix}
 1 & \beta \\
 0 & \beta
 \end{bmatrix}
 \otimes \Ib_d, \quad \Bb = -\frac{\eta}{\Lcal}
 \begin{bmatrix}
 1 & 1 \\
 1 & 1
 \end{bmatrix}
 \otimes \Ib_d.
\end{equation*}
Then $\xi_{k+1} = \Ab\xi_{k} + \Bb\omega_k$. For $0<\theta\leq 1$, we let $\tilde{\kappa} = \kappa/\theta$; and set $\eta$, $\beta$, $\lambda$, $\nu$ according to the theorem with $\rho$ satisfying $1 - 1/\sqrt{\tilde{\kappa}} \leq \rho \leq \sqrt{1-1/\tilde{\kappa}}$. We further define the column vectors
\begin{equation*}
U=(-\theta \mcal, - \theta \mcal\beta/\eta, 1,0)^\top, \quad V_1 = (\Lcal, \Lcal \beta/\eta, -1, 0)^\top, \quad V_2=(0, \Lcal(\eta+\beta-1)/\eta, -1, 0)^\top,
\end{equation*}
and define matrices
\begin{align*}
	\Xb_1 &=  -\left[ \frac{(1-\rho^2)}{2}(U V_1^\top + V_1 U^\top)+ \frac{\rho^2}{2} (U V_2^\top+ V_2 U^\top)\right] \otimes \Ib_d,\\[2pt]
	\Xb_2 &= \frac{\eta}{\Lcal}
	\begin{bmatrix}
		0 & 0 & 0 & \rho^2-1\\
		0 & 0 & 0 & -\beta\\
		0 & 0 & 0 & \eta/\Lcal\\
		\rho^2-1 & -\beta & \eta/\Lcal & \eta/\Lcal
	\end{bmatrix} \otimes \Ib_d, \quad\quad\quad \Pb =
	\begin{bmatrix}
		(1-\rho^2)^2 & \rho^2(1-\rho^2) \\
		\rho^2(1-\rho^2) & \rho^4
	\end{bmatrix} \otimes \Ib_d.
\end{align*}
By direct calculation, we have
\begin{equation} \label{eq:LME-NAM-20-2}
	\lambda
	\begin{bmatrix}
		\Ab^\top \Pb \Ab - \rho^2 \Pb & \Ab^\top \Pb \Bb \\
		\Bb^\top \Pb \Ab & \Bb^\top \Pb \Bb
	\end{bmatrix} - \Xb_1 - \lambda \Xb_2
	= - \nu U U^\top \otimes \Ib_d.
\end{equation}
Let $\tilde{x}_k=x_k -x_\star$. We have
\begin{subequations}
	\begin{align}
		&\hskip-0.4cm \nabla f_{k-1} \stackrel{\eqref{eq:NAM-20b}   }{=}\frac{\Lcal}{\eta}(\beta z_{k-1} - z_k)-\varepsilon_{k-1}, \label{pequ:4}\\[2pt]
		&\hskip-0.4cm x_{k} \stackrel{\eqref{eq:NAM-20a}}{=} y_k + \frac{\beta}{\eta} z_k \leftstackrel{\eqref{eq:NAM-20c}}{=} y_{k-1} + \frac{\eta+\beta}{\eta} z_{k}
		\stackrel{\eqref{eq:NAM-20b}}{=}  y_{k-1}+\frac{\eta \beta +\beta^2}{\eta} z_{k-1} - \frac{\eta+\beta}{\Lcal} (\nabla f_{k-1} + \varepsilon_{k-1}), \label{pequ:5}\\[2pt]
		&\hskip-0.4cm \tilde{x}_k - \tilde{x}_{k-1} = x_k - x_{k-1} \stackrel{\eqref{eq:NAM-20a}}{=} y_k + \frac{\beta}{\eta}z_k - y_{k-1} - \frac{\beta}{\eta}z_{k-1} \stackrel{\eqref{eq:NAM-20c}}{=}z_k + \frac{\beta}{\eta}(z_k - z_{k-1}). \label{pequ:6}
	\end{align}
\end{subequations}
Further, we obtain
\begin{align*}
	\Lcal (\tilde{x}_{k} - \tilde{x}_{k-1}) - (\nabla f_k - \nabla f_{k-1})
	\stackrel{\mathclap{\eqref{pequ:6}}}{=}& \;\; \Lcal z_k + \frac{\Lcal \beta }{\eta}(z_k-z_{k-1}) - \nabla f_k + \nabla f_{k-1}\\[2pt]
	\stackrel{\mathclap{\eqref{pequ:4}}}{=} &\;\; \Lcal z_k + \frac{\Lcal \beta }{\eta}(z_k-z_{k-1}) - \nabla f_k + \frac{\Lcal}{\eta}(\beta z_{k-1} - z_k)-\varepsilon_{k-1}  \\[2pt]
	= & \frac{\eta+\beta-1}{\eta} \Lcal z_k - \nabla f_k -\varepsilon_{k-1}.
\end{align*}
Combining the above expression with
Lemma \ref{lem:main-IQC}(c), 
for any $0< \theta\leq 1$,
we get
\begin{multline} \label{eq:IQC-stochastic}
	(\nabla f_k - \theta \mcal \tilde{x}_k)^\top \{\Lcal (\tilde{x}_{k} -  \tilde{x}_{k-1}) -  (\nabla f_k - \nabla f_{k-1} ) \} \\[2pt]
	= (\nabla f_k - \theta \mcal \tilde{x}_k)^\top \left( \frac{\eta+\beta-1}{\eta} \Lcal z_k - \nabla f_k - \varepsilon_{k-1}\right) \geq h_\theta(x_k) - h_\theta (x_{k-1}).
\end{multline}
In what follows, we aim to bound $\sum_{j=1}^J\lambda_j \EE [(\xi_k^\top, \omega_k^\top)\Xb_j(\xi_k^\top, \omega_k^\top)^\top]$ and further get the error~recursion via \eqref{eq:exp-dissipation-expectation}. First, we note that $\nu\geq 0$ and $\lambda\geq 0$ ensured by the condition on $\rho$. Second, we bound the first quadratic form $S_1$ in \eqref{eq:exp-dissipation-expectation}:
\begin{align}\label{pequ:7}
	\begin{bmatrix}
		\xi_k \\ \omega_k
	\end{bmatrix}^\top
	\Xb_1
	\begin{bmatrix}
		\xi_k \\ \omega_k
	\end{bmatrix} &=  -(1-\rho^2)(\xi_k^\top, \omega_k^\top) U\cdot (\xi_k^\top, \omega_k^\top)V_1 - \rho^2(\xi_k^\top, \omega_k^\top) U\cdot (\xi_k^\top, \omega_k^\top)V_2 \nonumber\\[2pt]
	&= -(1-\rho^2)\rbr{\nabla f_k - \theta\mcal (y_k - x_\star) - \frac{\theta\mcal\beta}{\eta}z_k}^\top \rbr{\Lcal(y_k - x_\star) + \frac{\Lcal\beta}{\eta}z_k - \nabla f_k} \nonumber\\[2pt]
	& \quad -\rho^2\rbr{\nabla f_k - \theta\mcal (y_k - x_\star) - \frac{\theta\mcal\beta}{\eta}z_k}^\top\rbr{\frac{\Lcal(\eta+\beta-1)}{\eta}z_k - \nabla f_k} \nonumber\\[2pt]
	& \stackrel{\mathclap{\eqref{eq:NAM-20a}}}{=} \;\; - (1-\rho^2)\rbr{\nabla f_k -\theta\mcal(x_k - x_\star) }^\top\rbr{\Lcal(x_k - x_\star) - \nabla f_k} \nonumber\\[2pt]
	& \ \quad -\rho^2\rbr{\nabla f_k -\theta\mcal(x_k - x_\star) }^\top\rbr{\frac{\Lcal(\eta+\beta-1)}{\eta}z_k - \nabla f_k} \nonumber\\[2pt]
	&\hskip-0.8cm\stackrel{\text{Lemma \ref{lem:main-IQC}(b)}}{\leq} -(1-\rho^2)h_\theta(x_k)  -\rho^2\rbr{\nabla f_k -\theta\mcal(x_k - x_\star) }^\top\rbr{\frac{\Lcal(\eta+\beta-1)}{\eta}z_k - \nabla f_k} \nonumber\\[2pt]
	&\stackrel{\mathclap{\eqref{eq:IQC-stochastic} } }{\leq} -(1-\rho^2)h_\theta(x_k) - \rho^2(h_\theta(x_k) - h_\theta(x_{k-1})) - \rho^2\rbr{\nabla f_k -\theta\mcal(x_k - x_\star) }^\top\varepsilon_{k-1} \nonumber\\[2pt]
	&=  - h_\theta(x_k) + \rho^2h_\theta(x_{k-1}) - \rho^2\rbr{\nabla f_k -\theta\mcal(x_k - x_\star) }^\top\varepsilon_{k-1}.
\end{align}
Recall that $\Gcal_k$ is the $\sigma$-algebra containing randomness $\{g_i\}_{i=0}^{k-1}$, which is generated by $\{z_j, y_j, x_j\}_{j=0}^k$ for \eqref{eq:NAM-20}. Thus, the last term may not mean zero. We deal with it as follows:
\begin{align*}
\EE [- & \rbr{\nabla f_k - \theta\mcal(x_k - x_\star) }^\top\varepsilon_{k-1} \mid \Gcal_{k-1} ] \stackrel{\eqref{pequ:1}}{=}  \EE[\theta\mcal x_k^\top\varepsilon_{k-1} - \nabla f_k^\top \varepsilon_{k-1}\mid \Gcal_{k-1}]\\[2pt]
& \stackrel{\eqref{pequ:5}}{=} \EE\bigg[\theta\mcal \rbr{\l_{k-1} - \frac{\eta+\beta}{\Lcal} \varepsilon_{k-1}}^\top\varepsilon_{k-1} + \cbr{\nabla f\rbr{l_{k-1}} - \nabla f\rbr{\l_{k-1} - \frac{\eta+\beta}{\Lcal} \varepsilon_{k-1}}  }^\top \varepsilon_{k-1}\\[2pt]
&\quad\quad\quad\quad  - \nabla f(l_{k-1})^\top\varepsilon_{k-1} \mid \Gcal_{k-1}\bigg],
\end{align*}
where $l_{k-1} = y_{k-1}+\frac{\eta \beta +\beta^2}{\eta} z_{k-1} - \frac{\eta+\beta}{\Lcal} \nabla f_{k-1}$ belongs to $\Gcal_{k-1}$. Applying Lemma \ref{lemma:smooth}(d), we have
\begin{equation*}
\cbr{\nabla f\rbr{l_{k-1}} - \nabla f\rbr{\l_{k-1} - \frac{\eta+\beta}{\Lcal} \varepsilon_{k-1}}  }^\top \varepsilon_{k-1} \leq (\beta+\eta)\|\varepsilon_{k-1}\|^2.
\end{equation*}
Combining the above two displays, we further obtain
\begin{align*}
	\EE [-  \rbr{\nabla f_k - \theta\mcal(x_k - x_\star) }^\top\varepsilon_{k-1} \mid \Gcal_{k-1} ] &\leq  \frac{(\Lcal - \theta\mcal)(\eta+\beta)}{\Lcal}\EE[\|\varepsilon_{k-1}\|^2\mid \Gcal_{k-1}]\\[2pt]
	& \stackrel{\mathclap{\eqref{eq:GC}}}{\leq}\;\; \frac{\delta(\Lcal - \theta\mcal)(\eta+\beta)}{\Lcal}\|\nabla f_{k-1}\|^2 + \frac{\sigma^2(\Lcal - \theta\mcal)(\eta+\beta)}{\Lcal}.
\end{align*}
Combining the above inequality with \eqref{pequ:7} and taking the full expectation, we have
\begin{multline}\label{pequ:8}
	\EE\sbr{\begin{bmatrix}
			\xi_k \\ \omega_k
		\end{bmatrix}^\top
		\Xb_1
		\begin{bmatrix}
			\xi_k \\ \omega_k
	\end{bmatrix}}\\[2pt]
	\leq -\EE h_\theta(x_k) + \rho^2 \rbr{\EE h_{\theta}(x_{k-1}) + \frac{\delta(\Lcal - \theta\mcal)(\eta+\beta)}{\Lcal} \EE\|\nabla f_{k-1}\|^2} + \frac{\rho^2\sigma^2(\Lcal - \theta\mcal)(\eta+\beta)}{\Lcal}.
\end{multline}
Third, we bound the second quadratic form $S_2$ in \eqref{eq:exp-dissipation-expectation}:
\begin{multline}\label{pequ:9}
	\EE\sbr{\begin{bmatrix}
			\xi_k \\ \omega_k
		\end{bmatrix}^\top
		\lambda \Xb_2
		\begin{bmatrix}
			\xi_k \\ \omega_k
	\end{bmatrix} } = \frac{\lambda\eta^2}{\Lcal^2}\EE[\|\varepsilon_{k}\|^2] \\[2pt]
	\stackrel{\eqref{eq:GC}}{\leq}  \frac{\delta\lambda\eta^2}{\Lcal^2}\EE\|\nabla f_k\|^2 + \frac{\lambda\eta^2\sigma^2}{\Lcal^2} = \delta(1-\nu)\EE\|\nabla f_k\|^2  + (1-\nu)\sigma^2,
\end{multline}
where the last equality uses the fact that $\lambda\eta^2/\Lcal^2 = 1-\nu$. Fourth, we bound the third quadratic~form of $-\nu UU^\top\otimes \Ib_d$, which is done by
\begin{align}\label{pequ:10}
	& \EE\sbr{\begin{bmatrix}
			\xi_k \\ \omega_k
		\end{bmatrix}^\top
		(-\nu UU^\top\otimes \Ib_d)
		\begin{bmatrix}
			\xi_k \\ \omega_k
	\end{bmatrix} } =  -\nu\EE\sbr{\nbr{\nabla f_k - \theta\mcal \rbr{(y_k-x_\star) + \frac{\beta}{\eta}z_k}}^2} \nonumber\\[2pt]
	&\stackrel{\mathclap{\eqref{eq:NAM-20a}}}{=} \;\; -\nu\EE\|\nabla f_k - \theta\mcal (x_k - x_\star)\|^2
	=  -\nu \rbr{\EE\|\nabla f_k\|^2 -2\theta\mcal \EE\nabla f_k^\top(x_k-x_\star) + \theta^2\mcal^2\EE\|x_k - x_\star\|^2} \nonumber\\[2pt]
	&\stackrel{\text{Lemma } \ref{lemma:smooth}(a)}{\leq} -\nu\rbr{\EE\|\nabla f_k\|^2 -2\theta\mcal \EE\|\nabla f_k\|\cdot\|x_k-x_\star\| + \frac{\theta^2\mcal^2}{\Lcal^2}\EE\|\nabla f_k\|^2 } \nonumber\\[2pt]
	&\stackrel{\text{Lemma } \ref{lemma:convex}(c)}{\leq} -\nu \rbr{1-2\theta + \frac{\theta^2\mcal}{\Lcal^2}}\EE\|\nabla f_k\|^2.
\end{align}
Combining \eqref{pequ:8}, \eqref{pequ:9}, \eqref{pequ:10}, defining the following potential function
\begin{align*}
	V_k = \lambda \xi_k^\top\Pb\xi_k + h_\theta(x_{k-1}) + \frac{\delta(\Lcal - \theta\mcal)(\eta+\beta)}{\Lcal}\|\nabla f_{k-1}\|^2,
\end{align*}
we arrive the final conclusion
\begin{align*}
	\EE V_{k+1}
	&\leq \rho^2\EE V_k + \cbr{\delta(1-\nu) - \nu \rbr{1-2\theta + \frac{\theta^2\mcal^2}{\Lcal^2}} + \frac{\delta(\Lcal - \theta\mcal)(\eta+\beta)}{\Lcal} }\EE\|\nabla f_k\|^2 \\[2pt]
	&\quad + \rbr{\frac{\rho^2(\Lcal - \theta\mcal)(\eta+\beta)}{\Lcal} + 1-\nu}\sigma^2\\[2pt]
	&\leq  \rho^2\EE V_k - \cbr{\nu \rbr{1 - 2\theta + \frac{\theta^2\mcal^2}{\Lcal^2}} - \delta(1 - \nu + \eta + \beta)}\EE\|\nabla f_k\|^2 + \rbr{\rho^2(\eta+\beta) + 1-\nu}\sigma^2.
\end{align*}

\subsection{Proof of Corollary \ref{cor:rmm}}\label{pf:cor:rmm}

We consider the following two cases.

\noindent\textbf{Case 1: $\delta = 0$}. We apply \eqref{eq:RMM-rec} and have
\begin{equation}\label{npequ:5}
	\EE V_{k+1}\leq \rho^2\EE V_k - \nu \rbr{1 - 2\theta + \frac{1}{\tilde{\kappa}^2}} \EE\|\nabla f_k\|^2 + \rbr{\rho^2(\eta+\beta) + 1-\nu}\sigma^2.
\end{equation}
It hence suffices to have $\nu(1-2\theta + 1/\tilde{\kappa}^2) = 0$. 
We consider two subcases to discuss the fastest~convergence rate that we can obtain.
\vskip 5pt
\noindent\textbf{Case 1a: $1-2\theta + 1/\tilde{\kappa}^2 = 0$}. Using the condition $\theta\in(0, 1]$, we obtain $\theta = 1/(1+\sqrt{1 - 1/\kappa^2}) <1$. Using the condition $1-1/\sqrt{\tilde{\kappa}}\leq \rho \leq \sqrt{1-1/\tilde{\kappa}}$, we know that the smallest $\rho$ we can obtain is
\begin{equation*}
\rho = 1 - 1/\sqrt{\tilde{\kappa}} = 1 -\sqrt{\theta}/\sqrt{\kappa}.
\end{equation*}

\noindent\textbf{Case 1b: $\nu = 0$}. 
By the setup of $\nu$ in Theorem \ref{thm:rmm}, $\nu = 0$ implies $\rho = 1 - 1/\sqrt{\tilde{\kappa}} = 1 - \sqrt{\theta}/\sqrt{\kappa}$. However, in this case, the only restriction on $\theta$ is that $\theta\in(0, 1]$.

Comparing the above two subcases, we should set $\theta = 1$ and $\rho = 1 - 1/\sqrt{\kappa}$ to achieve the fastest convergence rate. Moreover, by the setup of $\beta$ and $\eta$ in Theorem \ref{thm:rmm} and the condition $\rho \in [1-1/\sqrt{\tilde{\kappa}}, \sqrt{1-1/\tilde{\kappa}}]$ on $\rho$,
we know
\begin{equation}\label{npequ:2}
	\eta =  \frac{\kappa}{\theta}(1-\rho)(1-\rho^2) = \frac{\kappa}{\theta} (1-\rho^2)(1+\rho) \leq 1+\rho,\quad\quad
	\beta = \frac{\rho^2}{1 - \theta/\kappa} \rho \leq\rho.
\end{equation}
Thus, \eqref{npequ:5} leads to
\begin{equation}\label{npequ:6}
	\EE V_{k+1}\leq \rho^2\EE V_k  + \rbr{\rho^2(\eta+\beta) + 1-\nu}\sigma^2\leq \rho^2\EE V_k + (3\rho^2+1)\sigma^2, \quad\quad \forall k\geq 1.
\end{equation}

\noindent\textbf{Case 2: $\delta\in(0, 1/4)$.}
For any $c\in(0, 1)$, we let $\rho = 1 - c/\sqrt{2\kappa}$ be our target convergence rate. We aim to establish a relation between $c$ and $\delta$. Note from \eqref{eq:RMM-rec} that we require $\nu(1-2\theta + 1/\tilde{\kappa}^2) >0$ in this case, in order to tolerate the effect of the multiplicative noise. Let $c_1>0$ be a constant to be specified later based on $c$. To ensure that inequality $1 - 2\theta + \theta^2/\kappa^2 \geq c_1 $ has a solution of $\theta$ in $(0, 1]$, we require
\begin{equation}\label{npequ:7}
0<c_1 \leq \max_{\theta\in (0, 1]} 1-2\theta + \theta^2/\kappa^2<1.
\end{equation}
Furthermore, solving $1 - 2\theta + \theta^2/\kappa^2 \geq c_1$ within $(0, 1]$ leads to
\begin{equation}\label{pequ:cond:theta}
	\theta \leq \frac{1-c_1}{1+ \sqrt{1 - \frac{1-c_1}{\kappa^2}}} \Longleftarrow \theta = \frac{1-c_1}{2}.
\end{equation}
Under the setup on $\theta$ in \eqref{pequ:cond:theta}, we can set $\eta$, $\beta$, $\lambda$, $\nu$, $\Pb$ as in the Theorem \ref{thm:rmm}. To ensure that $\rho = 1 - c/\sqrt{2\kappa}$ satisfies the condition on $\rho$ in Theorem \ref{thm:rmm}, we need
\begin{align}\label{npequ:1}
	1 - \frac{\sqrt{\theta}}{\sqrt{\kappa}} \leq \rho \leq \sqrt{1 - \frac{\theta}{\kappa}}
	&\Longleftrightarrow 1 - \frac{\sqrt{\theta}}{\sqrt{\kappa}} \leq 1 - \frac{c}{\sqrt{2\kappa}} \leq \sqrt{1 - \frac{\theta}{\kappa}}\nonumber \\[2pt]
	&\leftstackrel{\eqref{pequ:cond:theta}}{\Longleftrightarrow} 1 - \frac{\sqrt{1-c_1}}{\sqrt{2\kappa}} \leq 1 - \frac{c}{\sqrt{2\kappa}} \leq \sqrt{1 - \frac{1-c_1}{2\kappa}}  \nonumber\\[2pt]
	&\Longleftrightarrow c_1\leq 1-c^2 \text{ and } c_1\geq c^2+1 - 2c\sqrt{2\kappa} \nonumber \\[2pt]
	&\Longleftarrow c^2 - 2\sqrt{2}c + 1\leq c_1\leq 1-c^2.
\end{align}
One can easily see that the last condition on $c_1$ is valid since $c^2 - 2\sqrt{2}c + 1\leq 1-c^2$ for $c\in(0, 1)$. Combining the conditions \eqref{npequ:7} and \eqref{npequ:1}, we only need to choose $c_1$ to satisfy
\begin{equation}\label{npequ:8}
\max\{c^2 - 2\sqrt{2}c +1 , 0\}< c_1 \leq 1-c^2.
\end{equation}
Moreover, we have
\begin{align*}
	\nu\rbr{1 - 2\theta + \frac{1}{\tilde{\kappa}^2}}\geq \delta(1-\nu + \eta + \beta)
	&\stackrel{\eqref{pequ:cond:theta}}{\Longleftarrow} c_1\nu \geq \delta(1-\nu + \eta + \beta) \Longleftarrow c_1\nu \geq \delta(1 + \eta +\beta)\\[2pt]
	&\stackrel{\eqref{npequ:2}}{\Longleftarrow} c_1\nu \geq 2\delta(1+\rho)
	\Longleftrightarrow 1 - \tilde{\kappa}(1-\rho)^2 \geq \frac{4\delta\rho}{c_1} \\[2pt]
	& \Longleftrightarrow \delta \leq \frac{c_1}{4\rho}\rbr{1 - \frac{\kappa}{\theta}\cdot \frac{c^2}{2\kappa}}
	\stackrel{\eqref{pequ:cond:theta}}{\Longleftarrow} \delta \leq \frac{c_1}{4}\rbr{1 - \frac{c^2}{1-c_1}}.
\end{align*}
Inspired by the last inequality, we can let $c_1 = 1 - c$, which satisfies \eqref{npequ:8}, and further simplify~the condition on $\delta$ as $\delta \leq (1-c)^2/4$. Given $\delta\in(0, 1/4)$, we can let $c = 1 - 2\sqrt{\delta}\in(0, 1)$ to make such a condition on $\delta$ hold. Thus, $\theta \stackrel{\eqref{pequ:cond:theta}}{=}(1-c_1)/2 = c/2 = 1/2-\sqrt{\delta}$ and $\rho = 1 - c/\sqrt{2\kappa} = 1 - 2\theta/\sqrt{2\kappa} = 1 - \sqrt{2}\theta/\sqrt{\kappa}$. This validates the setup of $\theta, \delta$ in Case 2. Finally, combining \eqref{npequ:2} and \eqref{eq:RMM-rec}, we immediately have
\begin{equation}\label{npequ:3}
	\EE V_{k+1} \leq \rho^2\EE V_k + (3\rho^2+1)\sigma^2, \quad\quad \forall k\geq 1,
\end{equation}
which is consistent with \eqref{npequ:6} in Case 1.

\vskip 5pt
\noindent\textbf{Iterate convergence}. 
Applying \eqref{npequ:3} (or \eqref{npequ:6}) recursively, we have
\begin{equation}\label{npequ:4}
	\EE\xi_{k+1}^\top\Pb\xi_{k+1}\leq \frac{1}{\lambda}\EE V_{k+1} \leq \frac{1}{\lambda}\cbr{\rho^{2k}\EE V_1 + \frac{3\rho^2+1}{1-\rho^2}\sigma^2},
\end{equation}
where the first inequality is due to Lemma \ref{lem:main-IQC}(a) and $\|\nabla f_k\|^2\geq 0$. By the definition of $\Pb$ in~Theorem \ref{thm:rmm}, we have
\begin{equation*}
	\xi_{k+1}^\top\Pb\xi_{k+1} = \|q_{k+1}\|^2, \quad\text{ where } q_{k+1}\coloneqq  (1-\rho^2)(y_{k+1}-x_\star) + \rho^2z_{k+1}
	\stackrel{\eqref{eq:NAM-20c}}{=}(y_{k+1}-x_\star) - \rho^2(y_k-x_\star),
\end{equation*}
which implies $y_{k+1}-x_\star = \rho^2(y_k-x_\star) + q_{k+1} = \rho^{2k}(y_1-x_\star) + \sum_{t=1}^{k}\rho^{2(k-t)}q_{t+1}$, $\forall k\geq 0$. Therefore,
\begin{align*}
	\|y_{k+1}-x_\star\|^2 
	& = \nbr{\rho^{2k}(y_1-x_\star) + \sum_{t=1}^{k}\rho^{2(k-t)}q_{t+1}}^2\leq 2\rho^{4k}\|y_1-x_\star\|^2 + 2\nbr{\sum_{t=1}^{k}\rho^{2(k-t)}q_{t+1}}^2\\[2pt]
	&\leq 2\rho^{4k}\|y_1-x_\star\|^2 + 2\sum_{t=1}^{k}\rho^{0.5(k-t)}\cdot \sum_{t=1}^k\rho^{0.5(k-t)}\rho^{1.5^2(k-t)^2}\|q_{t+1}\|^2\\[2pt]
	&\leq 2\rho^{4k}\|y_1-x_\star\|^2 + 2\sum_{t=1}^{k}\rho^{0.5(k-t)}\cdot \sum_{t=1}^k\rho^{(0.5 + 1.5^2)(k-t)}\|q_{t+1}\|^2\\[2pt]
	&\leq 2\rho^{4k}\|y_1-x_\star\|^2 + \frac{2}{1-\sqrt{\rho}}\sum_{t=1}^{k}\rho^{2.75(k-t)}\|q_{t+1}\|^2.
\end{align*}
We further obtain
\begin{align*}
	& \EE\|y_{k+1}-x_\star\|^2  \leq 2\rho^{4k}\EE\|y_1-x_\star\|^2 + \frac{2}{1-\sqrt{\rho}}\sum_{t=1}^{k}\rho^{2.75(k-t)}\EE\|q_{t+1}\|^2\\[2pt]
	&\leftstackrel{\mathclap{\eqref{npequ:4}}}{\leq} \;\; 2\rho^{4k}\EE\|y_1-x_\star\|^2 + \frac{2}{\lambda(1-\sqrt{\rho})}\sum_{t=1}^{k}\rho^{2.75(k-t)}\cbr{\rho^{2t}\EE V_1 + \frac{3\rho^2 + 1}{1-\rho^2}\sigma^2}\\[2pt]
	&\leq  2\rho^{4k}\EE\|y_1-x_\star\|^2 + \frac{2\rho^{2k}\EE V_1}{\lambda(1-\sqrt{\rho})}\sum_{t=1}^{k}\rho^{0.75(k-t)} + \frac{2(3\rho^2+1)\sigma^2}{\lambda(1-\sqrt{\rho})(1-\rho^2)}\sum_{t=1}^{k}\rho^{2.75(k-t)}\\[2pt]
	&\leq  2\rho^{4k}\EE\|y_1-x_\star\|^2 + \frac{2\rho^{2k}\EE V_1}{\lambda(1-\sqrt{\rho})(1-\rho^{0.75})}+ \frac{2(3\rho^2+1)\sigma^2}{\lambda(1-\sqrt{\rho})(1-\rho^2)(1-\rho^{2.75})}
	=  \Ocal \left( \rho^{2(k+1)} +\sigma^2 \right),
\end{align*}
which completes the proof.

\section{Proofs of Section \ref{sec:4}}\label{sup:pf:sec:4}

\subsection{Proof of Theorem \ref{thm:iADA-robust}}\label{pf:thm:iADA-robust}

We require the following preparation lemma.

\begin{lemma} \label{lemma:DAM}
Consider the scheme \eqref{eq:iADA} and let $L_k$ be defined in \eqref{eq:DA-ld2}. Then,
\begin{multline*}
A_k L_{k+1} - A_{k-1} L_{k} \geq a_k f(x_k) +  a_k \nabla f_k^\top (\nabla \phi_k^\star(z_{k+1})-x_k)
\\- a_k \varepsilon_{k}^\top(x_\star- \nabla \phi_k^\star(z_{k+1}) )  + \frac{m_k}{2}\| \nabla \phi_{k}^\star (z_{k+1}) - \nabla \phi_{k}^\star (z_{k})\|^2.
\end{multline*}
\end{lemma}

\begin{proof}
By Lemma \ref{lemma:lower-bd}(a), we have $\nabla h_{k-1}(\nabla \phi_{k-1}^\star(z_{k})) =0$. Thus,
\begin{align} \label{eq:hk-1}
h_{k-1}(\nabla \phi_k^\star(z_{k+1})) & - h_{k-1}(\nabla \phi_{k-1}^\star(z_{k})) \nonumber\\[2pt]
&\leftstackrel{\mathclap{\eqref{eq:Bregman-divergence}}}{=}\;\;\Delta_{h_{k-1}}(\nabla \phi_k^\star(z_{k+1}), \nabla \phi_{k-1}^\star(z_{k}))  +\nabla h_{k-1}(\nabla \phi_{k-1}^\star(z_{k}))^\top (\nabla \phi_k^\star(z_{k+1})-\nabla \phi_{k-1}^\star(z_{k})) \nonumber\\[2pt]
&\leftstackrel{(*)}{=}\;\;\Delta_{h_{k-1}}(\nabla \phi_k^\star(z_{k+1}),\nabla \phi_{k-1}^\star(z_{k}))
\geq \frac{m_{k-1}}{2} \|\nabla \phi_k^\star(z_{k+1})-\nabla \phi_{k-1}^\star(z_{k})\|^{2},
\end{align}
where $(*)$ is due to $\nabla h_{k-1}(\nabla \phi_{k-1}^\star(z_{k})) = 0$ and the last inequality follows by the fact that $h_{k-1}$ is $m_{k-1}$-strongly convex with $m_{k-1} = A_{k-1} \mcal$. Furthermore, by \eqref{eq:closed-form},
\begin{equation} \label{eq:phik}
\nabla \phi_{k}^\star (z_{k}) = \frac{m_{k-1}}{m_{k}} \nabla \phi_{k-1}^\star (z_{k}) + \frac{\mcal a_k}{m_{k}} x_k.
\end{equation}
Thus, using the convexity of $\|\cdot\|^2$, we have
\begin{align} \label{eq:ji}
\left\| \nabla \phi_k^\star(z_{k+1}) - \nabla \phi_k^\star(z_{k})  \right\|^{2} \;\;
&\leftstackrel{\mathclap{\eqref{eq:phik}}}{=}\;\; \left\| \left( \frac{m_{k-1}}{m_k} + \frac{a_k \mcal}{m_k}\right)\nabla \phi_k^\star(z_{k+1}) - \frac{m_{k-1}}{m_k} \nabla \phi_{k-1}^\star(z_{k}) - \frac{a_k \mcal}{m_k} x_k \right\|^{2} \nonumber\\[3pt]
& \leq \frac{m_{k-1}}{m_k} \|\nabla \phi_k^\star(z_{k+1}) - \nabla \phi_{k-1}^\star(z_{k})\|^{2} + \frac{a_k \mcal}{m_k} \|\nabla \phi_k^\star(z_{k+1})-x_k\|^2,
\end{align}
where the first inequality uses the fact that $m_k = m_{k-1} + a_k \mcal$. Using the fact that
\begin{equation}\label{eq:hk2}
h_k(x) \stackrel{\eqref{eq:hk}}{=}h_{k-1}(x) + a_k g_k^\top (x-x_k) +  \frac{a_k \mcal}{2} \|x-x_k\|^2
\end{equation}
and combining all inequalities above, we obtain
\begin{align} \label{eq:hk-hk-1}
& h_{k}(\nabla \phi_k^\star(z_{k+1}))-h_{k-1}(\nabla \phi_{k-1}^\star(z_{k})) \nonumber\\
&\leftstackrel{\eqref{eq:hk2}}{=} h_{k-1}(\nabla \phi_k^\star(z_{k+1}))-h_{k-1}(\nabla \phi_{k-1}^\star(z_{k})) +a_k g_k^\top (\nabla \phi_k^\star(z_{k+1}) -x_k) +\frac{a_k \mcal}{2} \|\nabla \phi_k^\star(z_{k+1}) -x_k\|^2 \nonumber\\
&\leftstackrel{\eqref{eq:hk-1}}{\geq}\frac{m_{k-1}}{2} \|\nabla \phi_k^\star(z_{k+1})-\nabla \phi_{k-1}^\star(z_{k})\|^{2} + \frac{a_k \mcal}{2} \|\nabla \phi_k^\star(z_{k+1})-x_k\|^2 + a_k g_k^\top (\nabla \phi_k^\star(z_{k+1})-x_k) \nonumber\\
& \leftstackrel{\eqref{eq:ji}}{\geq} \frac{m_k}{2}\| \nabla \phi_{k}^\star (z_{k+1}) - \nabla \phi_{k}^\star (z_{k})\|^2 + a_k g_k^\top (\nabla \phi_k^\star(z_{k+1})-x_k).
\end{align}
As a result, applying Lemma \ref{lemma:lower-bd}(b), (c), and having
\begin{align*}
& A_k L_{k+1} - A_{k-1} L_{k} \stackrel{\eqref{eq:DA-ld2}}{=}  a_k f(x_k) - a_k \varepsilon_k^\top (x_\star-x_k)  +  h_k(\nabla \phi_k^\star(z_{k+1})) -  h_{k-1}(\nabla \phi_{k-1}^\star(z_{k})) \\[2pt]
&\leftstackrel{\eqref{eq:hk-hk-1}}{\geq}a_k f(x_k) - a_k \varepsilon_k^\top (x_\star-x_k) + a_k g_k^\top (\nabla \phi_k^\star(z_{k+1})-x_k) + \frac{m_k}{2}\| \nabla \phi_{k}^\star (z_{k+1}) - \nabla \phi_{k}^\star (z_{k})\|^2 \\[2pt]
& \;\;= a_k f(x_k) +  a_k \nabla f_k^\top (\nabla \phi_k^\star(z_{k+1})-x_k) - a_k \varepsilon_{k}^\top(x_\star- \nabla \phi_k^\star(z_{k+1}) ) + \frac{m_k}{2}\| \nabla \phi_{k}^\star (z_{k+1}) - \nabla \phi_{k}^\star (z_{k})\|^2,
\end{align*}
we complete the proof.
\end{proof}

We are now ready to prove Theorem \ref{thm:iADA-robust}. The upper bound of $A_k U_{k+1} - A_{k-1} U_{k}$~comes~from~Lemma \ref{lemma:ADAM}(a), since iDAM+ and DAM+ have the same $x_k$, $z_{k+1}$ update. The lower bound of $ A_k L_{k+1} - A_{k-1} L_{k}$ comes from Lemma~\ref{lemma:DAM}. Thus, we have
\begin{align} \label{eq:iADA-pot}
A_k V_{k+1} & - A_{k-1} V_{k} \nonumber\\[2pt]
& =  (A_k U_{k+1} - A_{k-1} U_{k}) - (A_k L_{k+1} - A_{k-1} L_{k}) \nonumber\\[2pt]
& \leq A_{k}(f(y_{k+1})-f(x_{k}))+a_{k}\nabla f_k^{\top}(\nabla \phi_{k}^{\star}(z_{k})-x_{k}) -  a_k \nabla f_k^\top (\nabla \phi_k^\star(z_{k+1})-x_k) \nonumber\\[2pt]
& \quad + a_k \varepsilon_{k}^\top(x_\star- \nabla \phi_k^\star(z_{k+1}) ) - \frac{m_k}{2}\| \nabla \phi_{k}^\star (z_{k+1}) - \nabla \phi_{k}^\star (z_{k})\|^2 \nonumber\\[2pt]
&= A_{k}(f(y_{k+1})-f(x_{k})) - a_{k}\nabla f_k^{\top}(\nabla \phi_{k}^{\star}(z_{k+1})-\nabla \phi_{k}^{\star}(z_{k})) + a_k \varepsilon_{k}^\top(x_\star-\nabla \phi_k^\star(z_{k+1})) \nonumber\\[2pt]
&\quad - \frac{m_k}{2}\| \nabla \phi_{k}^\star (z_{k+1}) - \nabla \phi_{k}^\star (z_{k})\|^2  \nonumber\\[2pt]
& \leftstackrel{\mathclap{\eqref{eq:iADAc}}}{=} \;\;A_{k}\rbr{ f(y_{k+1})-f(x_{k}) - \nabla f_k^{\top}(y_{k+1}-x_k)} + a_k \varepsilon_{k}^\top(x_\star-\nabla \phi_k^\star(z_{k+1})) \nonumber\\[2pt]
& \quad - \frac{m_k}{2}\| \nabla\phi_{k}^\star (z_{k+1}) - \nabla\phi_{k}^\star (z_{k})\|^2 \nonumber \\[2pt]
&  \leftstackrel{\mathclap{\eqref{equ:class}}}{\leq}  \;\;\frac{A_k \Lcal}{2} \|y_{k+1}-x_k\|^2 - \frac{m_k}{2}\| \nabla\phi_{k}^\star (z_{k+1}) - \nabla\phi_{k}^\star (z_{k})\|^2 + a_k \varepsilon_k^\top(x_\star-\nabla \phi_k^\star(z_{k+1})) \nonumber\\[2pt]
& \leftstackrel{\mathclap{\eqref{eq:iADAc}}}{=} \;\; \left( \frac{a_k^2 \Lcal}{2A_k} - \frac{m_k}{2} \right)\| \nabla\phi_{k}^\star (z_{k+1}) - \nabla\phi_{k}^\star (z_{k})\|^2  + a_k \varepsilon_{k}^\top(x_\star-\nabla \phi_k^\star(z_{k+1})) \nonumber\\[2pt]
& \leftstackrel{\mathclap{\eqref{eq:closed-form}}}{=} \;\; \left( \frac{a_k^2 \Lcal}{2 A_k} - \frac{m_k}{2} \right)\frac{a_k^2}{m_k^2} \| g_k\|^2 + a_k \varepsilon_{k}^\top(x_\star-\nabla \phi_k^\star(z_{k+1})),
\end{align}
where the second inequality follows by Lemma~\ref{lemma:ADAM}(a) and Lemma~\ref{lemma:DAM}. Recall that $\Gcal_k$ is the $\sigma$-algebra that contains all randomness $\{g_i\}_{i=1}^{k-1}$, which is generated by $\{x_j, z_j, y_j\}_{j=1}^{k}$ for \eqref{eq:iADA}. Thus, $\EE[\varepsilon_k^\top x_i]=0$ for all $i \leq k$, $\EE[\varepsilon_k^\top g_i]=0$ for $i<k$, and
\begin{equation*}
\EE \left[a_{k} \varepsilon_{k}^\top (x_\star-\nabla \phi_k^\star(z_{k+1})) | \Gcal_k \right]
\stackrel{\eqref{eq:closed-form}}{=}  a_k \EE \left\{\varepsilon_{k}^\top\left[ x_\star - \frac{-\sum_{i=1}^k a_i g_i + \mcal \sum_{i=1}^{k} a_{i} x_i }{m_k} \right] \middle| \Gcal_k \right\}
=  \frac{a_{k}^{2}}{m_{k}} \EE [\|\varepsilon_{k}\|^{2}| \Gcal_{k}].
\end{equation*}
Taking full expectation and multiplying by $1/A_k$ on both sides of \eqref{eq:iADA-pot}, we finish the proof by
\begin{align*}
 \EE \left[ V_{k+1} - \frac{A_{k-1}}{A_k} V_{k} \right]
& \leq\left( \frac{a_k^2 \Lcal}{ A_k^2} - \mcal \right)\frac{a_k^2}{2 A_k^2 \mcal^2} \EE \| g_k \|^2
+ \frac{a_{k}^{2}}{A_k^2 \mcal} \EE \|\varepsilon_{k}\|^{2} \\[2pt]
& \leftstackrel{\eqref{pequ:1}}{=} - \left( \mcal - \frac{a_k^2 \Lcal}{ A_k^2}  \right)\frac{a_k^2}{2 A_k^2 \mcal^2} \EE \|\nabla f_k \|^2
+ \left(  \mcal  + \frac{a_k^2 \Lcal}{ A_k^2}  \right) \frac{a_k^2}{2 A_k^2 \mcal^2}  \EE \|\varepsilon_{k}\|^{2} \\[2pt]
& \leftstackrel{\eqref{eq:GC}}{\leq} - \left( (1 - \delta )  \mcal  -(1+\delta)\frac{a_k^2 \Lcal}{A_k^2}\right) \frac{a_k^2}{2 A_k^2 \mcal^2} \EE \|\nabla f_{k}\|^{2} +\frac{a_k^2}{2 A_k^2 \mcal^2} \left(  \mcal  + \frac{a_k^2 \Lcal}{ A_k^2}  \right)   \sigma^2.
\end{align*}

\subsection{Proof of Corollary \ref{cor:iADA}}\label{pf:cor:iADA}

Combining $a_k/A_k = \sqrt{(1-\delta)/((1+\delta)\kappa)}$ with Theorem~\ref{thm:iADA-robust}~and noting that
\begin{equation*}
\rho_k =  \frac{A_{k-1}}{A_k} = 1 - \frac{a_k}{A_k} = 1 - \sqrt{\frac{1-\delta}{(1+\delta)\kappa}},\ \
\left(  \mcal  + \frac{a_k^2 \Lcal}{ A_k^2}  \right) \frac{a_k^2}{2 A_k^2 \mcal^2} = \rbr{1 + \frac{1-\delta}{1+\delta}}\frac{1-\delta}{2(1+\delta)\kappa\mcal} = \frac{1-\delta}{(1+\delta)^2\Lcal},
\end{equation*}
we immediately obtain the result.

\end{APPENDICES}

\bibliographystyle{informs2014}
\bibliography{ref}

\begin{thebibliography}{41}
\providecommand{\natexlab}[1]{#1}
\providecommand{\url}[1]{\texttt{#1}}
\providecommand{\urlprefix}{URL }

\bibitem[{Assran \protect\BIBand{} Rabbat(2020)}]{Assran2020Convergence}
Assran M, Rabbat M (2020) On the convergence of {N}esterov’s accelerated
  gradient method in stochastic settings. \emph{International Conference on
  Machine Learning},
  \urlprefix\url{http://proceedings.mlr.press/v119/assran20a.html}.

\bibitem[{Aybat et~al.(2019)Aybat, Fallah, Gurbuzbalaban, \protect\BIBand{}
  Ozdaglar}]{Aybat2019Universally}
Aybat NS, Fallah A, Gurbuzbalaban M, Ozdaglar A (2019) A universally optimal
  multistage accelerated stochastic gradient method. \emph{Advances in Neural
  Information Processing Systems},
  \urlprefix\url{https://proceedings.neurips.cc/paper/2019/file/d630553e32ae21fb1a6df39c702d2c5c-Paper.pdf}.

\bibitem[{Aybat et~al.(2020)Aybat, Fallah, Gürbüzbalaban, \protect\BIBand{}
  Ozdaglar}]{Aybat2020Robust}
Aybat NS, Fallah A, Gürbüzbalaban M, Ozdaglar A (2020) Robust accelerated
  gradient methods for smooth strongly convex functions. \emph{{SIAM} Journal
  on Optimization} 30(1):717--751,
  \urlprefix\url{http://dx.doi.org/10.1137/19m1244925}.

\bibitem[{Bauschke \protect\BIBand{} Borwein(1997)}]{Bauschke1997Legendre}
Bauschke HG, Borwein JM (1997) Legendre functions and the method of random
  bregman projections. \emph{Journal of Convex Analysis} 4(1):27--67,
  \urlprefix\url{http://eudml.org/doc/227096}.

\bibitem[{Bottou et~al.(2018)Bottou, Curtis, \protect\BIBand{}
  Nocedal}]{Bottou2018Optimization}
Bottou L, Curtis FE, Nocedal J (2018) Optimization methods for large-scale
  machine learning. \emph{{SIAM} Review} 60(2):223--311,
  \urlprefix\url{http://dx.doi.org/10.1137/16m1080173}.

\bibitem[{Bubeck(2015)}]{Bubeck2015Convex}
Bubeck S (2015) Convex optimization: Algorithms and complexity.
  \emph{Foundations and Trends{\textregistered} in Machine Learning}
  8(3-4):231--357, \urlprefix\url{http://dx.doi.org/10.1561/2200000050}.

\bibitem[{Cevher \protect\BIBand{} V{\~{u}}(2018)}]{Cevher2018linear}
Cevher V, V{\~{u}} BC (2018) On the linear convergence of the stochastic
  gradient method with constant step-size. \emph{Optimization Letters}
  13(5):1177--1187,
  \urlprefix\url{http://dx.doi.org/10.1007/s11590-018-1331-1}.

\bibitem[{Chen et~al.(2012)Chen, Lin, \protect\BIBand{} Pena}]{Chen2012Optimal}
Chen X, Lin Q, Pena J (2012) Optimal regularized dual averaging methods for
  stochastic optimization. \emph{Advances in neural information processing
  systems} 25,
  \urlprefix\url{https://proceedings.neurips.cc/paper/2012/hash/274ad4786c3abca69fa097b85867d9a4-Abstract.html}.

\bibitem[{Cohen et~al.(2018)Cohen, Diakonikolas, \protect\BIBand{}
  Orecchia}]{Cohen2018Acceleration}
Cohen M, Diakonikolas J, Orecchia L (2018) On acceleration with noise-corrupted
  gradients. \emph{International Conference on Machine Learning},
  \urlprefix\url{http://proceedings.mlr.press/v80/cohen18a.html}.

\bibitem[{Cyrus et~al.(2018)Cyrus, Hu, Scoy, \protect\BIBand{}
  Lessard}]{Cyrus2018Robust}
Cyrus S, Hu B, Scoy BV, Lessard L (2018) A robust accelerated optimization
  algorithm for strongly convex functions. \emph{2018 Annual American Control
  Conference ({ACC})}, 1376--1381, IEEE ({IEEE}),
  \urlprefix\url{http://dx.doi.org/10.23919/acc.2018.8430824}.

\bibitem[{Diakonikolas \protect\BIBand{}
  Orecchia(2018)}]{Diakonikolas2018Accelerated}
Diakonikolas J, Orecchia L (2018) Accelerated extra-gradient descent: A novel
  accelerated first-order method. \emph{9th Innovations in Theoretical Computer
  Science Conference (ITCS 2018)}, Schloss Dagstuhl-Leibniz-Zentrum fuer
  Informatik (Schloss Dagstuhl - Leibniz-Zentrum fuer Informatik GmbH,
  Wadern/Saarbruecken, Germany),
  \urlprefix\url{http://dx.doi.org/10.4230/LIPICS.ITCS.2018.23}.

\bibitem[{Diakonikolas \protect\BIBand{}
  Orecchia(2019)}]{Diakonikolas2019Approximate}
Diakonikolas J, Orecchia L (2019) The approximate duality gap technique: A
  unified theory of first-order methods. \emph{{SIAM} Journal on Optimization}
  29(1):660--689, \urlprefix\url{http://dx.doi.org/10.1137/18m1172314}.

\bibitem[{Ghadimi \protect\BIBand{} Lan(2013)}]{Ghadimi2013Optimal}
Ghadimi S, Lan G (2013) Optimal stochastic approximation algorithms for
  strongly convex stochastic composite optimization, {II}: Shrinking procedures
  and optimal algorithms. \emph{{SIAM} Journal on Optimization}
  23(4):2061--2089, \urlprefix\url{http://dx.doi.org/10.1137/110848876}.

\bibitem[{Gower et~al.(2019)Gower, Loizou, Qian, Sailanbayev, Shulgin,
  \protect\BIBand{} Richt{\'a}rik}]{Gower2019SGD}
Gower RM, Loizou N, Qian X, Sailanbayev A, Shulgin E, Richt{\'a}rik P (2019)
  Sgd: General analysis and improved rates. \emph{International Conference on
  Machine Learning},
  \urlprefix\url{http://proceedings.mlr.press/v97/qian19b.html}.

\bibitem[{Grimmer(2019)}]{Grimmer2019Convergence}
Grimmer B (2019) Convergence rates for deterministic and stochastic subgradient
  methods without lipschitz continuity. \emph{{SIAM} Journal on Optimization}
  29(2):1350--1365, \urlprefix\url{http://dx.doi.org/10.1137/18m117306x}.

\bibitem[{Gürbüzbalaban et~al.(2015)Gürbüzbalaban, Ozdaglar,
  \protect\BIBand{} Parrilo}]{Guerbuezbalaban2015globally}
Gürbüzbalaban M, Ozdaglar A, Parrilo P (2015) A globally convergent
  incremental newton method. \emph{Mathematical Programming} 151(1):283--313,
  \urlprefix\url{http://dx.doi.org/10.1007/s10107-015-0897-y}.

\bibitem[{Hu \protect\BIBand{} Lessard(2017)}]{Hu2017Dissipativity}
Hu B, Lessard L (2017) Dissipativity theory for nesterov's accelerated method.
  \emph{International Conference on Machine Learning},
  \urlprefix\url{http://proceedings.mlr.press/v70/hu17a.html}.

\bibitem[{Hu et~al.(2020)Hu, Seiler, \protect\BIBand{}
  Lessard}]{Hu2020Analysis}
Hu B, Seiler P, Lessard L (2020) Analysis of biased stochastic gradient descent
  using sequential semidefinite programs. \emph{Mathematical Programming}
  187(1-2):383--408,
  \urlprefix\url{http://dx.doi.org/10.1007/s10107-020-01486-1}.

\bibitem[{Hu et~al.(2018)Hu, Wright, \protect\BIBand{}
  Lessard}]{Hu2018Dissipativity}
Hu B, Wright S, Lessard L (2018) Dissipativity theory for accelerating
  stochastic variance reduction: A unified analysis of svrg and katyusha using
  semidefinite programs. \emph{Proceedings of Machine Learning Research} 80,
  \urlprefix\url{http://proceedings.mlr.press/v80/hu18b.html}.

\bibitem[{Jain et~al.(2018{\natexlab{a}})Jain, Kakade, Kidambi, Netrapalli,
  \protect\BIBand{} Sidford}]{Jain2018Accelerating}
Jain P, Kakade SM, Kidambi R, Netrapalli P, Sidford A (2018{\natexlab{a}})
  Accelerating stochastic gradient descent for least squares regression.
  \emph{Conference On Learning Theory},
  \urlprefix\url{http://proceedings.mlr.press/v75/jain18a.html}.

\bibitem[{Jain et~al.(2018{\natexlab{b}})Jain, Kakade, Kidambi, Netrapalli,
  \protect\BIBand{} Sidford}]{Jain2018Parallelizing}
Jain P, Kakade SM, Kidambi R, Netrapalli P, Sidford A (2018{\natexlab{b}})
  Parallelizing stochastic gradient descent for least squares regression:
  Mini-batching, averaging, and model misspecification. \emph{Journal of
  Machine Learning Research} 18(223):1--42,
  \urlprefix\url{http://jmlr.org/papers/v18/16-595.html}.

\bibitem[{Jofr{\'{e}} \protect\BIBand{} Thompson(2018)}]{Jofre2018variance}
Jofr{\'{e}} A, Thompson P (2018) On variance reduction for stochastic smooth
  convex optimization with multiplicative noise. \emph{Mathematical
  Programming} 174(1-2):253--292,
  \urlprefix\url{http://dx.doi.org/10.1007/s10107-018-1297-x}.

\bibitem[{Kidambi et~al.(2018)Kidambi, Netrapalli, Jain, \protect\BIBand{}
  Kakade}]{Kidambi2018Insufficiency}
Kidambi R, Netrapalli P, Jain P, Kakade S (2018) On the insufficiency of
  existing momentum schemes for stochastic optimization. \emph{Information
  Theory and Applications Workshop ({ITA})}, 1--9, IEEE ({IEEE}),
  \urlprefix\url{http://dx.doi.org/10.1109/ita.2018.8503173}.

\bibitem[{Lessard et~al.(2016)Lessard, Recht, \protect\BIBand{}
  Packard}]{Lessard2016Analysis}
Lessard L, Recht B, Packard A (2016) Analysis and design of optimization
  algorithms via integral quadratic constraints. \emph{{SIAM} Journal on
  Optimization} 26(1):57--95,
  \urlprefix\url{http://dx.doi.org/10.1137/15m1009597}.

\bibitem[{Liu \protect\BIBand{} Belkin(2020)}]{Liu2020Accelerating}
Liu C, Belkin M (2020) Accelerating sgd with momentum for over-parameterized
  learning. \emph{International Conference on Learning Representations},
  \urlprefix\url{https://openreview.net/forum?id=r1gixp4FPH}.

\bibitem[{Mourtada(2022)}]{Mourtada2022Exact}
Mourtada J (2022) Exact minimax risk for linear least squares, and the lower
  tail of sample covariance matrices. \emph{The Annals of Statistics}
  50(4):2157--2178, \urlprefix\url{http://dx.doi.org/10.1214/22-aos2181}.

\bibitem[{Nemirovskij \protect\BIBand{} Yudin(1983)}]{Nemirovskij1983Problem}
Nemirovskij AS, Yudin DB (1983) \emph{Problem complexity and method efficiency
  in optimization} (Wiley-Interscience),
  \urlprefix\url{https://www2.isye.gatech.edu/~nemirovs/Nemirovskii_Yudin_1983.pdf}.

\bibitem[{Nesterov(2004)}]{Nesterov2004Introductory}
Nesterov Y (2004) \emph{Introductory Lectures on Convex Optimization},
  volume~87 (Springer {US}),
  \urlprefix\url{http://dx.doi.org/10.1007/978-1-4419-8853-9}.

\bibitem[{Nesterov(1983)}]{Nesterov1983method}
Nesterov YE (1983) A method for solving the convex programming problem with
  convergence rate $o(1/k^2)$. \emph{Dokl. akad. nauk Sssr}, volume 269,
  543--547, \urlprefix\url{http://www.ams.org/mathscinet-getitem?mr=0701288}.

\bibitem[{Nguyen et~al.(2019)Nguyen, Nguyen, Richt{\'{a}}rik, Scheinberg,
  Tak{\'{a}}c, \protect\BIBand{} van Dijk}]{Nguyen2019New}
Nguyen LM, Nguyen PH, Richt{\'{a}}rik P, Scheinberg K, Tak{\'{a}}c M, van Dijk
  M (2019) New convergence aspects of stochastic gradient algorithms.
  \emph{Journal of Machine Learning Research} 20:176:1--176:49,
  \urlprefix\url{http://jmlr.org/papers/v20/18-759.html}.

\bibitem[{Polyak(1964)}]{Polyak1964Some}
Polyak B (1964) Some methods of speeding up the convergence of iteration
  methods. \emph{{USSR} Computational Mathematics and Mathematical Physics}
  4(5):1--17, \urlprefix\url{http://dx.doi.org/10.1016/0041-5553(64)90137-5}.

\bibitem[{Rakhlin et~al.(2012)Rakhlin, Shamir, \protect\BIBand{}
  Sridharan}]{Rakhlin2012Making}
Rakhlin A, Shamir O, Sridharan K (2012) Making gradient descent optimal for
  strongly convex stochastic optimization. \emph{International Conference on
  Machine Learning}, \urlprefix\url{http://icml.cc/2012/papers/261.pdf}.

\bibitem[{Schmidt \protect\BIBand{} Roux(2013)}]{Schmidt2013Fast}
Schmidt M, Roux NL (2013) Fast convergence of stochastic gradient descent under
  a strong growth condition. \emph{arXiv preprint arXiv:1308.6370}
  \urlprefix\url{https://arxiv.org/abs/1308.6370}.

\bibitem[{Scoy et~al.(2018)Scoy, Freeman, \protect\BIBand{}
  Lynch}]{Scoy2018Fastest}
Scoy BV, Freeman RA, Lynch KM (2018) The fastest known globally convergent
  first-order method for minimizing strongly convex functions. \emph{{IEEE}
  Control Systems Letters} 2(1):49--54,
  \urlprefix\url{http://dx.doi.org/10.1109/lcsys.2017.2722406}.

\bibitem[{Solodov(1998)}]{Solodov1998Incremental}
Solodov M (1998) Incremental gradient algorithms with stepsizes bounded away
  from zero. \emph{Computational Optimization and Applications} 11(1):23--35,
  \urlprefix\url{http://dx.doi.org/10.1023/a:1018366000512}.

\bibitem[{Stich(2019)}]{Stich2019Unified}
Stich SU (2019) Unified optimal analysis of the (stochastic) gradient method.
  \emph{arXiv preprint arXiv:1907.04232}
  \urlprefix\url{https://arxiv.org/abs/1907.04232}.

\bibitem[{Tseng(1998)}]{Tseng1998Incremental}
Tseng P (1998) An incremental gradient(-projection) method with momentum term
  and adaptive stepsize rule. \emph{{SIAM} Journal on Optimization}
  8(2):506--531, \urlprefix\url{http://dx.doi.org/10.1137/s1052623495294797}.

\bibitem[{Vaswani et~al.(2019)Vaswani, Bach, \protect\BIBand{}
  Schmidt}]{Vaswani2019Fast}
Vaswani S, Bach FR, Schmidt M (2019) Fast and faster convergence of {SGD} for
  over-parameterized models and an accelerated perceptron. \emph{International
  Conference on Artificial Intelligence and Statistics},
  \urlprefix\url{http://proceedings.mlr.press/v89/vaswani19a.html}.

\bibitem[{Xiao(2010)}]{Xiao2010Dual}
Xiao L (2010) Dual averaging methods for regularized stochastic learning and
  online optimization. \emph{Journal of Machine Learning Research}
  11(Oct):2543--2596, \urlprefix\url{http://jmlr.org/papers/v11/xiao10a.html}.

\bibitem[{Xu et~al.(2018)Xu, Wang, \protect\BIBand{} Gu}]{Xu2018Continuous}
Xu P, Wang T, Gu Q (2018) Continuous and discrete-time accelerated stochastic
  mirror descent for strongly convex functions. \emph{International Conference
  on Machine Learning}
  \urlprefix\url{http://proceedings.mlr.press/v80/xu18g.html}.

\bibitem[{Zhou(2018)}]{Zhou2018fenchel}
Zhou X (2018) On the fenchel duality between strong convexity and lipschitz
  continuous gradient. \emph{arXiv preprint arXiv:1803.06573}
  \urlprefix\url{https://arxiv.org/abs/1803.06573}.

\end{thebibliography}

\end{document}